\documentclass[aop,seceqn,citesort,MSNbibl,dvips]{arximspdf}
\usepackage{mathbh}
\usepackage{graphicx}

\doi{10.1214/10-AOP637}
\volume{40}
\issue{2}
\pubyear{2012}
\firstpage{459}
\lastpage{534}

\makeatletter

\newcommand{\iitem}{\fontsize{8.36pt}{8.36pt}\selectfont{(i)}}
\newcommand{\ii}{\fontsize{8.36pt}{8.36pt}\selectfont{(ii)}}
\newcommand{\iii}{\fontsize{8.36pt}{8.36pt}\selectfont{(iii)}}

\newcommand{\Z}{\mathbb{Z}}
\newcommand{\N}{\mathbb{N}}
\newcommand{\R}{\mathbb{R}}

\newcommand{\E}{\mathbb{E}}
\newcommand{\PP}{\mathbb{P}}
\newcommand{\Var}{\operatorname{Var}}
\newcommand{\Cov}{\operatorname{Cov}}
\newcommand{\esssup}{\operatorname{ess}\sup}

\newtheorem{theorem}{Theorem}[section]
\newtheorem{lemma}[theorem]{Lemma}
\newtheorem{corollary}[theorem]{Corollary}
\newtheorem{proposition}[theorem]{Proposition}

\newproclaim{remark}[theorem]{Remark}

\newtheorem{conjecture}[theorem]{Conjecture}
\newtheorem{claim}[theorem]{Claim}

\newproclaim{definition}[theorem]{Definition}

\makeatother

\begin{document}
\begin{frontmatter}

\title{Quenched exit estimates and ballisticity conditions for
higher-dimensional random walk in~random environment}
\runtitle{Exit estimates and ballisticity conditions}

\begin{aug}
\author[A]{\fnms{Alexander} \snm{Drewitz}\corref{}\thanksref{t1,t2}\ead[label=e1]{alexander.drewitz@math.ethz.ch}}
and
\author[B]{\fnms{Alejandro F.} \snm{Ram\'{i}rez}\thanksref{t3}\ead[label=e2]{aramirez@mat.puc.cl}}
\runauthor{A. Drewitz and A. F. Ram\'{i}rez}
\affiliation{ETH Z\"{u}rich and Pontificia Unversidad Cat\'{o}lica de Chile}
\address[A]{Departement Mathematik\\
ETH Z\"{u}rich\\
R\"{a}mistrasse 101\\
8092 Z\"{u}rich\\
Switzerland\\
\printead{e1}}
\address[B]{Facultad de Matem\'{a}ticas\\
Pontificia Universidad Cat\'{o}lica de Chile\\
Vicu\~{n}a Mackenna 4860, Macul\\
Santiago 6904441\\
Chile\\
\printead{e2}}
\end{aug}

\thankstext{t1}{Supported in part by the International Research
Training Group
``Stochastic Models of Complex Processes.''}
\thankstext{t2}{Supported in part by the Berlin Mathematical School.}
\thankstext{t3}{Supported in part by Fondo Nacional de Desarrollo
Cient\'{\i}fico
y Tecnol\'{o}gico Grant 1100298.}

\received{\smonth{5} \syear{2010}}
\revised{\smonth{11} \syear{2010}}

%
\begin{abstract}
Consider a random walk in an i.i.d. uniformly elliptic environment in
dimensions larger than one. In 2002, Sznitman introduced for each
$\gamma\in(0,1)$ the ballisticity condition~$(T)_\gamma$ and the
condition~$(T')$ defined as the fulfillment of $(T)_\gamma$ for each
$\gamma\in(0,1)$. Sznitman proved that~$(T')$ implies a ballistic law
of large numbers. Furthermore, he showed that for all $\gamma\in
(0.5,1)$, $(T)_\gamma$ is equivalent to $(T')$. Recently, Berger has
proved that in dimensions larger than three, for each $\gamma\in
(0,1)$, condition $(T)_\gamma$ implies a ballistic law of large
numbers. On the other hand, Drewitz and Ram\'{\i}rez have shown that in
dimensions $d\ge2$ there is a constant $\gamma_d\in(0.366,0.388)$
such that for each $\gamma\in(\gamma_d,1)$, condition $(T)_\gamma$ is
equivalent to $(T')$. Here, for dimensions larger than three, we extend
the previous range of equivalence to all $\gamma\in(0,1)$. For the
proof, the so-called \textit{effective criterion} of Sznitman is
established employing a sharp estimate for the probability of atypical
quenched exit distributions of the walk leaving certain boxes. In this
context, we also obtain an affirmative answer to a conjecture raised by
Sznitman in 2004 concerning these probabilities. A~key ingredient for
our estimates is the multiscale method developed recently by
Berger.\looseness=-1
\end{abstract}

%

\end{frontmatter}

\section{Introduction and statement of the main results}
We continue our investigation of the interrelations between the
ballisticity conditions~$(T)_\gamma$ and~$(T')$
introduced by Sznitman in \cite{Sz-02} for random walk in random
environment (RWRE). In dimensions larger than or equal to four,
the
results we establish in\vadjust{\goodbreak} this paper amount to a considerable improvement
of what has been obtained in our work~\cite{DrRa-09b}.
To prove the corresponding results, we take advantage of techniques
recently developed
by Berger in \cite{Be-09}. We derive sharp estimates on the
probability of certain quenched exit
distributions of the RWRE and thereby provide an affirmative answer to
a slightly stronger version of a~conjecture announced by Sznitman in~\cite{Sz-04}.


We start by giving an introduction to the model, thereby fixing the
notation we employ.
Denote by $\mathcal{M}_d$ the space of probability measures on the set
$\{e \in\Z^d \dvtx\Vert e \Vert_1 = 1\}$
of canonical unit vectors
and set $\Omega:= (\mathcal{M}_d)^{\Z^d}$.
For each \textit{environment} $\omega=(\omega(x,\cdot))_{x \in\Z^d}
\in\Omega$,
we consider the Markov chain $(X_n)_{n\in\N}$ with transition
probabilities from $x$ to $x+e$
given by $\omega(x,e)$ for $\Vert e \Vert_1 = 1$, and~$0$ otherwise.
We denote by $P_{x,\omega}$ the law of this Markov chain conditioned
on $\{X_0 = x\}$.
Furthermore, let $\PP$ be a probability measure on $\Omega$
such that
the coordinates $(\omega(x,\cdot))_{x \in\Z^d}$ of the environment
$\omega$
are i.i.d. under~$\PP$. Then~$\PP$ is called \textit{elliptic}
if $\PP(\min_{\Vert e \Vert_1 =1} \omega(0,e)>0)=1$ while it is
called \textit{uniformly elliptic} if
there is a constant $\kappa>0$
such that
$
\PP(\min_{\Vert e \Vert_1 =1} \omega(0,e)\ge\kappa)=1.
$
We
call~$P_{x,\omega}$ the \textit{quenched law} of the RWRE starting from
$x$, and correspondingly
we define the \textit{averaged} (or \textit{annealed}) law of the RWRE by
$P_x:=\int_\Omega P_{x,\omega} \PP(d\omega)$.

Given a direction $l\in\mathbb S^d$, we say that the RWRE is \textit{transient
in the direction $l$} if
\[
P_{0}\Bigl(\lim_{n\to\infty} X_n\cdot l=\infty\Bigr)=1.
\]
Furthermore, we say that the RWRE is \textit{ballistic in the direction $l$}
if $P_0$-a.s.
\[
\liminf_{n\to\infty}\frac{X_n\cdot l}{n}>0.
\]
It is well known that in dimension one there exists
uniformly elliptic RWRE in i.i.d. environments which is
transient but not ballistic to the right.
It was also recently established that in dimensions larger than one
there exists elliptic RWRE in i.i.d. environments which is
transient but not ballistic in a given direction see Sabot and
Tournier in \cite{SaTo-09}.
Nevertheless, the following fundamental conjecture remains
open.\vspace*{-1pt}
%
\begin{conjecture} In dimensions larger than one, every uniformly elliptic
RWRE in an i.i.d. environment which is transient in a given direction is
necessarily ballistic in the same direction.\vspace*{-1pt}
\end{conjecture}

Some partial progress has been made toward the resolution of
this conjecture by studying transient RWRE satisfying some
additional assumptions introduced in \cite{Sz-02}, usually called
\textit{ballisticity
conditions}. For each $l\in\mathbb S^{d-1}$ and $L>0$, let us define
\[
T_L^l:=\inf\{n\ge0\dvtx X_n\cdot l>L\}.
\]
%

\begin{definition}
Let $\gamma\in(0,1)$ and $l\in\mathbb S^{d-1}$.
We say that \textit{condition $(T)_\gamma$} is satisfied with respect\vadjust{\goodbreak}
to $l$
[written $(T)_\gamma|l$ or $(T)_\gamma$] if for each $l'$ in a
neighborhood of $l$ and
each $b>0$ one has that
\[
\limsup_{L \to\infty} L^{-\gamma} \log P_0(T_L^{l'} > T_{bL}^{-l'})
< 0.
\]
We say that \textit{condition $(T')$} is satisfied with respect to $l$
[written~$(T')|l$ or~$(T')$], if for each $\gamma\in(0,1)$, condition
$(T)_\gamma|l$
is fulfilled.
\end{definition}

It is known that in dimensions $d\ge2$, condition
$(T')$ implies the existence of a deterministic $v \in\R^d \setminus
\{0\}$ such that
$P_0$-a.s. $\lim_{n \to\infty} \frac{X_n}{n} = v$, as well as a
central limit theorem for the
RWRE so that under the annealed law $P_0$,
\[
B^n_\cdot:= \frac{X_{{\lfloor\cdot n \rfloor}}-{\lfloor\cdot n
\rfloor}v}{n}
\]
converges in distribution to a Brownian motion in the Skorokhod space
$D([0,\infty)$, $\R^d)$ as $n \to\infty$; see, for instance, Theorem
4.1 in \cite{Sz-04}
for further details. Recently, in \cite{Be-09} the author has shown
that in
dimensions
larger than three, the above law of large numbers and central limit
theorem remain valid if condition $(T)_\gamma$ is
satisfied for some $\gamma\in(0,1)$. In addition, in \cite{Sz-04}
the author has proven
that if $\PP$ is uniformly elliptic, then in dimensions $d \geq2$,
for each $\gamma\in(0.5,1)$ and each $l \in\mathbb{S}^{d-1}$, condition
$(T)_\gamma\vert l$ is equivalent to $(T') \vert l$.
In \cite{DrRa-09b}, the authors pushed down this equivalence to
each $\gamma\in(\gamma_d,1)$, where $\gamma_d\in(0.366,0.388)$ is
decreasing with the dimension.
The first main result of
the present paper is a considerable improvement of these previous results
for dimensions larger than three.
%
\begin{theorem} \label{ballisticityCondEqThm}
Let $d \geq4$ and $\PP$ be uniformly elliptic.
Then for all $\gamma\in(0,1)$ and $l \in\mathbb{S}^{d-1}$, condition
$(T)_\gamma\vert l$ is equivalent to $(T')|l$.
\end{theorem}

\begin{figure}

\includegraphics{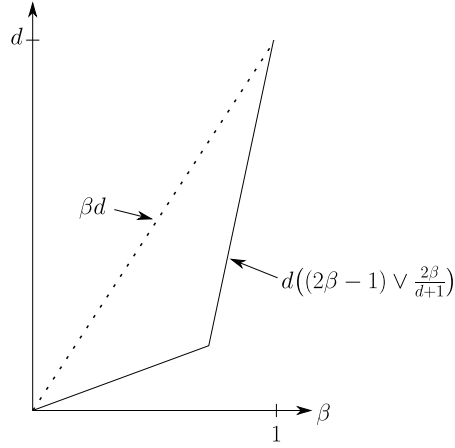}

\caption{Sketch of the known and conjectured bounds for $\alpha$.}
\label{fig1}
\end{figure}

The proof of Theorem
\ref{ballisticityCondEqThm} takes advantage of the effective criterion and
is therefore closely related to upper bounds for quenched probabilities
of atypical exit behavior of the RWRE.
To state the corresponding result, denote for
any subset $B \subset\Z^d$ its boundary by
\[
\partial B := \{x \in\Z^d \setminus B \dvtx\exists y \in B \mbox{ such
that } \Vert x-y \Vert_1 = 1\}
\]
and define the slab
\[
U_{\beta,l,L} := \{ x \in\Z^d \dvtx-L^\beta\leq x \cdot l \leq L\}.
\]
Furthermore, for the rest of this paper
we let
\[
T_B := \inf\{ n \in\N_0 \dvtx X_n \in B\}
\]
denote the first
hitting time. For $x \in\Z^d$ set $T_x := T_{\{x\}}$.
In terms of this notation, in \cite{Sz-04} the author conjectured
the following (cf. Figure \ref{fig1}).

\begin{conjecture} \label{conj:slabExitConj}
Let $d\geq2$, $\PP$ be uniformly elliptic and assume $(T') \vert l$ to
hold for some $l \in\mathbb{S}^{d-1}$.
Fix $c > 0$ and $\beta\in(0,1)$. Then for all $\alpha\in(0,\beta d)$,
\[
\limsup_{L \to\infty} L^{-\alpha} \log\PP\bigl(
P_{0,\omega}(X_{T_{\partial U_{\beta,l,L}}} \cdot l > 0) \leq
e^{-cL^\beta}
\bigr) < 0.\vadjust{\goodbreak}
\]
\end{conjecture}

Theorem 4.4 of \cite{Sz-04} states that the above conjecture holds true
for all positive $\alpha$ with
\[
\alpha< d \biggl( (2\beta-1) \vee\frac{2\beta}{d+1} \biggr).
\]

%
%

\begin{figure}

\includegraphics{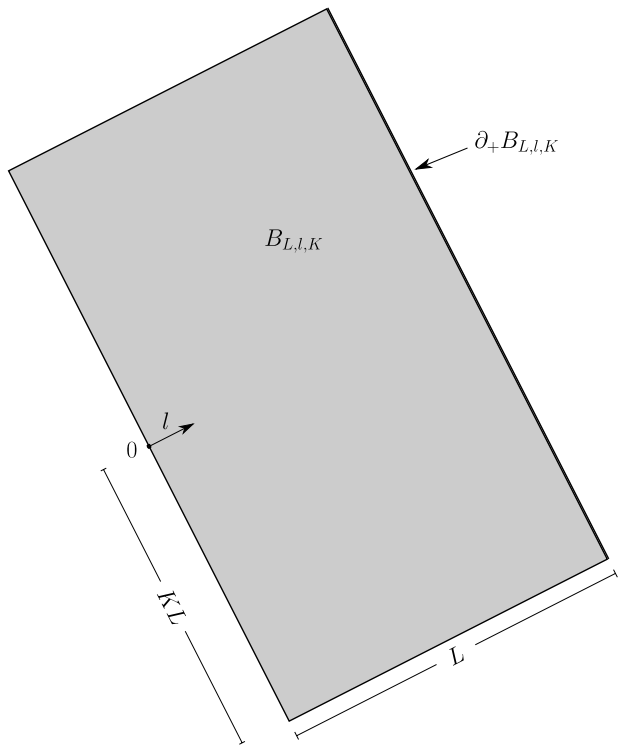}

\caption{The set $B_{L,l,K}$ and its boundary part $\partial_+
B_{L,l,K}$.}
\label{fig2}
\end{figure}

The second main result of the present
paper gives an affirmative answer to
a seemingly stronger statement than the one of
Conjecture \ref{conj:slabExitConj}.
For $l \in\mathbb{S}^{d-1}$,
denote by
\[
\pi_l\dvtx\R^d \ni x \quad\mapsto\quad(x \cdot l) l \in\R^d
\]
the orthogonal projection on the space $\{\lambda l \dvtx\lambda\in\R\}$
as well as by
\[
\pi_{l^\bot} \dvtx\R^d \ni x \quad\mapsto\quad x - \pi_l(x) \in\R^d
\]
the orthogonal projection on the orthogonal complement $\{\lambda l
\dvtx\lambda\in\R\}^\bot$.
Using this notation, for $K > 0$ we define the box
\[
B_{L,l,K} := \{x \in\Z^d \dvtx0 \leq x \cdot l \leq L, \Vert\pi
_{l^\bot}(x) \Vert_\infty\leq KL\}
\]
as well as
its right boundary part
%
%
\begin{equation} \label{eq:BPlusBdPart}
\partial_+ B_{L,l,K} := \{x \in\partial B_{L,l,K} \dvtx x \cdot l >
L\},
\end{equation}
see Figure \ref{fig2}.

%
%

We can now state the desired result.
%
\begin{theorem} \label{slabExitDistConjThm}
Let $d \geq4$, $\PP$ be uniformly elliptic and assume $(T)_\gamma
\vert l$ to hold for some
$\gamma\in(0,1)$, $l \in\mathbb{S}^{d-1}$. Fix $c > 0$ and $\beta
\in(0,1)$.
Then there exists a constant $K > 0$ such that for all $\alpha\in(0,
\beta d)$,
\[
\limsup_{L \to\infty} L^{-\alpha} \log\PP\bigl(
P_{0,\omega} (T_{\partial B_{L,l,K}} = T_{\partial_+ B_{L,l,K}}) \leq
e^{-cL^\beta} \bigr) < 0.
\]
\end{theorem}

\begin{remark}
\begin{longlist}[(a)]
\item[(a)]
The result we prove is slightly
stronger than the conjecture announced in \cite{Sz-04} since
we can dispose of the extent of the slab in direction $-l$ as well as
restrict the extent in directions
orthogonal to $l$. Scrutinizing the proof it will be clear that one can
improve this result
replacing the box $B_{L,l,K}$ by a~parabola-shaped set which grows in
the directions transversal to $\hat{v}$
at least like $L^\alpha$ for some $\alpha> 1/2$.

\item[(b)]
Note that this theorem is optimal in the sense that its conclusion will
not hold in general for
$\alpha>\beta d$. In fact, for plain nestling RWRE,
this can be shown by the use of so-called \textit{na\"{i}ve traps}
(see \cite{Sz-04}, page 244).

\item[(c)]
In both, Theorem \ref{ballisticityCondEqThm} as well as Theorem \ref
{slabExitDistConjThm},
the restriction to dimensions larger than three is caused by the
following: for a very large~set
of environments we need that the trajectories of two independent
$d$-dimensional random walks in this environment
intersect only very rarely; see equations~(\ref{eq:JNProb}) and~(\ref
{eq:JHittingProbBd}).
\end{longlist}
\end{remark}

The proof of Theorem \ref{slabExitDistConjThm} exploits heavily
a recent multiscale technique introduced in
\cite{Be-09} to study the slowdown upper bound for RWRE.
To explain this in more detail, note that from that source one also
infers that
every RWRE in a uniformly elliptic i.i.d. environment
which satisfies
condition~$(T)_\gamma$ for some $\gamma\in(0,1)$, has an asymptotic
speed $v \ne0$.
The main result of \cite{Be-09}
states that for every RWRE in a uniformly elliptic
i.i.d. environment satisfying condition $(T)_\gamma$, some $\gamma\in
(0,1)$, the following holds:
for each $a\ne v$ in the convex hull of
$0$ and $v$ as well as $\epsilon>0$ small enough,
and any $\alpha< d$ the inequality
\[
P_0\biggl(\biggl\|\frac{X_n}{n}-a\biggr\|_\infty< \epsilon\biggr)
\le\exp\{-(\log n)^\alpha\}
\]
holds for all $n$ large enough.
To prove the above result, Berger develops a~multiscale technique
which describes the behavior of the walk at the scale of the so called
na\"{i}ve traps, which at time $n$ are of radius of order $\log n$.
Here, we rely on such a multiscale technique to make explicit the
role of the regions of the same scale as the na\"{i}ve traps  to prove
Theorem \ref{slabExitDistConjThm}.

In Section \ref{section1}, we show how certain exit estimates from
boxes imply
Theorem~\ref{slabExitDistConjThm} and how in turn such a result
implies Theorem \ref{ballisticityCondEqThm}.
In Section \ref{sec:propProof},
we start with giving a heuristic explanation of a modified version of Berger's
multiscale technique and of how to deduce the aforementioned
exit estimates.
We then set up our framework of notation and auxiliary results before
making precise the previous heuristics by giving the corresponding
proofs. In the
\hyperref[sec:Appendix]{Appendix} we establish several specific results concerning local
limit theorem type results and estimates involving intersections of random
walks.

\section{Proofs of the main results}
\label{section1}
The proofs of Theorems \ref{ballisticityCondEqThm}
and \ref{slabExitDistConjThm} are based on a multiscale
argument and a semi-local limit theorem developed in~\cite{Be-09}
for RWRE in dimensions larger than or equal to four.

It is well known that if for some
$\gamma\in(0,1)$ and $l \in\mathbb{S}^{d-1}$,
condition $(T)_\gamma\vert l$ is fulfilled, then $P_0$-a.s. the limit
\[
\hat{v} := \lim_{n \to\infty} \frac{X_n}{\Vert X_n \Vert_2} \in
\mathbb{S}^{d-1}
\]
exists and is constant (cf., e.g., Theorem 1 in Simenhaus \cite
{Si-07}); it is called
the \textit{asymptotic direction}.

\begin{figure}

\includegraphics{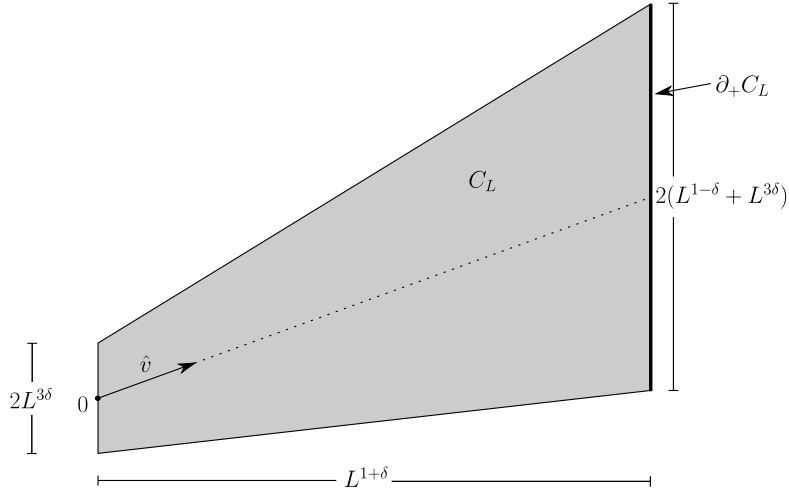}%
\vspace*{-3pt}
\caption{The set $C_L$ and its boundary part $\partial_+ C_L$.}
\label{fig3}
\vspace*{-3pt}
\end{figure}

Define for
a vector $e_j$ of the canonical basis of $\Z^d$ and
$l \in\mathbb{S}^{d-1}$ such that $l \cdot e_j \ne0$
the projection $\tilde{\pi}^j_{l}$ via
\[
\tilde{\pi}^j_{l}\dvtx\R^d \ni x \quad\mapsto\quad\frac{x \cdot e_j}{l
\cdot
e_j} l \in\R^d
\]
on the space $\{\lambda l \dvtx\lambda\in\R\}$
and
by $\tilde{\pi}^j_{l^\bot}$ the projection
\[
\tilde{\pi}^j_{l^\bot}\dvtx\R^d \ni x \quad\mapsto\quad x - \frac{x
\cdot
e_j}{l \cdot e_j} l \in\R^d
\]
on the space $\{\lambda e_j \dvtx\lambda\in\R\}^\bot$. In the case
$j=1$, we will abbreviate this notation by $\tilde{\pi}_{l}$ and $\tilde{\pi}_{l^\bot}$.
For $j \in\{1, \ldots, d\}$, $\delta> 0$ and $L > 0$, define the set
\[
C_L := \{ x \in\Z^d \dvtx0 \leq x \cdot e_j \leq L^{1+\delta},
\Vert\tilde{\pi}^j_{\hat{v}^\bot}(x) \Vert_\infty\leq
L^{3\delta} + x \cdot e_j L^{-2\delta} \};
\]
cf. Figure \ref{fig3}. In analogy to (\ref{eq:BPlusBdPart}),
we introduce the right boundary parts
\[
\partial_+ C_L := \{ x \in\partial C_L \dvtx x \cdot e_j > L^{1+\delta
}\}
\]
and $\partial_+ (x+C_L) := x+\partial_+ C_L$ for $x \in\Z^d$.

%
%

The proof of the following proposition will be deferred to Section \ref
{sec:propProof}.

\begin{proposition} \label{prop:quenchedConeTypeExit}
Let $d \geq4$, $\PP$ be uniformly elliptic and assume $(T)_\gamma
\vert l$ to hold for some
$\gamma\in(0,1)$, $l \in\mathbb{S}^{d-1}$. Without loss of
generality, let $e_j$ be a~vector of the
canonical basis such that
$\hat{v} \cdot e_j > 0$
and fix $\beta\in(0,1)$ as well as $\alpha\in(0,\beta d)$.

Then for all $\delta> 0$ small enough there exists
a sequence of events $(\Xi_{L})_{L \in\N}$
such that for all $L$ large enough we have
\[
\inf_{\omega\in\Xi_L}P_{0,\omega} ( T_{\partial C_{L}} =
T_{\partial_+ C_{L}}) \geq e^{-L^{\beta-\delta}}
\]
and
\[
\PP(\Xi_L^c) \leq e^{-L^\alpha}.
\]
\end{proposition}

For the sake of notational simplicity and
without loss of generality, we assume $j = 1$ from now on.

\subsection{\texorpdfstring{Proof of Theorem \protect\ref{slabExitDistConjThm}}{Proof of Theorem 1.5}}

We will show that Theorem \ref{slabExitDistConjThm} is a consequence of
Proposition~\ref{prop:quenchedConeTypeExit}.\vadjust{\goodbreak}
For this purpose, let the assumptions of Theorem~\ref
{slabExitDistConjThm} be fulfilled.
In particular, let $(T)_\gamma\vert l$ be fulfilled
(which implies $l\cdot\hat{v} > 0$; cf. Theorem 1.1 of \cite{Sz-02})
and
fix $c > 0$, $\beta\in(0,1)$ as well as $\alpha\in(0,\beta d)$.
Let $\delta> 0$ small enough such that the implication of Proposition
\ref{prop:quenchedConeTypeExit}
holds and $3\delta< \beta- \delta$. Choose
$\beta' \in(\beta-\delta,\beta)$
and define $x \in\Z^d$ to be one of the (possibly several) sites
closest to $L^{\beta'} l$.
Then the following property of the displaced set $x+C_L$ will
be used:
\begin{longlist}[(Exit)]
\item[(Exit)]\hypertarget{item:Exit}
Let $K$ be large enough and $\delta> 0$ small enough. Then for $L$
large enough,
if the walk starting in $x$ leaves $x+C_L$ through $\partial
_+(x+C_L)$, it
also leaves the box $B_{L,l,K}$ through $\partial_+ B_{L,l,K}$.
\end{longlist}
Now since the measure $\PP$ is uniformly elliptic, we know
that there exists a~constant $C$ depending on the dimension $d$, such
that for
all $L$ large enough and for $\PP$-a.a. $\omega$ the inequality
%
%
\begin{equation}
\label{eq:p1}
P_{0,\omega} (T_{\partial B_{L,l,K}} > T_x) \geq e^{-CL^{\beta'}}
\end{equation}
holds true.
By Proposition \ref{prop:quenchedConeTypeExit},
for $\alpha\in(0,\beta d)$ fixed, there
are subsets $\Xi_{L} \subset\Omega$ such that for $L$ large enough,
$\PP(\Xi_{L}) \geq1 - e^{-L^\alpha}$ and such that for $\omega\in
\Xi_{L}$ one has
\begin{eqnarray*}
&&P_{0,\omega} (T_{\partial B_{L,l,K}} = T_{\partial_+ B_{L,l,K}}) \\
&&\qquad\geq P_{0,\omega} \bigl(
T_{\partial B_{L,l,K}}> T_x, T_{\partial(x+C_L)}(\theta_{T_x}
(X_\cdot)) = T_{\partial_+ (x+C_L)}
(\theta_{T_x} (X_\cdot))\bigr)\\
&&\qquad\geq P_{0,\omega} (
T_{\partial B_{L,l,K}}> T_x)
P_{x,\omega} \bigl(T_{\partial(x+C_L)} = T_{\partial_+ (x+C_L)}
\bigr)\\
&&\qquad\geq e^{-CL^{\beta'}} e^{-L^{\beta- \delta}} = e^{-CL^{\beta}}
\end{eqnarray*}
for $L$ large enough,
where $\theta_n\dvtx(\Z^d)^{\N_0} \to(\Z^d)^{\N_0}$ denotes the
canonical $n$-fold left shift
and to obtain the first inequality we used property \hyperlink
{item:Exit}{(Exit)}.
In the second inequality, we have used the strong Markov property and
in the third one we employed inequality (\ref{eq:p1}) as well as
Proposition \ref{prop:quenchedConeTypeExit} in combination with the
translation invariance of the measure $\PP$.
This finishes the proof of the theorem.

\subsection{\texorpdfstring{Proof of Theorem \protect\ref{ballisticityCondEqThm}}{Proof of Theorem 1.3}}
In \cite{Sz-02}, the author introduces the so called \textit{effective
criterion},
which is a ballisticity condition equivalent to condition~$(T')$ and
which facilitates the explicit verification of condition~$(T')$.
The proof of Theorem \ref{ballisticityCondEqThm} will rest
on the fact that the effective criterion implies condition $(T')$.
Indeed, we will prove that $(T)_\gamma$ implies the effective criterion,
the main ingredient being Theorem~\ref{slabExitDistConjThm}.

For the sake of convenience, we recall here the effective criterion and its
features.
For positive numbers $L$, $L'$ and $\tilde{L}$ as well as a space
rotation $R$
around the origin we define the
\textit{box specification} ${\mathcal{B}}(R, L, L', \tilde{L})$
as the box
$B:= \{x\in\mathbb Z^d\dvtx x\in R((-L,L') \times(-\tilde{L}, \tilde
{L})^{d-1})\}$. Furthermore, let
\[
\rho_{\mathcal{B}}(\omega) := \frac{P_{0,\omega} (X_{T_{\partial
B}} \notin\partial_+ B)}{P_{0,\omega} (X_{T_{\partial B}} \in
\partial_+ B)}.
\]
Here,
$
\partial_+ B := \{x \in\partial B \dvtx R(e_1) \cdot x \geq L', \vert
R(e_j) \cdot x \vert< \tilde{L} \  \forall j \in\{2, \ldots, d\}\}.
$
We will sometimes write $\rho$ instead of $\rho_{\mathcal{B}}$ if
the box we refer to is clear from the context
and use $\hat{R}$ to label any rotation mapping $e_1$ to $\hat{v}$. Note
that due to the uniform ellipticity assumption, $\PP$-a.s. we have
$\rho\in
(0,\infty)$. Given $l\in\mathbb{S}^{d-1}$, we say that the \textit{effective
criterion with respect to $l$} is satisfied if
%
%
\begin{equation} \label{effectiveCritInf}
\inf_{{\mathcal{B}}, a} \biggl\{ c_{1}(d) \biggl(\log\frac{1}{\kappa
} \biggr)^{3(d-1)} \tilde{L}^{d-1} L^{3(d-1)+1} \E\rho_{\mathcal
{B}}^a \biggr\} < 1.
\end{equation}
Here, when taking the infimum,
$a$ runs over $[0,1]$ while ${\mathcal{B}}$ runs over the box-specifications
${\mathcal{B}}(R, L-2, L+2, \tilde{L})$ with $R$ a rotation such that
\mbox{$R(e_1) = l$},\break
$L \geq c_{2}(d)$, $3\sqrt{d} \leq\tilde{L} < L^3$. Furthermore,
$c_1(d)$ and $c_2(d)$ are dimension dependent constants.

The following result was proven in
\cite{Sz-02}.\vspace*{-2pt}
%
\begin{theorem} \label{Sz24}
$\!\!\!$For each $l\,{\in}\,\mathbb{S}^{d-1}$, the following
conditions are \mbox{equivalent}:

\begin{longlist}[(a)]
\item[(a)] The effective criterion with respect to $l$ is satisfied.

\item[(b)]
$(T')|l$ is satisfied.\vspace*{-2pt}
\end{longlist}
\end{theorem}

Due to this result, we can check condition $(T')$, which by nature of
its definition is asymptotic, by investigating the
local behavior of the walk only; indeed, to have the infimum on the
left-hand side of (\ref{effectiveCritInf})
smaller than $1$, it is sufficient to find
one box $\mathcal{B}$ and $a \in[0,1]$ such that the corresponding
inequality holds.

Recall that from Theorem 1.1 of \cite{Sz-02} we infer that for $l$
such that $l \cdot\hat{v} >0$, we
have that $(T)_\gamma\vert l$ implies $(T)_\gamma\vert\hat{v}$,
and $(T)' \vert\hat{v}$ implies
$(T)' \vert l$.
Thus,
because of (\ref{effectiveCritInf}) and Theorem \ref{Sz24}, in order
to prove Theorem \ref{ballisticityCondEqThm} it is then sufficient
to show that $(T)_\gamma\vert\hat{v}$ implies that
\begin{longlist}[(D)]
\item[(D)]\hypertarget{item:rhoExpDecay}
for every natural $n \in\N$, one has that
$
\E\rho^a = o(L^n)
$ as $L \to\infty$;
\end{longlist}
here, $\rho$ corresponds to a box specification $\mathcal{B} (\hat
{R},L-2, L+2, L^2)$.

To show the desired decay, we split
$
\E\rho^a
$
according to
%
%
\begin{equation}
\label{decomp}
\E\rho^a = \mathcal{E}_0 + \sum_{j=1}^{n-1} \mathcal{E}_j+\mathcal{E}_n,
\end{equation}
where $n=n(\gamma)$ is a natural number the choice of which will
depend on $\gamma$,
\begin{eqnarray*}
\mathcal{E}_0 &:=& \E\bigl( \rho^a, P_{0,\omega}(X_{T_{\partial B}}
\in\partial_+ B) > e^{-k_1 L^{\beta_1}} \bigr),
\\[-2pt]
\mathcal{E}_j &:=& \E\bigl( \rho^a, e^{-k_{j+1} L^{\beta_{j+1}}} <
P_{0,\omega}(X_{T_{\partial B}} \in\partial_+ B) \leq e^{-k_j
L^{\beta_j}} \bigr)
\end{eqnarray*}
for $j \in\{1, \ldots, n-1\}$ and
\[
\mathcal{E}_n := \E\bigl( \rho^a, P_{0,\omega}(X_{T_{\partial B}}
\in\partial_+ B) \le e^{-k_n L^{\beta_n}} \bigr)
\]
with parameters
\[
\gamma=: \beta_1 < \beta_2 < \cdots< \beta_n := 1,\vadjust{\goodbreak}
\]
$a = L^{-\varepsilon}$, $\varepsilon\in(0,1)$, as well as $k_n$
large enough and arbitrary
positive constants $k_1, k_2, \ldots, k_{n-1}$.
To bound $\mathcal{E}_0$, we employ the following lemma, which has
been proven in \cite{DrRa-09b}.
%
\begin{lemma}
\label{I} For all $L>0$,
\[
\mathcal{E}_0
\leq e^{k_1L^{\gamma- \varepsilon} - \delta_1 L^{\gamma-
\varepsilon}+o(L^{\gamma-\varepsilon})},
\]
where
\[
\delta_1:= -\limsup_{L \to\infty} L^{-\gamma} \log
P_0(X_{T_{\partial B}} \notin\partial_+ B) > 0.
\]
\end{lemma}

To deal with the middle summand in the right-hand side of (\ref
{decomp}), we
use the following lemma.
%
\begin{lemma}
\label{II} For all $L>0$, $j \in\{1, \ldots, n\}$ and $\varepsilon
>0$, we have that
\[
\mathcal{E}_j
\leq e^{k_{j+1}L^{\beta_{j+1} - \varepsilon} - L^{\beta_jd -
\varepsilon}+o(L^{\beta_jd-\varepsilon})}.
\]
\end{lemma}
\begin{pf}
Using Markov's inequality, for $j \in\{1, \ldots, n-1\}$ we obtain
the estimate
%
%
\begin{equation} \label{(II)Est}
\mathcal{E}_j\leq e^{k_{j+1} L^{\beta_{j+1} -\varepsilon}} \PP\bigl(
P_{0,\omega}(X_{T_{\partial B}} \in\partial_+ B) \leq e^{-k_j
L^{\beta_j}} \bigr).
\end{equation}
Due to Theorem
\ref{slabExitDistConjThm}, for $\varepsilon> 0$ fixed,
the outer probability on the right-hand side of (\ref{(II)Est}) can be
estimated from
above by
$
e^{-L^{\beta_j d - \varepsilon} + o(L^{\beta_j d-\varepsilon})}.
$
\end{pf}

For the term $\mathcal{E}_n$ in (\ref{decomp}), we have the following
estimate.
%
\begin{lemma}
\label{III} There exists a constant $C>0$ such that for any
$\varepsilon> 0$,
\[
\mathcal{E}_n
\leq e^{CL^{1- \varepsilon} - L^{\beta_nd - \varepsilon}+o(L^{\beta
_nd-\varepsilon})}.
\]
\end{lemma}
\begin{pf} Using the uniform ellipticity assumption, we
see that there is a~constant $C>0$ such that
%
%
\begin{equation} \label{(III)Est}
\mathcal{E}_n\leq e^{C L^{1- \varepsilon}} \PP\bigl( P_{0,\omega
}(X_{T_{\partial B}} \in\partial_+ B) \leq e^{-k_n L^{\beta_n}} \bigr).
\end{equation}
An application of Theorem \ref{slabExitDistConjThm} to estimate
the second factor of the right-hand side of inequality (\ref{(III)Est})
establishes the proof.
\end{pf}

From Lemmas \ref{I}, \ref{II} and \ref{III}, we deduce that
for $k_1 <\delta_1$, $n = n(\gamma)$ large enough, arbitrarily chosen
positive constants $k_2,\ldots,k_n$
as well as $\varepsilon$ and $\beta_1,\ldots,\beta_n$ satisfying
\begin{eqnarray*}
\beta_1 &=& \gamma, \qquad\varepsilon<\gamma,
\\
\beta_{j+1} &<& \beta_j d\vadjust{\goodbreak}
\end{eqnarray*}
for $j \in\{1, \ldots, n-1\}$, and
\[
1 < \beta_n d,
\]
all the terms $\mathcal{E}_0,\ldots, \mathcal{E}_n$ on the
right-hand side of (\ref{decomp})
decay stretched exponentially.
It is easily observed that the above choice of parameters is feasible,
which establishes the desired decay
in \hyperlink{item:rhoExpDecay}{(D)} and thus finishes the proof of
Theorem \ref{ballisticityCondEqThm}.

\section{\texorpdfstring{Proof of Proposition \protect\ref{prop:quenchedConeTypeExit} and auxiliary results}
{Proof of Proposition 2.1 and auxiliary results}} \label{sec:propProof}

The proof of Proposition~\ref{prop:quenchedConeTypeExit} is based
on a modified version of the multiscale argument developed in~\cite{Be-09}.
In general, in our construction, we will name the corresponding results
of the construction
in \cite{Be-09} in brackets in the corresponding places.

We start with giving a heuristic (and cursory) idea of the proof.
Afterward, we will set up all the necessary notation and auxiliary
results before providing a rigorous proof of Proposition \ref
{prop:quenchedConeTypeExit}.

\subsection{\texorpdfstring{Heuristics leading to Proposition \protect\ref{prop:quenchedConeTypeExit}}{Heuristics leading to Proposition 2.1}}
\label{subsec:heuristics}

The basic strategy of the proof is to
construct, for $\beta\in(0,1)$ and $\alpha<\beta d$ given, a
sequence of events $(G_L)_{L \in\N}$, each a subset of $\Omega$,
such that
for $L$ large enough
one has
%
%
\begin{equation}
\label{oneh}
\PP(G_L^c)\le e^{-L^{\alpha}}
\end{equation}
and at the same time
%
%
\begin{equation}
\label{twoh}
\inf_{\omega\in G_L} P_{0,\omega}(T_{\partial_+ C_L} = T_{\partial
C_L})\ge e^{-cL^{\beta}},
\end{equation}
where $c$ is a constant that changes values various times throughout
this subsection.
In order to define $G_L$, for each of finitely many scales, we
cover the box $C_L$ with boxes of that certain scale. Boxes of the
first scale
have extent roughly $L^{2\psi}$ in direction $\hat{v}$, and extent
marginally larger than
$L^\psi$ in directions orthogonal to $\hat{v}$. Here, $\psi> 0$ is
much smaller than $\beta$. The boxes of larger scale more or less have
$\psi$ replaced by larger numbers
[see (\ref{def:PBox}), (\ref{eq:scales}) and (\ref{eq:scaledPs})].
Given an environment, we declare a box to be \textit{good} if
within this box and with respect to the given environment, the quenched
random walk behaves very much like the
annealed one.
Otherwise, it is called \textit{bad}.

We then define
$G_L$ as the event that there are not significantly more than~$L^{\alpha}$ bad boxes
of each scale contained in $C_L$. Using Proposition \ref{prop:closenessbase},
which states that the probability of a box being bad decays faster than
polynomially as a function in $L$, by large deviations for
binomially distributed variables one shows that
the probability of the complement of this event is smaller than
$e^{-L^{\alpha}}$, so that (\ref{oneh}) is satisfied (cf. Lemma \ref
{lem:GProbEstimate}).

It remains to show that on $G_L$, inequality (\ref{twoh}) is satisfied.
For this purpose, we associate to the walk a ``current scale''
that slowly increases as the $e_1$-coordinate of the walk increases.
We will then require the walk to essentially leave
in $e_1$-direction (i.e., through their right boundary parts) all the\vadjust{\goodbreak}
boxes of its current scale it traverses; this ensures that
it leaves $C_L$ through~$\partial_+ C_L$.
Since the
probability that the random walk exits a good box through the right
boundary part is relatively large,
one can essentially bound the probability of leaving $C_L$ through~$\partial_+ C_L$ from below by
the cost the walk incurs when traversing bad boxes.

Now each time the walk finds itself in a bad box of its current scale,
it will instead move in boxes of smaller scale that contain its current
position, and leave these boxes through their right boundary
parts.
Each time this happens, it has to ``correct'' the errors incurred by
moving in such boxes through some deterministic steps,
the cost of which will not exceed $e^{-cL^{2\psi}}$; in a certain way,
these corrections make the walk look as if it has been leaving
a box of its current scale through its right boundary part.
Thus, we can roughly bound the probability of leaving $C_L$ through
$\partial_+ C_L$ by
%
%
\begin{equation} \label{eq:badCosts}
e^{-cNL^{2\psi}},
\end{equation}
where $N$ is the number of bad boxes that the walk visits.
%

Now in order to obtain a useful upper bound for $N$, we can force the
random walk to have CLT-type fluctuations
in directions transversal to~$\hat{v}$
at constant cost in each box (see \textit{random direction event},
Section \ref{subsec:RDE}).
By means of this random direction event, one can then infer the
existence of a~direction (depending on the environment)
such that, if the CLT-type fluctuations
of the walk essentially center around this direction, then the walk
encounters a~little less than $L^\beta$ bad boxes of each scale
on its way through~$C_L$.
From~(\ref{eq:badCosts}), we deduce that the probability for the
walker to leave $C_L$ through~$\partial_+ C_L$
can then be bounded from below by $e^{-cL^\beta}$.
This suggests that~(\ref{twoh}) holds.

\subsection{Preliminaries}

We first recall an equivalent formulation of condition~$(T)_\gamma$ and
introduce the basic notation that will be used throughout
the rest of this paper.

We will use $C$ to denote a generic constant that may
change from one side to the other of the same inequality.
This constant may usually depend on various parameters,
but in particular does not depend on the variable $L$ nor $N$
(recall that $L$ is the variable which makes the slabs and boxes grow,
and $N$ will play a similar role
in general results).
In ``general lemmas,'' we will usually denote the corresponding
probability measure
and expectation by $P$ and $E$, respectively. Furthermore,
when considering stopping times without mentioning the process they
apply to,
then they will usually refer to the RWRE $X$.

Not all auxiliary results will appear in the order in which they are
employed. In fact,
in order to improve readability, we defer the majority of them to the
\hyperref[sec:Appendix]{Appendix}.\vspace*{8pt}

\textit{In addition, we assume the conditions of Proposition} \ref
{prop:quenchedConeTypeExit}
\textit{to be fulfilled for the rest of this paper without further
mentioning.}\vadjust{\goodbreak}

We first introduce the regeneration times in direction $e_1$.
Setting $\tau_0 := 0$, we define the \textit{first regeneration time}
$\tau_1$
as the first time $X_n \cdot e_1$ obtains a~new maximum and never falls
below that maximum again,
that is,
\[
\tau_1:=\inf\Bigl\{n \in\N\dvtx\sup_{0\le k\le n-1}X_k\cdot e_1 <
X_n\cdot e_1 \mbox{ and }
\inf_{k\ge n} X_k\cdot e_1 \ge X_n \cdot e_1 \Bigr\}.
\]
Now define recursively in $n$ the \textit{$(n+1)$st regeneration time}
$\tau_{n+1}$
as the first time after $\tau_n$ that $X_n\cdot l$ obtains a new maximum
and never goes below that maximum again, that is, $\tau_{n+1} := \tau
_{1}(X_{\tau_n + \cdot})$.
For $n\in\N$, we define
the \textit{radius} of the $n$th regeneration as
\[
X^{*(n)}:=\sup_{\tau_{n-1} \leq k \leq\tau_n} \Vert X_k - X_{\tau
_{n-1}} \Vert_1.
\]
This notation gives rise to the following equivalent formulation of
$(T)_\gamma$ proven in~\cite{Sz-02}, Corollary 1.5.
%
\begin{theorem}
\label{equiv} Let $\gamma\in(0,1)$ and $l\in\mathbb S^{d-1}$.
Then the following are equivalent:
\begin{longlist}
\item Condition $(T)_\gamma|l$ is satisfied.
\item\hypertarget{item:integralCond}
$P_0(\lim_{n\to\infty} X_n\cdot l=\infty)=1$ and $E_0 \exp\{
c(X^{*(1)})^\gamma\}<\infty$
for some \mbox{$c>0$}.
\end{longlist}
\end{theorem}
%
\begin{remark}
Note in particular that, similarly to Proposition 1.3
of Sznitman and Zerner \cite{SzZe-99}, condition \hyperlink
{item:integralCond}{(ii)} implies $E_0 \exp\{ c(X^{*(n)})^\gamma\}
<\infty$ for any $n \ge2$.
\end{remark}

\begin{figure}

\includegraphics{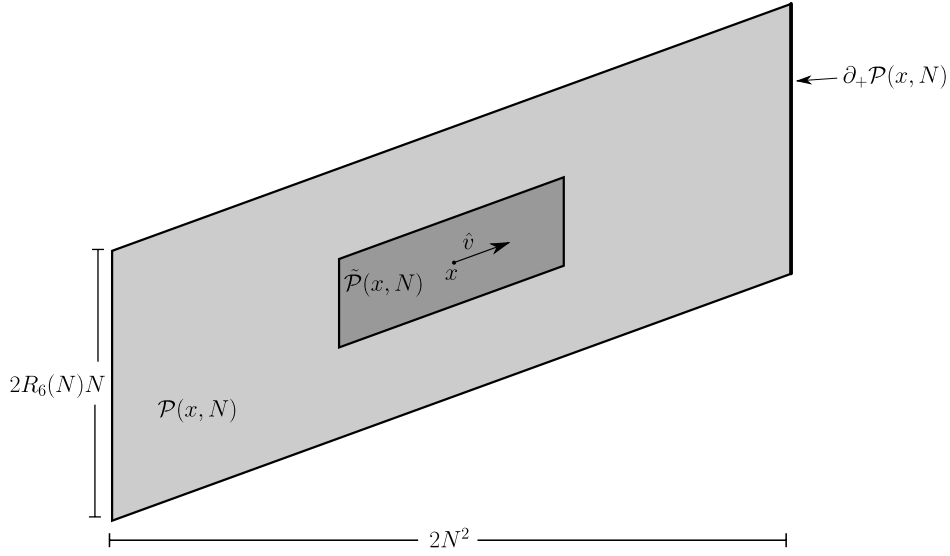}

\caption{The set $\mathcal{P}(x,N)$ and its right boundary part
$\partial_+ \mathcal{P}(x,N)$ as well as its middle third
$\tilde{\mathcal{P}}(x,N)$.}
\label{fig4}
\end{figure}

We will repeatedly use the above equivalence.
Now for each natural~$k$ and~$N$ we define the scales
\[
R_k(N) := \lceil\exp\{ (\log\log N)^{k+1}\} \rceil.
\]
Note that for every natural $n,N$ and $k$ one
has that
\[
R_k^n(N) \in o(R_{k+1}(N)) \quad\mbox{and}\quad R_{k}(N) \in o(N).
\]
Define for each natural $N$ the
sublattice
\[
\mathcal{L}_N := N^2 \Z\times\biggl\lfloor\frac{R_6(N)N}{4}
\biggr\rfloor\Z^{d-1}
\]
of $\Z^d$.
Furthermore, for each $N$ and $x \in\Z^d$ we define the blocks
%
%
\begin{equation} \label{def:PBox}
\qquad{\mathcal{P}}(0,N) := \{y \in\Z^d \dvtx-N^2 <y \cdot e_1 < N^2,
\Vert\tilde{\pi}_{\hat{v}^\bot}(y) \Vert_\infty
< R_6(N)N\}
\end{equation}
and
\[
{\mathcal{P}}(x,N) := x + {\mathcal{P}}(0,N)
\]
as well as their middle thirds
\[
\tilde{{\mathcal{P}}}(0,N) := \{y \in\Z^d \dvtx-N^2/3 <y \cdot
e_1 < N^2/3,
\Vert\tilde{\pi}_{\hat{v}^\bot}(y) \Vert_\infty
< R_6(N)N / 3\}\vadjust{\goodbreak}
\]
and
\[
\tilde{{\mathcal{P}}}(x,N) := x + \tilde{{\mathcal{P}}}(0,N).
\]
Note that this construction ensures that
for each $x\in N^2\Z\times\Z^{d-1}$ there exists a $z\in\mathcal
{L}_N$ such that
$x\in{\tilde{\mathcal{P}}}(z,N)$.
Furthermore, define its right boundary part
\[
\partial_+ \mathcal{P}(x,N) := \{ y \in\partial\mathcal
{P}(x,N) \dvtx(y-x) \cdot e_1 = N^2 \}.
\]
See Figure \ref{fig4} for an illustration.

%
%

For $N\ge1$, define the event
\[
A_N(X) := \bigl\{X^{*(n)} < R_2(N) \  \forall n \in\{1, \ldots, 2N^2\}
\bigr\},
\]
where at times we write $A_N$ instead of $A_N(X)$ if the corresponding
process~$X$ is clear from the context.
Using Markov's inequality, the following lemma is a consequence of
Theorem \ref{equiv}.
%
\begin{lemma} \label{lem:ANEst}
There exists a constant $C>0$ such that for each $N\ge1$,
%
%
\begin{equation} \label{eq:ANEst}
P_0(A_N^c) \leq Ce^{-C^{-1}R_2(N)^\gamma}
\end{equation}
and, defining the event
\[
\mathcal{A}_N := \bigcap_{x \in\tilde{\mathcal{P}}(0,N)} \bigl\{
P_{x,\omega}(A_N^c) \leq e^{-R_1(N)^\gamma} \bigr\},
\]
which is contained in the Borel-$\sigma$-algebra of $\Omega$,
one has
\[
\PP(\mathcal{A}_N^c) \leq Ce^{-C^{-1} R_2(N)^\gamma}.\vadjust{\goodbreak}
\]
\end{lemma}

We define the set of rapidly decreasing sequences as
\[
\mathcal S(\N):= \Bigl\{ (a_n)_{n \in\N} \in\R^\N\dvtx\sup_{n \in
\N}|n^ka_n|<\infty\ \forall k \in\N\Bigr\}
\]
and note that due to Lemma \ref{lem:ANEst} we have that $N \mapsto P_0(A_N^c)$
and $N \mapsto\PP(\mathcal{A}_N^c)$ are contained in
$\mathcal{S}(\N)$.\vspace*{-2pt}

\subsection{Berger's semi-local limit theorem and scaling}
As a first step in the scaling, we introduce a classification of
blocks. We need to define some
parameters which will remain fixed throughout this paper.
For $\beta$ and $\alpha$ as in the assumptions of Proposition \ref
{prop:quenchedConeTypeExit},
choose $\delta$ such that
\[
0 < \delta< \frac{\beta d - \alpha}{12d}.
\]
%
%

Furthermore, fix
%
%
\begin{equation} \label{eq:psiDef}
\psi\in\biggl(2 \delta, \frac{20 \delta}{9} \biggr)
\end{equation}
%
%
%
and $\chi$ such that
\[
0 < \chi<(\beta- 6\delta)/2 \wedge\psi/4\wedge 6/(d-1).
\]
%
%
%
%
%
%
From now on let $L \in\N$,
define $L_1:={\lfloor L^\psi\rfloor}$ and recursively
in $k$ the scales
%
%
\begin{equation} \label{eq:scales}
L_{k+1}:=L_k {\lfloor L^\chi\rfloor}.
\end{equation}
Define
$\iota$ to be the smallest $k$ such that $L_k^2> L^{1+\delta}$.
For $k \in\{1, \ldots, \iota\}$ and $x \in C_L \cap\mathcal{L}_{L_k}$,
we call a block
%
%
\begin{equation} \label{eq:scaledPs}
{\mathcal{P}}(x,L_{k})
\end{equation}
\textit{good} with respect to the environment $\omega$ if
the following three
properties are
satisfied
for $\vartheta:= \chi$ and all $z \in\tilde{{\mathcal{P}}}(x,L_k)$:
\begin{longlist}
\item
%
%
\begin{equation}
 \label{eq:quenProbRightExitEst}
P_{z, \omega} \bigl( T_{\partial{\mathcal{P}}(x,L_k)} \not=
T_{\partial_+ {\mathcal{P}}(x,L_k)} \bigr) \leq
e^{-R_1(L_k)^\gamma}.
\end{equation}

\item
%
%
\begin{eqnarray}
\label{eq:quenAnnExpClose}
&&\bigl\Vert E_{z,\omega} \bigl( X_{T_{\partial{\mathcal{P}}(x,L_k)}}
\vert T_{\partial{\mathcal{P}}(x,L_k)} = T_{\partial_+
{\mathcal{P}}(x,L_k)} \bigr) \nonumber\\[-8pt]\\[-8pt]
&&\hspace*{0pt}\qquad{} - E_{z} \bigl( X_{T_{\partial{\mathcal
{P}}(x,L_k)}} \vert
T_{\partial{\mathcal{P}}(x,L_k)} = T_{\partial_+ {\mathcal
{P}}(x,L_k)} \bigr)
\bigr\Vert_1 \leq R_4(L_k).\nonumber
\end{eqnarray}

\item
%
%
\begin{eqnarray}
\label{eq:quenAnnProbClose}
&&\max_Q
\bigl\vert P_{z,\omega} \bigl( X_{T_{\partial{\mathcal{P}}(x,L_k)}}
\in Q \vert
T_{\partial{\mathcal{P}}(x,L_k)} = T_{\partial_+ {\mathcal
{P}}(x,L_k)} \bigr)\nonumber\\
&&\qquad\hspace*{0pt}{} - P_{z} \bigl( X_{T_{\partial{\mathcal
{P}}(x,L_k)}} \in Q
\vert T_{\partial{\mathcal{P}}(x,L_k)} = T_{\partial_+ {\mathcal
{P}}(x,L_k)} \bigr)
\bigr\vert\\
&&\qquad\quad< L_k^{(\vartheta- 1)(d-1) - \vartheta
({d-1})/({d+1})},\nonumber
\end{eqnarray}
where the maximum in $Q$ is taken over all
$(d-1)$-dimensional hypercubes $Q \subset\partial_+ {\mathcal{P}}(x,L_k)$
of side length $\lceil{L_k^\vartheta} \rceil$.\vadjust{\goodbreak}
\end{longlist}
Otherwise, we say that the block $\mathcal P(x,L_k)$ is
\textit{bad}.
For $k \in\{1, \ldots, \iota\}$ we will usually refer to boxes of the
form ${\mathcal{P}}(x,L_k)$ as a \textit{box of scale $k$.}

The following result is essentially Proposition 4.5 of \cite{Be-09}, which
can be understood as a semi-local central limit theorem
for RWRE. For the sake of completeness, we will give its proof in the
\hyperref[sec:Appendix]{Appendix}.\vspace*{-3pt}
%
\begin{proposition}[(Proposition 4.5 of \cite{Be-09})] \label
{prop:closenessbase}
Assume that $(T)_\gamma\vert l$ is satisfied and fix $\vartheta\in
(0, \frac{6}{d-1} \wedge 1)$. Then
there exists a sequence of events $(G_L)_{L \in\N} \subset\Omega$
such that $\PP(G_\cdot^c) \in\mathcal S(\N)$
and for all $\omega\in G_L$ and $k \in\{1, \ldots, \iota\}$:\vspace*{-1pt}
\begin{longlist}
\item\hypertarget{item:quenchedBadExitEst} display (\ref
{eq:quenProbRightExitEst}),

\item\hypertarget{item:quenchedAnnealedExpDiff} display (\ref
{eq:quenAnnExpClose}) and

\item\hypertarget{item:quenchedAnnealedProbDiff} display (\ref
{eq:quenAnnProbClose})\vspace*{-1pt}
\end{longlist}
are satisfied
for all $x \in C_L \cap\mathcal{L}_{L_k}$, $z \in\tilde{{\mathcal
{P}}}(x,L_k)$
and the chosen $\vartheta$.

In particular, due to the translation invariance of the environment,
we have that $\PP(\mathcal P(x, \cdot) \mbox{ is bad} ) \in\mathcal
S(\N)$ for any $x \in\Z^d$.\vspace*{-3pt}
\end{proposition}
%
\begin{remark} \label{rem:PropProofStrat}
For the sake of notational simplicity, we will prove the proposition by showing
that there exist sequences $G_L^{\mbox{\hyperlink
{item:quenchedBadExitEst}{\iitem}}}$,
$G_L^{\mbox{\hyperlink{item:quenchedAnnealedExpDiff}{\ii}}}$ and
$G_L^{\mbox{\hyperlink{item:quenchedAnnealedProbDiff}{\iii}}}$, $L \in\N
$, of subsets of $\Omega$ such that
\[
\PP\bigl({G_L^{\mbox{\hyperlink{item:quenchedBadExitEst}{\iitem}}}}^c \bigr
),\qquad
\PP\bigl({G_L^{\mbox{\hyperlink{item:quenchedAnnealedExpDiff}{\ii
}}}}^c\bigr)
\quad\mbox{and}\quad
\PP\bigl({G_L^{\mbox{\hyperlink{item:quenchedAnnealedProbDiff}{\iii
}}}}^c\bigr)\vspace*{-1pt}
\]
are contained in $\mathcal{S}(\N)$ as functions in $L$ and such that for
$\omega$ contained
in these sets, $x= 0$, and $z \in\tilde{\mathcal{P}}(0,L)$, displays
(\ref{eq:quenProbRightExitEst}),
(\ref{eq:quenAnnExpClose}) and
(\ref{eq:quenAnnProbClose}), respectively, are fulfilled for $L$
instead of $L_k$.
The required result then follows by observing that $\PP$ is translation
invariant
and using $\vert C_L \vert\leq CL^{2d}$ in combination with a standard
union bound.\vspace*{-3pt}
\end{remark}

We next give an upper bound on the probability that an environment has
many bad blocks.
For this purpose, set
%
%
\begin{eqnarray}
\label{eq:ThetaDef}
&&\Theta_L:= \bigl\{ \omega\in\Omega\dvtx\vert\{ x \in C_L \cap
\mathcal{L}_{L_k} \dvtx{\mathcal{P}}(x,L_k)\nonumber\hspace*{-35pt}\\[-10pt]\\[-10pt]
&&\hspace*{77.3pt}
\mbox{ is bad with respect to } \omega\} \vert
\leq L^{\alpha+ \delta}\ \forall k \in\{1, \ldots, \iota\} \bigr\}
.\nonumber\vspace*{-1pt}\hspace*{-35pt}
\end{eqnarray}
Furthermore, observe that $\mathcal{L}_L$ can be represented as the
disjoint union of
$2 \cdot8^{d-1}$ (translated) sublattices of $\Z^d$ such that for any
sublattice $\mathcal{L}$
of these and $z_1, z_2 \in\mathcal{L}$, we have $\mathcal{P}(z_1,L)
\cap\mathcal{P}(z_2,L) = \varnothing$.\vspace*{-2pt}
%
\begin{lemma}[(Lemma 5.1 of \cite{Be-09})] \label{lem:GProbEstimate}
For $L$ large enough,
\[
\PP(\Theta_L^c) \leq e^{-L^\alpha}.\vspace*{-3pt}
\]
\end{lemma}
\begin{pf}
For $k \in\{1, \ldots, \iota\}$,
set
\[
J_{L_k}(\omega) := \vert\{ z \in C_L \cap\mathcal{L}_{L_k} \dvtx
{\mathcal{P}}(z,L_k) \mbox{ is bad with respect to } \omega\}
\vert\vspace*{-1pt}
\]
and note that
%
%
\begin{equation}\label{eq:boundk}
\PP(\Theta_L^c) \leq\sum_{k=1}^\iota\PP(J_{L_k} > L^{\alpha+
\delta}).\vspace*{-1pt}\vadjust{\goodbreak}
\end{equation}
As in \cite{Be-09} we can write $J_{L_k} = J_{L_k}^{(1)} + \cdots+
J_{L_k}^{(2 \cdot8^{d-1})}$ with
$
J_{L_k}^{(m)}
$
distributed binomially with parameters $D({L_k})$ and $p({L_k})$
for $m \in\{1, \ldots, 2\cdot8^{d-1}\}$.
Here,
$p(L) := \PP( \mathcal{P}(0,L) \mbox{ is bad})$, that is, in
particular, due to Proposition \ref{prop:closenessbase},
%
%
\begin{equation}\label{eq:pk}
p \in\mathcal{S}(\N),
\end{equation}
and $D(L_k)$ is the maximal number of intersection points any of the
above-mentioned translated sublattices has
with $C_L$, that is, in particular
%
%
\begin{equation}\label{eq:dk}
D(L_k)\le CL^{2d}
\end{equation}
for some constant $C$ and all $L$.
Now for $m \in\{1, \ldots, 2 \cdot8^{d-1}\}$,
we have
%
%
\begin{equation}\label{eq:jk}
\PP\biggl( J_{L_k}^{(m)} > \frac{L^{\alpha+\delta}}{2 \cdot8^{d-1}}
\biggr) \leq
\exp\biggl\{-\frac{L^{\alpha+\delta}}{2 \cdot8^{d-1}} \biggr\} \E
\exp\bigl\{ J_{L_k}^{(m)}\bigr\}
\end{equation}
with
%
%
\begin{eqnarray}\label{eq:expk}
\E\exp\bigl\{ J_{L_k}^{(m)}\bigr\}
&\leq&\sum_{j=0}^{D(L_k)} \pmatrix{D({L_k}) \cr j} (ep({L_k}))^j
\bigl(1-ep({L_k})\bigr)^{D({L_k})-j}\nonumber\\[-8pt]\\[-8pt]
&&\hspace*{22.6pt}{}\times\biggl( \frac{1-p({L_k})}{1-ep({L_k})}
\biggr)^{D({L_k})-j},\nonumber
\end{eqnarray}
and from (\ref{eq:pk}) and (\ref{eq:dk}) we conclude that
\[
\lim_{L\to\infty}\biggl(\frac{1-p({L_k})}{1-ep({L_k})}
\biggr)^{D({L_k})-j} =1
\]
uniformly in $j \in\{0, \ldots, D(L_k)\}$.
Substituting this back into displays (\ref{eq:expk}), (\ref{eq:jk})
and (\ref{eq:boundk}),
we conclude the proof.
\end{pf}

We now need to recall the concept of \textit{closeness} between two
probability measures introduced in \cite{Be-09}.
Here and in the following, if $Z$ is a $d$-dimensional random variable
defined on a probability space
with probability measure~$\mu$, we write $E_\mu Z := \int Z \,d\mu$
and if
$\mu$ is a measure on $\R^d$, then we write $E_\mu:= \int x \,d\mu$,
whenever the integrals are well defined.
Furthermore, we define its variance
via $\Var Z := E \Vert Z - EZ \Vert_1^2$ whenever this expression is
well defined
and correspondingly for a
probability measure $\mu$ on $\R^d$ with appropriate integrability conditions
we write $\Var_\mu$.
%
\begin{definition} \label{def:closeness}
Let $\mu_1$ and $\mu_2$ be two probability
measures on $\Z^d$. Let $\lambda\in[0,1)$ and $K$ be a natural
number. We say that $\mu_2$ is \textit{$(\lambda, K)$-close}
to $\mu_1$ if there exists a coupling $\mu$ of three random variables
$Z_1$, $Z_2$ and $Z_0$ such that:
\begin{longlist}[(a)]
\item[(a)]\hypertarget{item:C1}
$\mu\circ Z_j^{-1} = \mu_j$ for $j\in\{1,2\}$,

\item[(b)]\hypertarget{item:C2} $\mu(Z_1 \not= Z_0)
\leq\lambda$,

\item[(c)]\hypertarget{item:C3} $\mu(\Vert Z_0 - Z_2 \Vert_1 \leq
K) = 1$,\vadjust{\goodbreak}

\item[(d)]\hypertarget{item:C4} $E_\mu Z_1 = E_\mu Z_0$,

\item[(e)]\hypertarget{item:C5} $\sum_x \Vert x - E_\mu Z_1 \Vert
_1^2 \cdot\vert\mu(Z_1 = x) - \mu(Z_0 = x) \vert\leq\lambda\Var Z_1$.
\end{longlist}
\end{definition}
%
\begin{remark} \label{rem:couplingSpace}
Assume given a random variable $X$ that is distributed according to
some distribution
which is $(\lambda,K)$-close to some other distribution.
Then the corresponding coupling which establishes this closeness
can be defined on an extension of the probability space
$X$ is defined on, with $X$ playing the role of $Z_2$.
We will therefore assume this property to be fulfilled
from now on without further mentioning when dealing with such
couplings.
\end{remark}

\subsection{General auxiliary results}
The following lemma is a sort of remedy for the fact that
\[
\Biggl( \tilde{\pi}_{\hat{v}^\bot} \Biggl(\sum_{j=1}^n (X_{\tau_j}
- X_{\tau_{j-1}}) \Biggr) \Biggr)_{n \in\{ 2, \ldots, 2L^2 \}}
\]
with respect to $P_0( \cdot \vert A_L)$,
due to the conditioning on $A_L$, is not a martingale.
To state the result, set
%
%
\begin{equation} \label{eq:vLDef}
\hat{v}_L := \frac{E_0 (X_{\tau_2} - X_{\tau_1})\mathbh
{1}_{A_L}}{\Vert E_0 (X_{\tau_2} - X_{\tau_1}) \mathbh{1}_{A_L}
\Vert_2}.
\end{equation}
We start with showing that for $L$ large, $\hat{v}_L$ hardly deviates
from the asymptotic direction $\hat{v}$.
%
\begin{lemma} \label{lem:directionsDistance}
\[
\Vert\hat{v} - \hat{v}_{\cdot} \Vert_2 \in\mathcal{S}(\N).
\]
\end{lemma}

\begin{pf}
Note that
%
%
\begin{eqnarray}\label{vAvDiff}
\Vert\hat{v} - \hat{v}_{L} \Vert_2 &=&
\biggl\Vert\frac{E_0(X_{\tau_2} - X_{\tau_1})}{\Vert E_0(X_{\tau_2}
- X_{\tau_2}) \Vert_2}
- \frac{ E_0 (X_{\tau_2} - X_{\tau_1}, A_L)}{ \Vert E_0 (X_{\tau_2}
- X_{\tau_1}, A_L) \Vert_2} \biggr\Vert_2\nonumber\\
&=& \bigl\Vert E_0(X_{\tau_2} - X_{\tau_1}) \Vert E_0
(X_{\tau_2} - X_{\tau_1}, A_L) \Vert_2\nonumber\\[-8pt]\\[-8pt]
&&\hspace*{2.4pt}{} - E_0 (X_{\tau_2}-X_{\tau_1}, A_L) \Vert E_0 (X_{\tau_2} -
X_{\tau_1}) \Vert_2 \bigr\Vert_2\nonumber\\
&&\hspace*{0pt}{}\times\bigl(\Vert E_0(X_{\tau_2}-X_{\tau_1}) \Vert_2 \Vert E_0( X_{\tau
_2}-X_{\tau_1}, A_L) \Vert_2\bigr)^{-1}.
\nonumber
\end{eqnarray}
Inserting a productive $0$, the numerator evaluates to
\begin{eqnarray*}
\hspace*{-6pt}&&\bigl\Vert E_0(X_{\tau_2}\!-\!X_{\tau_1}) \Vert E_0
(X_{\tau_2}\!-\!X_{\tau_1}, A_L) \Vert_2\!-\!E_0(X_{\tau_2}\!-\!X_{\tau_1})
\Vert E_0(X_{\tau_2}\!-\!X_{\tau_1}) \Vert_2 \\
\hspace*{-6pt}&&\quad\hspace*{0pt}{}\!+\!E_0(X_{\tau_2},\!-\!X_{\tau_1})
\Vert E_0(X_{\tau_2}\!-\!X_{\tau_1}) \Vert_2\!-\!E_0 (X_{\tau_2}\!-\!X_{\tau_1}, A_L)
\Vert E_0 (X_{\tau_2}\!-\!X_{\tau_1}) \Vert_2 \bigr\Vert_2\\
\hspace*{-6pt}&&\qquad\!\leq\!\Vert E_0(X_{\tau_2}\!-\!X_{\tau_1}) \Vert_2
\bigl\vert\Vert E_0(X_{\tau_2}\!-\!X_{\tau_1}, A_L) \Vert_2\!-\!\Vert
E_0(X_{\tau_2}\!-\!X_{\tau_1}) \Vert_2 \bigr\vert\\
\hspace*{-6pt}&&\qquad\quad{}\!+\!\Vert E_0(X_{\tau_2}\!-\!X_{\tau_1}, A_L^c) \Vert_2 \Vert
E_0(X_{\tau_2}\!-\!X_{\tau_1}) \Vert_2\\
\hspace*{-6pt}&&\qquad\!\leq\!2 \Vert E_0( X_{\tau_2}\!-\!X_{\tau_1}) \Vert_2 \Vert E_0
(X_{\tau_2}\!-\!X_{\tau_1}, A_L^c) \Vert_2,
\end{eqnarray*}
where the last inequality follows from the reverse triangle inequality.
But Cauchy--Schwarz's inequality and Lemma \ref{lem:ANEst} yield
\begin{eqnarray*}
\Vert E_0 (X_{\tau_2} - X_{\tau_1}, A_L^c) \Vert_2 &\leq&
E_0( \Vert X_{\tau_2} - X_{\tau_1} \Vert_2, A_L^c) \\[-2pt]
&\leq& E_0( \Vert X_{\tau_2} - X_{\tau_1} \Vert_2^2)^{1/2}
P_0(A_L^c)^{1/2}\\[-2pt]
&\leq& Ce^{-C^{-1}R_2(L)^\gamma/2},
\end{eqnarray*}
whence (\ref{vAvDiff}) is contained in $\mathcal{S}(\N)$ as a
function in $L$.
\end{pf}

Therefore,\vspace*{1pt}
$( \tilde{\pi}_{\hat{v}^\bot} (\sum_{j=1}^n (X_{\tau_j} -
X_{\tau_{j-1}})) )_{n \in\{ 2, \ldots, 2L^2 \} }$
is nearly a mean-zero martingale with respect to $P_0( \cdot \vert A_L)$
and this is what we will exploit in the proof of the next lemma.
%
\begin{lemma} \label{lem:annealedPosExitProbEst}
For $L$ and $x \in\tilde{\mathcal{P}}(0,L)$, define the event
\begin{eqnarray*}
&&F_{x,L}:= \bigl\{ \exists n \in\{0, \ldots, T_{L^2}\}
\dvtx\Vert\tilde{\pi}_{\hat{v}^\bot} (X_n - x) \Vert_\infty\geq
R_3(L)L \\[-2pt]
&&\hspace*{128.3pt}\mbox{ or } (X_n -x)\cdot e_1 < -R_2(L) \bigr\}.
\end{eqnarray*}
Then there exists a constant $C > 0$ such that for all $L$,
\[
\max_{x \in\tilde{\mathcal{P}}(0,L)}
P_x ( F_{x,L}) \leq C e^{-C^{-1}R_2(L)^\gamma}.
\]
In particular,
%
%
\begin{equation} \label{eq:supAnnealedNotFrontExitEst}
\max_{x \in\tilde{{\mathcal{P}}}(0,L)} P_x\bigl(X_{T_{\partial{\mathcal
{P}}(0,L)}} \notin\partial_+ {\mathcal{P}}(0,L)\bigr)
\leq C e^{-C^{-1}R_2(L)^\gamma}.
\end{equation}
\end{lemma}
\begin{pf}
Setting
$
F_{x,L}' := \{ \exists n \in\{0, \ldots, T_{L^2}\}
\dvtx\Vert\tilde{\pi}_{\hat{v}^\bot} (X_n - x) \Vert_\infty\geq
R_3(L)L \},
$
we have
%
%
\begin{equation}\label{eq:FxLSplit}
P_x ( F_{x,L}) \leq P_x ( F_{x,L}', A_L ) + P_x(A_L^c).
\end{equation}
Note that
$( \tilde{\pi}_{\hat{v}_L^\bot} (\sum_{j=1}^n (X_{\tau_j} -
X_{\tau_{j-1}})) )_{n \in\{2, \ldots, 2L^2 \}}$
is a $(d-1)$-dimensional martingale with respect to $P_x( \cdot
\vert A_L)$.
Furthermore, observe that due to Lem\-ma~\ref{lem:directionsDistance},
in particular we have
\[
\sup_{y \in{\mathcal{P}}(0,L)} \Vert\tilde{\pi}_{\hat{v}_L^\bot
} (y) - \tilde{\pi}_{\hat{v}^\bot} (y) \Vert_\infty\leq R_2(L)
\]
for $L$ large enough.
Therefore, Azuma's inequality applied to the coordinates yields
\begin{eqnarray*}
&&P_x ( F_{x,L},A_L )\\[-2pt]
&&\qquad\leq P_x \bigl( \exists n \in\{\tau_{1}, \ldots, \tau_{2L^2}\}
\dvtx\Vert\tilde{\pi}_{\hat{v}^\bot} (X_n - X_{\tau_1})\Vert_\infty
\geq R_3(L)L -2R_2(L) \vert A_L \bigr) \\[-2pt]
&&\qquad\leq P_x \bigl( \exists n \in\{\tau_{1}, \ldots, \tau_{2L^2}\}
\dvtx\Vert\tilde{\pi}_{\hat{v}_L^\bot} (X_n - X_{\tau_1})\Vert
_\infty\geq R_3(L)L - 3R_2(L) \vert A_L \bigr) \\[-2pt]
&&\qquad\leq P_x \bigl( \exists n \in\{ 1, \ldots, 2L^2\}
\dvtx\Vert\tilde{\pi}_{\hat{v}_L^\bot} (X_{\tau_n} - X_{\tau
_1})\Vert_\infty\geq R_3(L)L - 4R_2(L) \vert A_L \bigr) \\[-2pt]
&&\qquad\leq2 (d-1) \sum_{j = 1}^{2L^2} \exp\biggl\{ -\frac
{(R_3(L)L/2)^2}{2j R_2(L)^2} \biggr\} \\[-2pt]
&&\qquad\leq4(d-1) L^2
\exp\biggl\{ -\frac{(R_3(L)L/2)^2}{4L^2 R_2(L)^2}\biggr\}\\[-2pt]
&&\qquad\leq\exp\{ -R_3(L) \}
\end{eqnarray*}
for $L$ large enough.
In particular, in combination with (\ref{eq:FxLSplit}) and (\ref{eq:ANEst})
this reasoning finishes the proof of the first part. Equality
(\ref{eq:supAnnealedNotFrontExitEst}) is an immediate
consequence.\vspace*{-3pt}
\end{pf}

The following lemma, which
we will
prove in Section \ref{sec:furtherAuxRes} (see page \pageref
{pr:lem:lbound}),
provides lower bounds on certain exit probabilities.\vspace*{-3pt}
%
\begin{lemma} \label{lem:lbound}
Let $C'$ be a positive constant. Then there exists a positive constant
$c$ such that
for all $L$ large enough, and all $x \in\tilde{\mathcal{P}}(0,L)$,
$y \in\partial_+ \mathcal{P}(0,L)$ with
$
\Vert\tilde{\pi}_{\hat{v}^\bot} (y - x) \Vert_1 < C' L,
$
we have
\[
P_x (X_{T_{\partial\mathcal{P}(0,L)}}= y )\geq cL^{1-d}.\vspace*{-3pt}
\]
\end{lemma}

Let now $x \in\mathcal{L}_{L}$ and
$z \in\tilde{{\mathcal{P}}}(x,L)$. Then, following \cite{Be-09},
for $\omega\in\Omega$
we define~$\mu_{z,x,\omega}^L$ to be the distribution
of $X_{T_{\partial{\mathcal{P}}(x,L)}}$ with respect to
$
P_{z,\omega}( \cdot \vert T_{\partial{\mathcal{P}}(x,L)} =
T_{\partial_+ {\mathcal{P}}(x,L)}).
$
Similarly, we define
$\mu_{z,x}^L$ to be the distribution of
$
X_{T_{\partial{\mathcal{P}}(x,L)}}
$
with respect to
$
P_z( \cdot \vert T_{\partial{\mathcal{P}}(x,L)} =
T_{\partial_+ {\mathcal{P}}(x,L)}).
$

We now get the following bounds for $\Var_{\mu_{x,0}^L}$, which will
turn out to be useful in the proof of Corollary \ref{cor:closeness} below.\vspace*{-3pt}
%
\begin{lemma} \label{lem:VarBd}
There exists a constant $C$ such that for all $x \in\tilde{\mathcal
{P}}(0,L)$ and all~$L$,
\[
C^{-1} L^2 \leq\Var_{\mu_{x,0}^L} \leq C L^2.\vspace*{-3pt}
\]
\end{lemma}

\begin{pf}
The lower bound is a consequence of Lemma \ref{lem:lbound}.\vadjust{\goodbreak}

To prove the
upper bound, note that $S_n := \sum_{k=1}^n X_{\tau_k} - X_{\tau_{k-1}}
- E_0 (X_{\tau_k} - X_{\tau_{k-1}})$ is a martingale in $n$ with
respect to $P_0$.
We define the stopping time
\[
T:= \inf\Biggl\{ n \in\N\dvtx\Biggl( S_n + \sum_{k=1}^n E_0 (X_{\tau_k}
- X_{\tau_{k-1}}) \Biggr) \cdot e_1 \geq L^2 \Biggr\}
\]
and note that in particular $(S_{n \wedge T} \cdot e_j)_{n \in\N}$ is
a martingale for any $j \in\{2, \ldots, d\}$.
The independence of the increments yields that so is
\[
\bigl( (S_{n \wedge T} \cdot e_j)^2 - \bigl( E(S_m \cdot e_j)^2
\bigr)_{m = n \wedge T} \bigr)_{n \in\N}.
\]
Since for $n=0$ the martingale equals $0$, we have, noting that
\[
E(S_m \cdot e_j)^2 = \sum_{k=1}^m E\bigl( \bigl(X_{\tau_k} - X_{\tau
_{k-1}} - E_0 (X_{\tau_k} - X_{\tau_{k-1}})\bigr) \cdot e_j
\bigr)^2\vadjust{\goodbreak}
\]
as well as
$T \leq L^2$,
that
\[
E (S_T \cdot e_j)^2 \leq C L^2.
\]
Taking into consideration Lemma \ref{lem:ANEst} and Lemma \ref
{lem:annealedPosExitProbEst}, this implies the upper bound.\vspace*{-3pt}
\end{pf}

For $x \in\Z^d$ and $k \in\Z$, we will use the
%
%
\begin{equation} \label{eq:hyperPlane}
H_k:=\{x \in\Z^d\dvtx x \cdot e_1 =k\}
\end{equation}
from the following proof onward.

In \cite{Be-09}, the author derived a result similar to the following
corollary of
Proposition \ref{prop:closenessbase}.\vspace*{-3pt}
%
\begin{corollary} \label{cor:closeness}
Fix $\vartheta\in(0,5/8]$ and let $L$ be large enough.
Furthermore, let
$k \in\{1, \ldots, \iota\}$,
$x \in C_L \cap\mathcal{L}_{L_k}$
and $\omega\in\Omega$
such that (\ref{eq:quenAnnExpClose})
and (\ref{eq:quenAnnProbClose}) are fulfilled for
this choice of $\vartheta$ and all $z \in\tilde{\mathcal{P}}(x,L_k)$.

Then
$\mu_{z,x,\omega}^{L_k}$ is
$( L_k^{-\vartheta({d-1})/({2(d+1)})}, 2dL_k^\vartheta)$-close
to
$\mu_{z,x}^{L_k}$.\vspace*{-3pt}
\end{corollary}
\begin{pf}
For fixed $k,x, \omega$ and $z$ as in the assumptions, we will show
the desired closeness for~$L$
large enough. Observing that this lower bound on~$L$ holds uniformly in
the admissible choices
of $k, x, \omega$ and $z$ then finishes the proof.

We will construct the coupling of Definition
\ref{def:closeness} in the case $x=0$, the remaining cases being
handled in exactly the same manner.
Cover $\partial_+{\mathcal{P}}(0,L)$
by
at most $n=\lceil2R_6(L_k)L_k^{1-\vartheta} \rceil^{d-1}$
disjoint cubes
$Q_1,Q_2,\ldots, Q_n$ of side length $\lceil L_k^\vartheta\rceil$.
Consider an i.i.d. sequence $(Y_j)_{j \in\N}$ of random variables
defined on a probability space with probability measure $P^{{*}}$
(the space should also be large enough to accommodate the random
variables we will define in the remaining part of this proof)
such that
\[
P^{*} (Y_j =x)=\frac{\mu_{z,0}^{L_k} (\{x\})}{\mu_{z,0}^{L_k}
(Q_j)},\qquad x \in Q_j,
\]
and $P^{*} (Y_j =x)=0$ if $x \notin Q_j$; set
\[
Y:=\sum_{j=1}^n Y_j \mathbh{1}_{\{X_{T_{\partial P(0,L_k)}}\in Q_j\}}
\]
and
\[
\mathbf{P}_{z,\omega} := P_{z,\omega}\bigl( \cdot \vert T_{\partial
_+ \mathcal{P}(0,L_k)} = T_{\partial\mathcal{P}(0,L_k)}\bigr) \otimes
P^{*}.
\]
Clearly, $\mathbf{P}_{z,\omega}$-a.s., $\Vert X_{T_{\partial P(0,L_k)}}-Y
\Vert_1 \leq(d-1) \lceil L_k^\vartheta\rceil$ and consequently we have
\[
\Vert\mathbf{E}_{z,\omega} Y- \mathbf{E}_{z,\omega} X_{T_{\partial
P(0,L_k)}} \Vert_1 < (d-1) \lceil L_k^\vartheta\rceil.
\]
Display (\ref{eq:quenAnnExpClose}) yields
\[
\bigl\Vert\mathbf{E}_{z,\omega} X_{T_{\partial P(0,L_k)}}- E_{z}
\bigl(X_{T_{\partial P(0,L_k)}} \vert
T_{\partial_+ \mathcal{P}(0,L_k)} = T_{\partial\mathcal{P}(0,L_k)}\bigr)
\bigr\Vert_1 \leq R_4({L_k})\vadjust{\goodbreak}
\]
and
thus
\[
\bigl\Vert\mathbf{E}_{z,\omega} Y - E_{z}
\bigl(X_{T_{\partial P(0,L_k)}} \vert T_{\partial_+ \mathcal
{P}(0,L_k)} = T_{\partial\mathcal{P}(0,L_k)}\bigr) \bigr\Vert_1 <d
L_k^\vartheta
\]
for $L$ large enough.
Let now $U$ be an $H_0$-valued random variable
defined~on the same probability space as the sequence $(Y_j)_{j \in\N
}$ (and choose $U$ to be~inde\-pendent of everything else) such that
$\mathbf{P}_{z,\omega}$-a.s. we have $\Vert U \Vert_1 \leq
dL_k^\vartheta$ as well~as
\[
\mathbf{E}_{z,\omega} U=E_{z} \bigl(X_{T_{\partial P(0,L_k)}} \vert
T_{\partial_+ \mathcal{P}(0,L_k)} = T_{\partial\mathcal{P}(0,L_k)}\bigr)
- \mathbf{E}_{z,\omega} Y.
\]
Define
\[
Z_0: = Y + U
\]
and
\[
Z_2: = X_{T_{\partial P(0,L_k)}}.
\]
Then taking $\mathbf{P}_{z,\omega}$ as the $\mu$ of Definition \ref
{def:closeness}, part \hyperlink{item:C3}{(c)} of that definition
is fulfilled for $K=2dL_k^\vartheta$ and $L$ large enough.
To show the remaining parts, we first note that for $y \in\Z^d$ we have
\[
\mathbf{P}_{z,\omega} (Z_0 = y)=\sum_{u \dvtx\Vert u \Vert_1 \leq
dL_k^\vartheta} \mathbf{P}_{z,\omega}(U=u) \mathbf{P}_{z,\omega} (Y=y-u).
\]
Since furthermore $\mathbf{P}_{z,\omega} \circ Y^{-1}$ is supported on
$\partial_+ {\mathcal{P}}(0,{L_k})$, we get
%
%
\begin{eqnarray} \label{eq:fmqnc}\qquad
&&\sum_{y \in\Z^d} \vert\mathbf{P}_{z,\omega} (Z_0 = y)-\mu
_{z,0}^{L_k} (y) \vert\nonumber\\
&&\qquad\leq2\sum_{u \dvtx\Vert u \Vert_1 \leq d L_k^\vartheta} \mathbf
{P}_{z,\omega}(U=u)\\
&&\hspace*{86.3pt}{}\times\sum_{y\dvtx y -u \in\partial_+{\mathcal{P}}(0,{L_k})} \vert\mathbf
{P}_{z,\omega}(Y=y-u)-\mu_{z,0}^{L_k}(y) \vert.
\nonumber
\end{eqnarray}
By heat-kernel estimates to be proven later [cf. part \hyperlink
{item:firstDer}{(b)} of Lemma \ref{lem:annExitDistDerivatives}],
for each $y \in H_{L_k^2}$ and every $u$ such that
$\Vert u \Vert_1 \leq dL_k^\vartheta$,
%
\[
\vert\mu_{z,0}^{L_k} (y-u) - \mu_{z,0}^{L_k} (y) \vert\leq
CL_k^\vartheta\cdot L_k^{-d}=CL_k^{\vartheta-d}.
\]
In combination with (\ref{eq:fmqnc}) and the validity of (\ref
{eq:quenAnnProbClose}), this
yields
%
%
\begin{eqnarray}
\label{eq:C2Est}
\hspace*{14pt}&&\sum_{y \in\Z^d} \vert\mathbf{P}_{z,\omega} (Z_0 = y)- \mu
_{z,0}^{L_k}(y) \vert\nonumber\\
\hspace*{14pt}&&\qquad\leq2\sum_{y \in\partial_+{\mathcal{P}}(0,{L_k})} \bigl( \vert
\mathbf{P}_{z,\omega}(Y = y) - \mu_{z,0}^{L_k} (y) \vert+C
L_k^{\vartheta-d}\bigr)\nonumber\\
\hspace*{14pt}&&\qquad\leq2\Biggl[ \Biggl(\sum_{j=1}^n \vert\mu_{z,0,\omega}^{L_k}
(Q_j)-\mu_{z,0}^{L_k}(Q_j) \vert
\Biggr)+(2R_6({L_k}))^{d-1}L_k^{d-1}\cdot CL_k^{\vartheta-d} \Biggr]\\
\hspace*{14pt}&&\qquad\leq CR_6({L_k})^{d-1}L_k^{\vartheta-1} + \lceil
2R_6(L_k)L_k^{1-\vartheta} \rceil^{d-1}
\cdot
L_k^{(\vartheta-1)(d-1)-\vartheta({d-1})/({d+1})}
\nonumber\\
\hspace*{14pt}&&\qquad\leq CR_6({L_k})^{d-1} \bigl(L_k^{\vartheta-1}+L_k^{-\vartheta
({d-1})/({d+1})} \bigr)
\leq R_7({L_k})L_k^{-\vartheta({d-1})/({d+1})}\nonumber\\
\hspace*{14pt}&&\qquad<L_k^{-\vartheta({d-1})/({2(d+1)})},\nonumber
\end{eqnarray}
$L$ large enough;
here, the second inequality is obtained by noting that the sign of
$
\vert\mathbf{P}_{z,\omega}(Y = y)-\mu_{0,z}^{L_k}(y) \vert
$
is constant as $y$ varies over $Q_j$ for fixed $j$,
while the penultimate inequality takes advantage of $\vartheta\leq5/8$
and $d \geq4$.
Thus, due to (\ref{eq:C2Est}), there exists a random variable $Z_1$
defined on the probability space with
probability measure $\mathbf{P}_{z,\omega}$ such that $\mathbf
{P}_{z,\omega
} \circ Z_1^{-1} = \mu_{z,0}$ and
$\mathbf{P}_{z,\omega}(Z_1 \ne Z_0) \leq L_k^{-\vartheta ({d-1})/({2(d+1)})}$.
This establishes \hyperlink{item:C1}{(a)}, \hyperlink{item:C2}{(b)}
and \hyperlink{item:C4}{(d)} of Definition \ref{def:closeness}
for $\lambda= L_k^{-\vartheta({d-1})/({2(d+1)})}$.

To see \hyperlink{item:C5}{(e)},
observe that
\[
\mathbf{Var}_{x,\omega} Z_1 = \Var_x \bigl(X_{T_{\partial\mathcal{P}(0,L)}}
\vert T_{\partial\mathcal{P}(0,L)} = T_{\partial_+ \mathcal{P}(0,L)}\bigr).
\]
Now note that
the support of $\mu_{0,z}^{L_k}(\cdot) - \mathbf{P}_{z,\omega} (Z_0 =
\cdot)$ is contained in
\[
\{y \in H_{L_k^2} \dvtx\exists z \in\partial_+ {\mathcal
{P}}(0,{L_k}) \mbox{ such that } \Vert y-z \Vert_1 \leq
dL_k^\vartheta\}.
\]
Thus,\vspace*{1pt} for any $y$ in the support of $\mu_{z,0}^{L_k}(\cdot) - \mathbf
{P}_{z,\omega} (Z_0 = \cdot)$
we get as a consequence of \hyperlink{item:C2}{(b)} in combination
with the penultimate line of (\ref{eq:C2Est}) that
\begin{eqnarray*}
&&\sum_{x}
\bigl\Vert x - E_z \bigl(X_{T_{\partial P(0,L_k)}} \vert T_{\partial
_+ \mathcal{P}(0,L_k)} = T_{\partial\mathcal{P}(0,L_k)}\bigr) \bigr\Vert
_1^2\\
&&\quad\hspace*{0pt}{} \times\vert\mathbf{P}_{z,\omega} (Z_1 =x) - \mathbf{P}_{z,\omega}
(Z_0 = x) \vert\\
&&\qquad\leq4 d^2({L_k}R_6^2({L_k}))^2\sum_{x} \vert \mu
_{z,0}^{L_k}(x) - \mathbf{P}_{z,\omega} (Z_0 = x) \vert\\
&&\qquad\leq4d^2R_6^2({L_k})L_k^2\cdot R_7({L_k})L_k^{-\vartheta
({d-1})/({d+1}) },
\end{eqnarray*}
where the last inequality holds for $L$ large enough.
In combination with Lem\-ma~\ref{lem:VarBd}, we deduce that the
right-hand side is bounded from above
by $\lambda\Var Z_1$ for $L$ large enough, which finishes the proof.
\end{pf}

\subsection{Auxiliary walk} \label{subsec:AuxWalk}

As a preparation to prove Proposition \ref{prop:quenchedConeTypeExit},
for each environment,
we introduce a refinement $(Y_n)$ of the finite-time auxiliary
random walk defined in \cite{Be-09}. In blocks $\mathcal{P}(x,L_k)$
where the environment is such that
the quenched RWRE $(X_n)$ behaves similarly to the annealed
one, the quenched walk $(Y_n)$ will behave quite like $(X_n)$. In
blocks where the quenched and annealed behavior
of $(X_n)$ differ significantly, the quenched walk $(Y_n)$ will make up
for this deviation by
corrections, in order to more or less mimick the annealed behavior of $(X_n)$.
As a consequence, the quenched walk $(Y_n)$ starting in $0$ will leave
$C_L$ through $\partial_+ C_L$ with a
probability not too small, with respect to sufficiently many environments.
Note that its construction will depend on a couple of parameters 
and in particular will be done for
each $L > 0$ separately. For
the sake of notational simplicity, we do not explicitly name these
dependencies in the notation $(Y_n)$.
In order to facilitate understanding for the reader familiar with \cite
{Be-09}, we
stick to the notation of that paper wherever appropriate.

On a heuristic level, the construction of the auxiliary walk $(Y_n)$
can be described as follows.
Let $L$ and $\omega$ be given. In order to leave $C_L$
through~$\partial_+ C_L$,
the walker starts with performing a few deterministic
steps in positive $e_1$-direction.

Then, starting a recursive step,
there is associated a natural scale $k' \in\{1, \ldots, \iota\}$ to
the current position of the walker
(this scale is roughly given by the largest $k \in\{1, \ldots, \iota
\}$ for which $L_k^2$ divides
the current $e_1$-coordinate of the walker);
the walker then looks for good boxes of the form $\mathcal{P}(x,L_k)$,
such that $k \in\{1, \ldots, k'\}$,
$x \in\mathcal{L}_{L_k}$ and such that his current position is
contained in $\tilde{\mathcal{P}}(x,L_k)$.
We now distinguish cases:
\begin{itemize}
\item
If such a box exists,
then the walker picks the largest of these boxes and moves according to
a random walk in the corresponding environment, conditioned
on leaving this box through its right boundary part. If this box is of
the form
$\mathcal{P}(x,L_k)$ for some $k < k'$, then before starting the
recursion step from a position with
natural scale $k'$ again,
the walker will perform a correction, making up for having moved in boxes
smaller than the ones corresponding to its natural scale.

\item
If no such good box exists, the walker performs some deterministic
steps in positive $e_1$-direction again and
then returns to the start of
the recursive step.
\end{itemize}

To formally construct our process, we need some auxiliary results.
The following lemma will be proved in Section \ref
{sec:furtherAuxRes} (see page \pageref{pr:bddFluctExp}).
%
\begin{lemma} \label{lem:bddFluctExp}
There exists a finite constant $C$ such that for all $L$ and $x \in
\tilde{\mathcal{P}}(0,L)$,
\[
\biggl\Vert E_{\mu_{x,0}^L} -x - \frac{L^2 -x\cdot e_1}{\hat{v} \cdot
e_1} \hat{v} \biggr\Vert_1 \leq C R_2(L)
\]
and
%
%
\begin{equation} \label{eq:bddFluctExp}
\biggl\Vert E_x X_{T_{\partial{\mathcal{P}}(0,L)}} -x - \frac{L^2
-x\cdot e_1}{\hat{v} \cdot e_1} \hat{v} \biggr\Vert_1 \leq C R_2(L).
\end{equation}
\end{lemma}

In order to state further auxiliary results, for $x \in\Z^d$ such
that $x \cdot e_1 \in L_k^2 \N$, define
$z(x,k)$ to be an element $z \in\mathcal{L}_{L_k}$
such that $x\cdot e_1 = z \cdot e_1$ and $x \in\tilde{{\mathcal
{P}}}(z, L_k)$.
Furthermore, for $x$ such that $x \cdot e_1 \notin L_k^2 \N$ set
$z(x,k) := 0$.
In addition, abbreviate for $j,k \in\N$ the hitting times
\[
T_{k}(j):= \inf\{n \in\N\dvtx Y_n \cdot e_1 = jL_k^2\}.
\]

\begin{lemma} \label{lem:goodExitOfBigBox}
Let $k \in\{1, \ldots, \iota-1\}$, $\Delta_0 \in H_0 \cap\tilde
{\mathcal{P}}(0,L_{k+1})$ deterministic
and $(\Delta_i)_{i\in\{1, \ldots, {\lfloor L^\chi\rfloor}^2\}}$
be random variables.
Set $S_j:=\sum_{i=0}^j \Delta_i$ and
assume furthermore
that for every $i$,
conditioned on $\Delta_1, \ldots, \Delta_{i-1}$, the variable
$\Delta_i$ takes values in $\partial_+ \mathcal{P}(z(S_{i-1},k),L_k)
- z(S_{i-1},k)$ only, with
%
%
\begin{equation} \label{eq:expCloseAss}
\quad\Vert E(\Delta_i \vert \Delta_1, \ldots, \Delta_{i-1}) -
( E_{\mu_{S_{i-1},z(S_{i-1}, k)}^{L_k}} - S_{i-1} )
\Vert_1 \leq R_4(L_k)
\end{equation}
a.s.
Then for $L$ large enough and $t \geq R_5(L_{k})L_{k+1}$,
%
%
\begin{eqnarray}
\label{eq:goodExitOfBigBoxEq}
&&P \bigl(\exists j \in \{1, \ldots, {\lfloor L^\chi\rfloor}^2\} \dvtx
\Vert\tilde{\pi}_{\hat{v}^\bot} ( S_j - \Delta_0) \Vert_\infty
\geq t\bigr)\nonumber\\[-8pt]\\[-8pt]
&&\qquad \leq2(d-1) {\lfloor L^\chi\rfloor}^2 \exp\biggl\{ -\frac
{t^2}{72L_{k+1}^2R_6(L_k)^2} \biggr\}.
\nonumber
\end{eqnarray}
\end{lemma}
\begin{pf}
Noting that $\tilde{\pi}_{\hat{v}^\bot} (\lambda\hat{v}) = 0$ for
all $\lambda\in\R$,
the triangle inequality yields
%
%
\begin{equation} \label{eq:splitIntoZs}
\Vert\tilde{\pi}_{\hat{v}^\bot}
(S_j - \Delta_0)
\Vert_\infty
\leq\bigr\Vert Z_j^{(1)} \bigl\Vert_\infty
+ \bigl\Vert Z_j^{(2)} \bigr\Vert_\infty
+ \bigl\Vert Z_j^{(3)} \bigr\Vert_\infty,
\end{equation}
where
\begin{eqnarray*}
Z_j^{(1)} &:=& \sum_{i=1}^{j} \tilde{\pi}_{\hat{v}^\bot} \bigl(
\Delta_i
- E(\Delta_i \vert \Delta_1, \ldots, \Delta_{i-1}) \bigr),\qquad
j \in\{1, \ldots, {\lfloor L^\chi\rfloor}^2\},\\
Z_0^{(1)} &:=& 0,
\\
Z_j^{(2)}&:=&
\sum_{i=1}^j \tilde{\pi}_{\hat{v}^\bot} \bigl( E(\Delta_i
\vert \Delta_1, \ldots, \Delta_{i-1})
- ( E_{\mu^{L_k}_{S_{i-1},z(S_{i-1},k)}} - S_{i-1} )\bigr)
\end{eqnarray*}
and
\[
Z_j^{(3)}:=
\sum_{i=1}^j \tilde{\pi}_{\hat{v}^\bot} \biggl( ( E_{\mu
^{L_k}_{S_{i-1},z(S_{i-1},k)}} - S_{i-1} )
- \frac{L_k^2}{\hat{v} \cdot e_1} \hat{v}\biggr).
\]
Due to (\ref{eq:expCloseAss}), a.s.
%
%
\begin{equation}\label{eq:Z2Bd}
\bigl\Vert Z_j^{(2)} \bigr\Vert_\infty\leq jR_4(L_k),
\end{equation}
while Lemma \ref{lem:bddFluctExp} results in
%
%
\begin{equation}\label{eq:Z3Bd}
\bigl\Vert Z_j^{(3)} \bigr\Vert_\infty\leq CjR_2(L_k).
\end{equation}
Using (\ref{eq:splitIntoZs}) to (\ref{eq:Z3Bd})
and because of $t \geq R_5(L_{k})L_{k+1}$, for $L$ large enough
the probability in (\ref{eq:goodExitOfBigBoxEq})
can be bounded from above by
\[
P \bigl( \exists j \in\{1, \ldots, {\lfloor L^\chi\rfloor}^2\} \dvtx
\bigl\Vert
\tilde{\pi}_{\hat{v}^\bot} \bigl(Z_j^{(1)}\bigr) \bigr\Vert_\infty
\geq t /3 \bigr).
\]
Now
with respect to
$P$,
the sequence $( \tilde{\pi}_{\hat{v}^\bot} (Z_j^{(1)})
)_{j \in\{0, \ldots, {\lfloor L^\chi\rfloor}^2\}}$
is a $(d-1)$-dimensional mean zero martingale
such that $\Vert\tilde{\pi}_{\hat{v}^\bot} (Z_{j+1}^{(1)})-\tilde
{\pi}_{\hat{v}^\bot} (Z_j^{(1)}) \Vert_\infty\leq2L_kR_6(L_k)$
for all $j \in\{0, \ldots, {\lfloor L^\chi\rfloor}^2\}$. Thus, Azuma's
inequality yields
\begin{eqnarray*}
&&P\bigl( \exists j \in \{1, \ldots, {\lfloor L^\chi\rfloor}^2 \} \dvtx
\bigl\Vert
\tilde{\pi}_{\hat{v}^\bot} \bigl(Z_j^{(1)}\bigr) \bigr\Vert_\infty\geq t/3
\bigr) \\
&&\qquad\leq2(d-1){\lfloor L^\chi\rfloor}^2 \exp\biggl\{ -\frac{t^2}{72
{\lfloor L^\chi\rfloor}^2 (L_k R_6(L_k))^2} \biggr\}\\
&&\qquad= 2(d-1) {\lfloor L^\chi\rfloor}^2 \exp\biggl\{ -\frac{t^2}{72 L_{k+1}^2
R_6(L_k)^2} \biggr\}.
\end{eqnarray*}
\upqed\end{pf}

\begin{figure}

\includegraphics{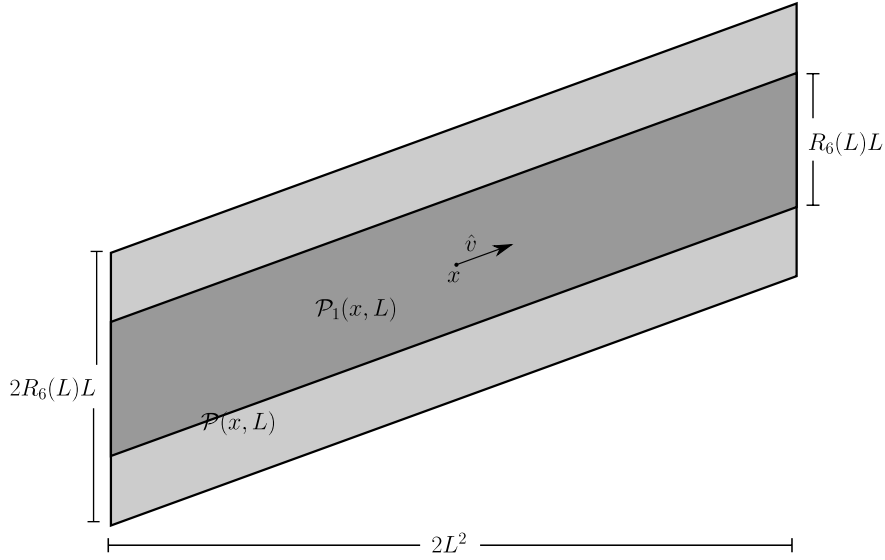}

\caption{The set $\mathcal{P}_1(x,L)$ contained in $\mathcal{P}(x,L)$.}
\label{fig5}
\end{figure}

We now introduce some quantities that will play an important role in
the remaining part of this paper.
For $k \geq2$, let
\[
\lambda_k := R_{9+k}(L) L_{1}^{-\vartheta({d-1})/({2(d+1)})}
\]
%
and
\[
K_k :=4 {\lfloor L^\chi\rfloor}^2 d \lceil L_{k-1}^\vartheta\rceil,
\]
where $\vartheta:= \chi$ as in the definition of good blocks.
Furthermore, define the boxes
\[
\mathcal{P}_1(0,L) := \{y \in\mathcal{P}(0,L) \dvtx\Vert\tilde
{\pi}_{\hat{v}^\bot} (y) \Vert_\infty\leq R_6(L)L/2 \}
\]
and
\[
\mathcal{P}_1(x,L) := x + \mathcal{P}_1(0,L)
\]
as well as its right boundary part
\[
\partial_+ \mathcal{P}_1(x,L) := \{ y \in\partial\mathcal
{P}_1 (x,L) \dvtx(y -x)\cdot e_1 = L^2 \};
\]
cf. Figure \ref{fig5}.

%
%

From now on, we will occasionally emphasize the process to which a certain
random time refers by writing it as a superscript to the corresponding
random time (as, e.g., $T^S_{\partial_+ \mathcal{P}_1(0,L_{k+1})}$ in the
following lemma).
%
\begin{lemma} \label{lem:condCloseness}
$\!\!\!$Consider $S_{{\lfloor L^\chi\rfloor}^2}$ of Lemma \ref{lem:goodExitOfBigBox}
and assume that the distri\-bution $P\!\circ\!S_{{\lfloor L^\chi\rfloor}^2}^{-1}$
of $S_{{\lfloor L^\chi\rfloor}^2}$ with respect to $P$ is
$(2k \lambda_{k+1}, 2k K_{k+1})$-close~to~$\mu_{x,0}^{L_{k+1}}$
for some $x \in\tilde{\mathcal{P}}(0,L_{k+1})$.
Then for $L$ large enough, the distribution of $S_{{\lfloor L^\chi
\rfloor}^2}$
with respect to
$P ( \cdot \vert T^S_{\partial_+ \mathcal
{P}_1(0,L_{k+1})} = T^S_{\partial\mathcal{P}_1(0,L_{k+1})} )$
is $((2k+1) \lambda_{k+1}, (2k+\break 1) K_{k+1})$-close to $\mu_{x,0}^{L_{k+1}}$
for all admissible choices of $\Delta_0$, $k$ and~$x$.
\end{lemma}
\begin{pf}
Let $\mu$ be the coupling of $Z_0,Z_1$ and $Z_2$ as in the definition of
$
(2k \lambda_{k+1}, 2kK_{k+1})$-closeness,
that is, such that:
\begin{longlist}[(a)]
\item[(a)] $\mu\circ Z_1^{-1} = \mu_{x,0}^{L_{k+1}}$ and
$\mu\circ Z_2^{-1} = P \circ S_{{\lfloor L^\chi\rfloor}^2}^{-1}$,

\item[(b)] $\mu(Z_1 \not= Z_0) \leq2k \lambda_{k+1}$,

\item[(c)] $\mu(\Vert Z_0 - Z_2 \Vert_1 \leq2k K_{k+1}) = 1$,

\item[(d)] $E_\mu Z_1 = E_\mu Z_0$,

\item[(e)] $\sum_x \Vert x - E_\mu Z_1 \Vert_1^2 \cdot\vert\mu
(Z_1 = x) - \mu(Z_0 = x) \vert\leq2k \lambda_{k+1} \Var Z_1$.
\end{longlist}

Then we look for $Z_0', Z_1'$ and $Z_2'$ such that:
\begin{longlist}[(a$'$)]
\item[(a$'$)] \hypertarget{item:D1'} $\mu\,{\circ}\,{Z'_1}^{-1}\,{=}\,\mu
_{x,0}^{L_{k+1}}$ and
$\mu\,{\circ}\,{Z'_2}^{-1}\,{=}\,P ( \cdot \vert T^S_{\partial_+
\mathcal{P}_1(0,L_{k+1})}\,{=}\,T^S_{\partial\mathcal{P}_1(0,L_{k+1})}
)\,{\circ}\,S_{{\lfloor L^\chi\rfloor}^2}^{-1}$,

\item[(b$'$)]\hypertarget{item:D2'} $\mu(Z'_1 \not= Z'_0) \leq
(2k+1) \lambda_{k+1}$,

\item[(c$'$)] \hypertarget{item:D3'} $\mu(\Vert Z'_0 - Z'_2 \Vert_1
\leq(2k+1) K_{k+1}) = 1$,

\item[(d$'$)] \hypertarget{item:D4'} $E_\mu Z'_1 = E_\mu Z'_0$,

\item[(e$'$)]\hypertarget{item:D5'} $\sum_x \Vert x - E_\mu Z'_1
\Vert_1^2 \cdot\vert\mu(Z'_1 = x) - \mu(Z'_0 = x) \vert\leq
(2k+1) \lambda_{k+1} \Var Z'_1$.
\end{longlist}

For this purpose and due to Remark \ref{rem:couplingSpace}, we can assume
\[
Z_2 = S_{{\lfloor L^\chi\rfloor}^2}
\]
without loss of generality.
Set
\[
Z_1' := Z_1
\]
and
\[
Z_2' := Z_2 \mathbh{1}_{T^S_{\partial_+ \mathcal{P}_1(0,L_{k+1})} =
T^S_{\partial\mathcal{P}_1(0,L_{k+1})}}
+ Z_2^* \mathbh{1}_{T^S_{\partial_+ \mathcal{P}_1(0,L_{k+1})} \not=
T^S_{\partial\mathcal{P}_1(0,L_{k+1})}},
\]
where $ Z_2^* $ is independent of the remaining random variables and
distributed as $S_{{\lfloor L^\chi\rfloor}^2}$ with respect to
\[
P \bigl( \cdot \vert T^S_{\partial_+ \mathcal
{P}_1(0,L_{k+1})} = T^S_{\partial\mathcal{P}_1(0,L_{k+1})}\bigr).
\]
Furthermore, set
\[
Z_0^* := Z_0 \mathbh{1}_{T^S_{\partial_+ \mathcal{P}_1(0,L_{k+1})} =
T^S_{\partial\mathcal{P}_1(0,L_{k+1})}}
+ Z_2' \mathbh{1}_{T^S_{\partial_+ \mathcal{P}_1(0,L_{k+1})} \not=
T^S_{\partial\mathcal{P}_1(0,L_{k+1})}}.
\]
Now
as $EZ_1' = EZ_0$ and since due to Lemma
\ref{lem:goodExitOfBigBox} we have that
\[
\max_{\Delta_0} P \bigl( T^S_{\partial_+ \mathcal{P}_1(0,L_{k+1})}
\not= T^S_{\partial\mathcal{P}_1(0,L_{k+1})} \bigr)
\]
is contained in $\mathcal{S}(\N)$ as a function in $L$
(where the maximum is taken over all admissible choices of $\Delta_0$;
see the assumptions of
Lemma \ref{lem:goodExitOfBigBox}), it follows that
$\Vert EZ_0^* - EZ_1'\Vert_1$ is contained in $\mathcal{S}(\N)$ as a
function in $L$.
Thus, there exists a random variable $U$ taking values in $H_0$ such that
$
P( \Vert U \Vert_1 \leq K_{k+1})= 1,
$
$
P(U \not= 0)
$
is contained in $\mathcal{S}(\N)$ as a function of $L$, and
such that $EZ_0^* + EU = EZ_1'$. Set
\[
Z_0' := Z_0^* + U.
\]

Then
\hyperlink{item:D1'}{(a$'$)}, \hyperlink{item:D3'}{(c$'$)} and
\hyperlink{item:D4'}{(d$'$)} are fulfilled. Furthermore,
\begin{eqnarray*}
P(Z_0' \not= Z_1') &\leq& 2k \lambda_{k+1} + P\bigl(T^S_{\partial_+
\mathcal{P}_1(0,L_{k+1})} \not= T^S_{\partial\mathcal{P}_1(0,L_{k+1})}\bigr)
+ P(U \not= 0) \\
&\leq&(2k+1) \lambda_{k+1}
\end{eqnarray*}
for $L$ large enough, which establishes \hyperlink{item:D2'}{(b$'$)}.
With respect to the variance bound we obtain
\begin{eqnarray*}
&&\sum_{y}\Vert y - E_{\mu_{x,0}^{L_{k+1}}} \Vert_1^2 \cdot\vert\mu
_{x,0}^{L_{k+1}}(y) - P(Z_0'=y) \vert\\
&&\qquad =\sum_{y} \Vert y - E_{\mu_{x,0}^{L_{k+1}}} \Vert_1^2
\\
&&\qquad\quad\hspace*{13.6pt}{}\times\vert\mu_{x,0}^{L_{k+1}}(y)
- P(Z_0 \mathbh{1}_{T^S_{\partial_+ \mathcal{P}_1(0,L_{k+1})} =
T^S_{\partial\mathcal{P}_1(0,L_{k+1})}} \\
&&\qquad\quad\hspace*{94.3pt}{} + Z_2' \mathbh{1}_{T^S_{\partial_+ \mathcal{P}_1(0,L_{k+1})}
\not= T^S_{\partial\mathcal{P}_1(0,L_{k+1})}} + U = y) \vert\\
&&\qquad \leq(3dL_{k+1} R_6(L_{k+1}))^2 \bigl(P\bigl(T^S_{\partial_+ \mathcal
{P}_1(0,L_{k+1})} \not= T^S_{\partial\mathcal{P}_1(0,L_{k+1})}\bigr) +
P(U \not= 0) \bigr) \\
&&\qquad\quad{} + \sum_{y} \Vert y - E_{\mu_{x,0}^{L_{k+1}}} \Vert_1^2 \cdot
\vert\mu_{x,0}^{L_{k+1}}(y) - P(Z_0 = y) \vert\\
&&\qquad\leq(2k+1) \lambda_{k+1} \Var_{\mu_{x,0}^{L_{k+1}}}
\end{eqnarray*}
for $L$ large enough.
Since the above computations are uniform in the admissible choices of
$\Delta_0$, $k$ and $x$,
the result follows.
\end{pf}
%
\begin{lemma} \label{lem:transCloseness}
Let $k \in\{1, \ldots, \iota\}$ and $x \in\tilde{\mathcal
{P}}(0,L_k) \cap H_0$.
Furthermore,
let a distribution $\nu$ be given which is supported on $\partial_+
\mathcal{P}(0,L_k)$ and
$((2k-1)\lambda_k,(2k-1)K_k)$-close to $\mu_{x,0}^{L_k}$.
Then for $L$ large enough,
$
\nu(\cdot+ x)
$
is $(2k\lambda_k,\break2kK_k)$-close to $\mu_{0,0}^{L_k}$
for all admissible choices of $k$ and $x$.
\end{lemma}

\begin{pf}
If $\nu$ is $((2k-1)\lambda_k,(2k-1)K_k)$-close to $\mu
_{x,0}^{L_k}$, then there exist~$Z_0$, $Z_1$ and $Z_2$
fulfilling the requirements of Definition \ref{def:closeness}, where
we denote the coupling measure
by $P$.

We set $Z_2':= Z_2 -x$ and will construct $Z_0'$ and $Z_1'$ such that
the corresponding points
of Definition \ref{def:closeness} are satisfied.
First of all, note that (as a~consequence of Lemma \ref
{lem:annealedPosExitProbEst}
and a decomposition into regenerations)
there exist random variables $Z_1^*$ and $V$ taking values in $\partial
_+ \mathcal{P}(0,L_k)$ and
$\{0,1\}$,
respectively, and such that
$P(V= 0) \in\mathcal{S}(\N)$ as a function in $L$ and
\[
Z_1' := (Z_1 - x) \mathbh{1}_{Z_0 \in\partial_+ \mathcal
{P}_1(x,L_k), V =1}
+ Z_1^* \mathbh{1}_{\{ Z_0 \notin\partial_+ \mathcal{P}_1(x,L_k) \}
\cup\{V= 0\}}
\]
is distributed according to $\mu_{0,0}^{L_k}$.
Let furthermore
\[
Z_0^* := (Z_0 -x) \mathbh{1}_{Z_0 \in\partial_+ \mathcal
{P}_1(x,L_k)} + Z_2' \mathbh{1}_{Z_0 \notin\partial_+ \mathcal{P}_1(x,L_k)}.
\]
As a consequence, there exists an $H_0$-valued random variable
independent from everything else such
that $P(\Vert U \Vert_1 \leq K_k)= 1$, $P(U \not= 0)$ is contained
in~$\mathcal{S}(\N)$ as a function
in $L$, and $E(Z_0^* + U) = EZ_1'$.
Set
\[
Z_0' := Z_0^* + U.
\]

Then, since $P(Z_0 \not= Z_1) \leq(2k-1)\lambda_k$ by assumption, we get
\begin{eqnarray*}
P(Z_0' \not= Z_1') &\leq& P(Z_0 \not= Z_1) + P\bigl( Z_0 \notin\partial_+
\mathcal{P}_1(x,L_k)\bigr) + P(U \not= 0)
+ P(V=0)\\
&\leq& 2k \lambda_k
\end{eqnarray*}
for $L$ large enough.
Furthermore, $P(\Vert Z_0' - Z_2' \Vert_1 \leq2kK_k) = 1$.
To check the remaining variance condition, note that
\begin{eqnarray*}
&&\sum_{y} \Vert y - E_{\mu_{0,0}^{L_k}} \Vert_1^2 \cdot
\vert P(Z_1' =y) - P(Z_0'=y) \vert\\
&&\qquad =\sum_{y} \Vert y - E_{\mu_{0,0}^{L_k}} \Vert_1^2 \\
&&\qquad\quad\hspace*{0pt}{} \times\bigl( \bigl\vert P \bigl( (Z_1 - x) \mathbh{1}_{Z_0 \in
\partial_+ \mathcal{P}_1(x,L_k), V =1}
+ Z_1^* \mathbh{1}_{\{ Z_0 \notin\partial_+ \mathcal{P}_1(x,L_k) \}
\cup\{V= 0\}} = y \bigr) \\
&&\qquad\quad\hspace*{40.7pt}{} - P\bigl( (Z_0 -x) \mathbh{1}_{Z_0 \in\partial_+ \mathcal
{P}_1(x,L_k)}
+ Z_2' \mathbh{1}_{Z_0 \notin\partial_+ \mathcal{P}_1(x,L_k)} + U =
y \bigr) \bigr\vert\bigr)\\
&&\qquad \leq(dL_k R_6(L_k))^2 \bigl( P\bigl(Z_0 \notin\partial_+ \mathcal
{P}_1(x,L_k)\bigr) + P(U \not= 0) + P(V = 0) \bigr) \\
&&\qquad\quad{} + \sum_{y} \Vert y - E_{\mu_{x,0}^{L_k}} \Vert_1^2
\cdot\vert\mu_{x,0}^{L_k}(y) - P(Z_0 = y) \vert\\
&&\qquad\leq2k\lambda_k \Var_{\mu_{0,0}^{L_k}}
\end{eqnarray*}
for $L$ large enough, where to obtain the last inequality we employed
the $((2k-1)\lambda_k, (2k-1)K_k)$-closeness of $\nu$ to $\mu
_{x,0}^{L_k}$ as well as Lemma
\ref{lem:VarBd}. Again, since the above computations are uniform in
the admissible choices
of~$k$ and~$x$, this yields the result.
\end{pf}

In order to construct the auxiliary walk, we need the
following result which guarantees that if boxes on a certain scale are
left in some way close
to the annealed distribution conditioned on leaving through the right
boundary part
of the boundary, then the same applies to the containing box on the
larger scale as well.
Essentially, this is Lemma 4.16 of \cite{Be-09}.
%
\begin{lemma} \label{lem:sumOfClose}
Let $\lambda\in(0,1)$, $L$ be large enough and $n \in\N$ such that
$ n \leq\lambda L$.
Furthermore, let $(\Delta_i)_{i=1}^n$ be random variables such that
for every $i$, the variable
$\Delta_i$ takes values in $\partial_+ \mathcal{P}(0,L)$ only, and,
conditioned on $\Delta_1,\ldots,\Delta_{i-1}$, the distribution of
$\Delta_i$
is $(\lambda,K)$-close to $\mu_{0,0}^L$.
In addition, assume $R_3(L) \leq K \leq L$.

Then for $S_n:=\sum_{i=1}^n \Delta_i$, the distribution of $S_n$ is
$(R_{9}(L)\lambda,4nK)$-close to~$\mu_{0,0}^{\sqrt{n}L}$.
\end{lemma}

The proof of this crucial lemma can be found from page \pageref
{pr:lem:sumOfClose} onward.

Now we rigorously construct the auxiliary walk $(Y_n)$ in environment
$\omega$ starting in $0$,
and denote the corresponding probability measure
by $P_{0,\omega}$ also.
For $k \in\{1, \ldots, \iota-1\}$ we
recursively define
%
%
\begin{eqnarray}
\label{eq:MDef}
&\displaystyle M_k \mbox{ as the smallest integer larger than or
equal to}&\nonumber\\[-10pt]\\[-10pt]
&\displaystyle L^{\beta-6\delta}+L^{2\chi} \mbox{ such that }L_{k+1}^2\mbox{ divides }
\sum_{j=1}^k M_j L_j^2.&
\nonumber
\end{eqnarray}
%
Note that $L_{k+1} = L_k {\lfloor L^\chi\rfloor}$ implies that $M_k
\leq
\lceil L^{\beta- 6\delta} \rceil+ 2\lceil L^{2\chi} \rceil$,
and that for every $k\in\{2,\ldots,\iota\}$, from
$x\cdot e_1-\sum_{j=1}^{k-1}M_jL_j^2\in L_k^2{\mathbb N}_0$ we can
infer that
$x\cdot e_1\in L_k^2{\mathbb N}$.
Define ${\mathcal{P}}^{(k)}(x) := {\mathcal{P}}(z(x,k), L_k)$,
%
%
\begin{eqnarray} \label{eq:kDef}
&&k(x) := \max\Biggl\{ k \in\{1, \ldots, \iota\} \dvtx
x\cdot e_1 - \sum_{j=1}^{k-1} M_j L_j^2 \in L_k^2 \N_0 \nonumber\\[-10pt]\\[-10pt]
&&\hspace*{155pt}\mbox{and }
{\mathcal{P}}^{(k)}(x) \mbox{ is good} \Biggr\}\nonumber
\end{eqnarray}
and
%
%
\begin{equation} \label{eq:kPrimeDef}
k'(x) := \max\Biggl\{ k \in\{1, \ldots, \iota\} \dvtx
x \cdot e_1 - \sum_{j=1}^{k-1} M_j L_j^2 \in L_k^2 \N
\Biggr\}
\end{equation}
with the maximum of the empty set defined to be $0$.
We now define the auxiliary random\vadjust{\goodbreak} walk $(Y_n)$ and a corresponding
sequence of stopping times~$(\zeta_n)$ recursively.
For $z \in\partial_+ {\mathcal{P}}(0,L_1)$ chosen according to $\mu
_{0,0}^L$,
fix $Y_0, \ldots, Y_{l_1}$
to be an arbitrary nearest-neighbor path (independent of $\omega$) of
shortest length connecting
0 with $z$ such that
$\{Y_0, \ldots, Y_{{l_1}-1}\} \subset{\mathcal{P}}(0,L_1)$.
Furthermore,
set $\zeta_1:= \zeta_1' := T^Y_{\partial_+ {\mathcal{P}}(0,L_1)} = l_1$.
Next, we define the recursive step of the construction.

\begin{enumerate}[(R)]
\item[(R)]\hypertarget{item:Recursion}
Assume that the walk is defined up to time $\zeta_n'$ and
set $x:= Y_{\zeta_n'}$.
\begin{enumerate}[$\bullet$]
\item[$\bullet$]
If $k(x) > 0$, then
choose $Y_{\zeta'_n+\cdot}$ according to the law of $X_\cdot$ with
respect to
\[
P_{x,\omega} \bigl( \cdot \vert
T_{\partial{\mathcal{P}}^{(k(x))}(x)} = T_{\partial_+ {\mathcal
{P}}^{(k(x))}(x)} \bigr),
\]
up to time $\zeta'_n + l_n$, where
\[
l_n:= T^{Y_{ \zeta'_n + \cdot}}_{\partial{\mathcal{P}}^{(k(x))}(x)}.
\]

\item[$\bullet$]
Otherwise, if $k(x)=0$, then
similarly to the start of the construction,
we choose
\[
\{Y_{\zeta'_n}, \ldots, Y_{\zeta'_n + l_n} \}
\]
to be
a nearest-neighbor path of shortest length connecting
$x$ with $z$, where $z$ is chosen according to
$\mu_{x,z(x,1)}^{L_1}$ and such that this path leaves $\mathcal
{P}^{(1)}(x)$ in its last step only.
\end{enumerate}
In both cases, set
$\zeta_{n+1} := \zeta'_n + l_n$.
If $Y_{\zeta_{n+1}} \cdot e_1 > L^{1+\delta}$, then we stop the
construction of $Y$.

If $1 \vee k(x) = k'(Y_{\zeta_{n+1}})$, then set $Y_{\zeta_{n+1}+1}
:= Y_{\zeta_{n+1}} + e_1$,
$Y_{\zeta_{n+1}+2} := Y_{\zeta_n}$ as well as $\zeta_{n+1}' := \zeta
_{n+1}+2$ and repeat step
\hyperlink{item:Recursion}{(R)}.

Otherwise, if $1 \vee k(x) < k'(Y_{\zeta_{n+1}})$,
given $\zeta_1, \ldots, \zeta_{n+1}$ and $(Y_i)_{i \in\{0, \ldots,
\zeta_{n+1}\}}$,
define for each $k \in\{(1 \vee k(x)) + 1, \ldots, k'(Y_{\zeta
_{n+1}})\}$ the number $j(k) :=
Y_{\zeta_{n+1}} \cdot e_1 /L_k^2$.
Furthermore,
define for $j,k \in\N$ the stopping time $T_k'(j)$ equal to $\zeta
_m'$ if there exists $m \leq n+1$ such that
$\zeta_m = T_k(j)$, and equal to $T_k(j)$ otherwise.

Now
for $k \in\{(1 \vee k(x))+1, \ldots, k'(Y_{\zeta_{n+1}})\}$ with
increasing order we iteratively
perform the following step, where $\zeta_{n+1}^{(1 \vee k(x))} :=
\zeta_{n+1}$:
\begin{enumerate}[(B)]
\item[(B)]\hypertarget{item:goodExitOfBigBox}
Conditioned\label{aaaa} on $Y_{T_{k}'(j(k)-1)}$,\vspace*{2pt}
by construction
(and as a consequence of Corollary \ref{cor:closeness} and Lemma \ref
{lem:sumOfClose}), the distribution of
the variable
\[
Y_{\zeta_{n+1}^{(k-1)}} - z\bigl(Y_{T_{k}'(j(k)-1)},k\bigr)
\]
is
$(2(k-1)\lambda_k, 2(k-1) K_k)$-close
to
\[
\mu_{Y_{T_{k}'(j(k)-1)} - z(Y_{T_{k}'(j(k)-1)},k) ,0}^{L_k}.
\]
We now condition the variable
\[
Y_{\zeta_{n+1}^{(k-1)}} - z\bigl(Y_{T_{k}'(j(k)-1)},k\bigr)
\]
on the event
%
%
\begin{equation} \label{eq:DkCond}
D_k:= \bigl\{
T^{Y_{T_{k}'(j(k)-1) + \cdot}}_{\partial_+ \mathcal
{P}_1(z(Y_{T_{k}'(j(k)-1)},k),L_k)}
= T^{Y_{T_{k}'(j(k)-1) + \cdot}}_{\partial\mathcal
{P}_1(z(Y_{T_{k}'(j(k)-1)},k),L_k)}
\bigr\}.
\end{equation}
In combination with Lemma \ref{lem:condCloseness}
we may infer that for $L$ large enough, the distribution of this
conditioned random variable
still is $((2k-1)\lambda_k, (2k-1) K_k)$-close to
\[
\mu_{Y_{T_{k}'(j(k)-1)} - z(Y_{T_{k}'(j(k)-1)},k) ,0}^{L_k}.
\]
Thus, Lemma \ref{lem:transCloseness} implies that
%
%
\begin{eqnarray} \label{eq:2kKkCloseness}
&&\mbox{the
distribution of the variable }
Y_{\zeta_{n+1}^{(k-1)}} - Y_{T_{k}'(j(k)-1)}\nonumber\\[-10pt]\\[-10pt]
&&\mbox{is } (2k\lambda_k, 2kK_k)\mbox{-close to } \mu
_{0 ,0}^{L_k}.
\nonumber
\end{eqnarray}
Set
$
\zeta_{n+1}^{(k)} := \zeta_{n+1}^{(k-1)} + \Vert\beta_{k,j(k)}
\Vert_1,
$
where the $\beta_{k,j(k)}$ defined
below take values in $H_0$ and play a correcting role.
Furthermore, let
$
Y_{\zeta_{n+1}^{(k-1)}}, \ldots,\allowbreak Y_{\zeta_{n+1}^{(k)}}
$
be a nearest-neighbor path of shortest length from
$
Y_{\zeta_{n+1}^{(k-1)}}
$
to
$
Y_{\zeta_{n+1}^{(k-1)}} + \beta_{k,j(k)}.
$
Note that from the conditioning in (\ref{eq:DkCond}) in combination
with Remark \ref{rem:betaBound} below,
we may infer that $Y_{\zeta_{n+1}^{(k)}}-Y_{T_{k}'(j(k)-1)}$ takes
values in $\partial_+ \mathcal{P}^{(k)} (0)$ only.
If $k < k'(Y_{\zeta_{n+1}})$, then repeat step \hyperlink
{item:goodExitOfBigBox}{(B)} for $k+1$;
if $k = k'(Y_{\zeta_{n+1}})$, continue below.
\end{enumerate}
Set
$Y_{\zeta_{n+1}^{(k)} + 1} := Y_{\zeta_{n+1}^{(k)}} + e_1$
as well as
$Y_{\zeta_{n+1}^{(k)} + 2} := Y_{\zeta_{n+1}^{(k)}}$
and $\zeta_{n+1}' := \zeta_{n+1}^{(k)} + 2$.
Now we continue the construction at the recursion step
\hyperlink{item:Recursion}{(R)}.
\end{enumerate}

It remains to define the variables $\beta_{k,j}$.
Set $\beta_{1,j} = 0$ for all $j$.
For any $n \in\N$, we will define those $\beta_{k,j}$, $k \in\{2,
\ldots, \iota\}$, for which $Y_{\zeta_{n}} \in H_{jL_k^2}$,
using only the environment $\omega$, the auxiliary walk $Y$ up to time
$\zeta_{n}$
as well as
the values of $\{\beta_{\hat{k}, \hat{j}} \dvtx\hat{k} \in\{2, \ldots
, k-1\} \mbox{ and } \hat{j}L_{\hat{k}}^2 = jL_k^2\}$.
We define $\beta_{k,j}$ to be $0$ in the following cases:
\begin{itemize}
\item If there is no $n \in\N$ such that $\zeta_n = T_k(j-1)$, then
$\beta_{k,j} = 0$.

\item Otherwise, let $n$ be such that $\zeta_n = T_k(j-1)$. If
${\mathcal{P}}^{(k)}(Y_{\zeta_n'})$ is good, then $\beta_{k,j} = 0$.

\end{itemize}
Thus, assume now that $\zeta_n = T_k(j-1)$ such that $\mathcal
{P}^{(k)}(Y_{\zeta_n'})$ is bad.
Let $x := Y_{\zeta_n'}$ and let $\mu^{k,j,Y}_{x,\omega}$ be the
distribution\vspace*{-1pt} of the variable
$
Y_{\zeta_{n+1}^{(k-1)}} - x,
$
which due to~(\ref{eq:2kKkCloseness}) is $(2k\lambda_k, 2kK_k)$-close
to $\mu_{0,0}^{L_k}$.
Thus, we find $(Z_0,Z_1,Z_2)$ defined on the same probability space as
$(Y_n)$ (which without
loss of generality is assumed to be large enough)
such that $Z_2$ equals $Y_{\zeta_{n+1}^{(k-1)}} - x$ (cf. Remark~\ref
{rem:couplingSpace}),\vadjust{\goodbreak}
such that
$Z_1 \sim\mu_{0,0}^{L_k}$,
and such that furthermore the requirements of $(2k\lambda_k,
2kK_k)$-closeness (cf. Definition~\ref{def:closeness}) are satisfied.
Now define
%
%
\begin{equation} \label{eq:betakjDef}
\beta_{k,j} := Z_0 - Z_2,
\end{equation}
and note that $\beta_{k,j} \in H_0$ a.s.
This completes the definition of $\beta_{k,j}$.

The deterministic corrections caused by the variables $\beta_{k,j}$
are not too
big, that is, not too expensive in terms of probability. This is made
precise in the following remark.\vspace*{-2pt}
%
\begin{remark} \label{rem:betaBound}
By construction of $(Y_n)$, for every $k \in\{2, \ldots, \iota\}$ and
$j \in\N$ such that $\beta_{k,j}$ has been defined above,
with probability $1$,
\[
\beta_{k,j} < 2\iota K_\iota\leq L^{4\chi}
\]
for $L$ large.\vspace*{-2pt}
\end{remark}
%
\begin{remark} \label{rem:YProp}
Observe that by construction
we infer that
$T^Y_{\partial C_L} = T^Y_{\partial_+ C_L}$ a.s., with $C_L$ denoting
the set of Proposition \ref{prop:quenchedConeTypeExit}.\vspace*{-2pt}
\end{remark}

\subsection{Random direction event} \label{subsec:RDE}
As in \cite{Be-09}, we will introduce a so-called \textit{random
direction event}
in order to ensure that, in most environments, the walker does not hit
too many bad boxes.
For this purpose, for $k\in\{1,\ldots,\iota\}$ set
%
%
\begin{equation} \label{eq:BDef}
B_{k} := \frac{\sum_{j=1}^{k-1} M_j L_j^2}{L^2_{k}}.
\end{equation}
For $w \in[-1,1]^{d-1}$, $k \in\{2, \ldots, \iota\}$ and
$j \in\{B_k +1, \ldots, M_k\}$, define
\[
W_k^{(w)}(j):=
\bigl\{ \bigl\Vert Y_{T_k'(j)} - Y_{T_k'(B_k)} - (j-B_k) \bigl( E_{\mu
_{0,0}^{L_k}} - L_k (0,w) \bigr)
\bigr\Vert_\infty< L_k \bigr\},
\]
where in a slight abuse of notation we write $L_k (0,w)$ to denote the vector
$(0,L_k w_1, \ldots, L_k w_{d-1}) \in\R^d$.
Furthermore, define
\[
W_k^{(w)} := \bigcap_{j = B_k + 1}^{B_k + M_k} W_k^{(w)}(j)
\]
as well as the \textit{random direction event}
\[
W^{(w)} := \bigcap_{k=1}^\iota W_k^{(w)}.
\]
To obtain a lower bound for the probability of this event, we have to
establish some auxiliary results
first.\vspace*{-2pt}
%
\begin{claim} \label{claim:AWCloseness}
For all $L$ large enough and all $k \in\{1, \ldots, \iota\}$, $j\in
\{B_{k}+1, \ldots, B_k + M_k\}$ and
$\omega\in\Omega$, one has that
$P_{0,\omega}( \cdot \vert Y_1, \ldots, Y_{T'_k(j-1)})$-a.s.
the distribution
of $Y_{T'_k(j)}-Y_{T'_k(j-1)}$
is $(2k \lambda_k, 2kK_k)$-close to $\mu_{0,0}^{L_k}$.\vspace*{-2pt}
\end{claim}

\begin{pf}
Similarly to Lemma 6.6 of \cite{Be-09}, this result is a consequence
of the construction of the auxiliary walk.\vadjust{\goodbreak}
In fact, if $\mathcal{P}^{(k)}(Y_{T_k'(j-1)})$ is
good, then the statement follows from the first part of step \hyperlink
{item:Recursion}{(R)}
in the construction of the auxiliary walk in combination
with Corollary \ref{cor:closeness}.

Otherwise, if $\mathcal{P}^{(k)}(Y_{T_{k}'(j-1)})$ is bad, it follows
from step \hyperlink{item:goodExitOfBigBox}{(B)} of that construction.
\end{pf}

We now get the following corollary.
%
\begin{corollary}[(Corollary 6.7 of \cite{Be-09})]\label{cor:entropy}
There exists a constant $\rho>0$ such that for all $L$ large enough,
$\omega\in\Omega$,
all
$k,j$ as\vspace*{-1pt} in Claim \ref{claim:AWCloseness},
$\overline{Y} := Y_{T'_k(j-1)} + E_{\mu_{0,0}^{L_k}}$, and
for
all $x\in H_{jL_k^2}$ such that $\Vert\overline{Y} - x \Vert_1 < 4L_k$,
one has
%
%
\begin{equation}\label{eq:ent}
P_{0,\omega} \bigl(
\bigl\Vert Y_{T^\prime_k(j)} - x \bigr\Vert_1 < L_k
\vert Y_1, \ldots, Y_{T^\prime_k(j-1)} \bigr) > \rho.
\end{equation}
\end{corollary}
\begin{pf}
This follows from Claim \ref{claim:AWCloseness} in combination with Lemmas
\ref{lem:bddFluctExp} and~\ref{lem:lbound}.
\end{pf}
%
\begin{lemma}[(Lemma 7.1 of \cite{Be-09})] \label{RDEprobBdLem}
There exists $\rho> 0$ such that for
all $L$ large enough, $\omega\in\Omega$, all $w \in[-1,1]^{d-1}$,
as well as $k,j$ as in Claim \ref{claim:AWCloseness},
one has
\[
P_{0,\omega} \bigl(
W_k^{(w)}(j) \vert W_1^{(w)}, \ldots, W_{k-1}^{(w)},
W_k^{(w)}(B_k +1), \ldots, W_k^{(w)}(j-1)
\bigr) > \rho
\]
[with $W_k^{(w)}(B_k) := \Omega$].
\end{lemma}
\begin{pf}
On the event
\[
W_1^{(w)} \cap\cdots\cap W_{k-1}^{(w)} \cap W_k^{(w)} (B_k+1) \cap
\cdots\cap W_k^{(w)}(j-1)
\]
one has
\[
\bigl\Vert Y_{T_k'(j-1)} - Y_{T'_k(B_k)} - (j-1-B_k) \bigl( E_{\mu_{0,0}^{L_k}}
- L_k(0,w) \bigr) \bigr\Vert_\infty
< L_k
\]
and thus
\begin{eqnarray*}
&&\bigl\Vert\underbrace{Y_{T'_k(B_k)} + (j-B_k) \bigl( E_{\mu
_{0,0}^{L_k}} + L_k(0,w) \bigr)}_{=: x}
-\bigl(Y_{T'_k(j-1)} + E_{\mu_{0,0}^{L_k}}\bigr) \bigr\Vert_\infty\\
&&\qquad<2L_k.
\end{eqnarray*}
Corollary \ref{cor:entropy} now yields the desired result.
\end{pf}

Departing from this result we obtain the desired lower bound on the probability
of the random direction event.
%
\begin{lemma} \label{lem:RDETotalProbBdLem}
There exists a constant $C > 0$ such that for
all $L$ large enough as well as all $\omega\in\Omega$ and $w \in
[-1,1]^{d-1}$,
\[
P_{0,\omega}\bigl(W^{(w)}\bigr) \geq e^{-CL^{\beta- 6\delta}}.
\]
\end{lemma}
\begin{pf}
We compute
\begin{eqnarray*}
&&P_{0,\omega} \bigl(W^{(w)}\bigr) \\
&&\qquad= \prod_{k=1}^\iota\prod_{j= B_k +1 }^{B_k + M_k}
P_{0,\omega} \bigl( W_k^{(w)} (j) \vert W_1^{(w)}, \ldots,
W_{k-1}^{(w)},\\
&&\hspace*{149pt}\hspace*{-42.7pt} W_{k}^{(w)}(B_k + 1), \ldots, W_k^{(w)} (j-1) \bigr)\\
&&\qquad\geq\rho^{\sum_{k=1}^\iota M_k} \geq e^{ - CL^{\beta- 6\delta}}
\end{eqnarray*}
for $C > 0$ large enough, where the first inequality is a consequence
of Lem\-ma~\ref{RDEprobBdLem} while the
second follows from
the bound $M_k \leq2 \lceil L^{\beta- 6\delta} \rceil$ for $L$
large enough; see directly after (\ref{eq:MDef}).
\end{pf}

We now want to bound from above the probability that the auxiliary walk
hits too many bad boxes.
For this purpose, we start with the following auxiliary result.
%
\begin{lemma}[(Lemma 7.4 of \cite{Be-09})] \label
{lem:boundIntersectionProbRestr}
For all $L$ large enough, $\omega\in\Omega$, $k \in\{1, \ldots,
\iota-1\}$, $j \in\{ B_{k+1} {\lfloor L^\chi\rfloor}^2, \ldots,
{\lfloor
L^{1+\delta} / L_k^2 \rfloor} \}
$
and $z \in\mathcal{L}_{L_k} \cap H_{jL_k^2}$ one has
%
%
\begin{eqnarray} \label{eq:boundIntersectionProbRestr}\qquad
&&\int_{[-1,1]^{d-1}} P_{0,\omega} \bigl( \bigl\{Y_n\dvtx n \in\{1, \ldots,
T_{L^{1+\delta}}^Y\}\bigr\} \cap{\mathcal{P}}(z,L_k)
\not= \varnothing \vert W^{(w)} \bigr) \,dw \nonumber\\[-8pt]\\[-8pt]
&&\qquad\leq L^{(-\beta+ 6\delta+ 2\chi)(d-1)}.
\nonumber
\end{eqnarray}
\end{lemma}
\begin{pf}
Choose $k'$ to be the number out of $\{k, \ldots, \iota-1\}$ such
that $B_{k'}L_{k'}^2 \leq jL_k^2 < B_{k'+1} L^2_{k'+1}$.
We start with noting that for fixed $w \in[-1,1]^{d-1}$,
with probability $1$ with respect to $P_0( \cdot \vert
W^{(w)})$, the walk $Y$ is located in a $(d-1)$-dimensional hypercube
of side length $\sum_{j=1}^{k'-1} L_j \leq\iota L_{k'-1}$ at time~%
$T'_{k'}(B_{k'})$. Letting $w$ vary over $[-1,1]^{d-1}$,
the union of all appearing hypercubes covers a hypercube of side length
at least $M_{k'-1}L_{k'-1} \geq\lceil L^{\beta-6\delta} \rceil L_{k'-1}$.

Now let $\{y_1, \ldots, y_r\} \subset\mathcal{L}_{L_{k'}}$ be the
set of all elements $y_j \in\mathcal{L}_{L_{k'}}$
such that $\mathcal{P}(z,L_k) \cap\mathcal{P}(y_j,L_{k'}) \not=
\varnothing$
for all $j \in\{1, \ldots, r\}$,
and note that, due to a reasoning
similar to the observation just before
Lemma \ref{lem:GProbEstimate},
$r$ is bounded from above by $3\cdot15^{d-1}$.

From steps \hyperlink{item:Recursion}{(R)} and \hyperlink
{item:goodExitOfBigBox}{(B)} in the construction of the auxiliary walk $Y$,
it follows that if there exists $\zeta_n'$ such that $z(Y_{\zeta
_n'},k') = y_j$, then $Y$ leaves~$\mathcal{P}(y_j,L_{k'})$
through $\partial_+ \mathcal{P}(y_j,L_{k'})$.
Therefore, we conclude that
\[
\bigl\{ \{Y_n \dvtx n \in\N\} \cap{\mathcal{P}}(z,L_k) \not=
\varnothing\bigr\}
\subset\bigcup_{j=1}^r \bigl\{ \{Y_n \dvtx n \in\N\} \cap{\mathcal
{P}}(y_j,L_{k'}) \not= \varnothing\bigr\}.
\]
But due to the above reasoning,
there exists a constant $C$ such that the right-hand side
can have positive probability
with respect\vadjust{\goodbreak} to $P_{0,\omega}( \cdot \vert W^{(w)})$ only if
$w$ lies in a certain $(d-1)$-dimensional hypercube of side length
\[
\frac{CR_6(L_{k'})L_{k'}}{L^{\beta-6\delta} L_{k'-1}} \leq
CL^{-\beta+ 6\delta+ 3\chi/2}.
\]
This establishes (\ref{eq:boundIntersectionProbRestr}).
\end{pf}

Now adopt the notation
\[
{\mathcal{D}}_{k,\omega} := \{ x \in C_L \cap\mathcal{L}_{L_k} \dvtx x
\cdot e_1 \geq B_{k} L_{k}^2 \mbox{ and }
{\mathcal{P}}(x,L_k) \mbox{ is bad with respect to }\omega\}
\]
and
%
%
\begin{equation} \label{eq:NrBadBoxesHit}
{\mathcal{B}}_{k,\omega} := \bigl\vert\bigl\{ x \in{\mathcal
{D}}_{k,\omega} \dvtx
\bigl\{Y_n \dvtx n \in\{1, \ldots, T_{L^{1+\delta}}^Y\}\bigr\}
\cap{\mathcal{P}}(x, L_k) \not= \varnothing\bigr\} \bigr\vert.
\end{equation}
We are interested in the distribution of the variable ${\mathcal
{B}}_{k,\omega}$. Recall that
$\Theta_L$ has been defined in (\ref{eq:ThetaDef}).
%
\begin{lemma}[(Lemma 7.5 of \cite{Be-09})]\label{lem:NrBadPointsVisitedBound}
For all $L$ large enough and all $k \in\{1, \ldots, \iota-1\}$ as
well as $\omega\in\Theta_L$,
\[
\int_{[-1,1]^{d-1}} E_{0,\omega} \bigl( {\mathcal{B}}_{k,\omega}
\vert W^{(w)} \bigr) \,dw
\leq15^d L^{\beta-6\delta}.
\]
\end{lemma}
\begin{pf}
With the same reasoning as in the proof of Lemma \ref
{lem:boundIntersectionProbRestr},
steps~\hyperlink{item:Recursion}{(R)} and \hyperlink
{item:goodExitOfBigBox}{(B)} of the construction of the auxiliary walk $Y$
imply that $P_0( \cdot \vert W^{(w)})$-a.s. we have
%
%
\begin{eqnarray} \label{eq:firstBoxesEst}
&&\bigl\vert\bigl\{x \in C_L \cap\mathcal{L}_{L_k}\dvtx B_k L_k^2 \leq x
\cdot e_1 < B_{k+1} L_{k+1}^2,\nonumber\\
&&\hspace*{17pt}\mathcal{P}(x, L_k) \cap\bigl\{Y_n \dvtx n \in\{1, \ldots, T_{L^{1+\delta
}}^Y\}\bigr\} \not= \varnothing\bigr\} \bigr\vert\\
&&\qquad\leq3 \cdot15^{d-1} M_k \leq(15^d-1) L^{\beta-6\delta}.
\nonumber
\end{eqnarray}
Now consider $x \in\mathcal{D}_{k,\omega}$ with $x \cdot e_1 \geq
B_{k+1} L_{k+1}^2$.
Then by Lemma \ref{lem:boundIntersectionProbRestr},
%
%
\begin{eqnarray} \label{eq:secondBoxesEst}\quad
&&\int_{[-1,1]^{d-1}} P_{0,\omega} \bigl(
\bigl\{Y_n\dvtx n \in\{1, \ldots, T_{L^{1+\delta}}^Y\}\bigr\}\cap{\mathcal
{P}}(x,L_k) \not= \varnothing \vert W^{(w)}
\bigr) \,dw\nonumber\\[-8pt]\\[-8pt]
&&\qquad\leq L^{(-\beta+ 6\delta+ 2\chi)(d-1)}
\nonumber
\end{eqnarray}
for $L$ large enough.
Therefore, (\ref{eq:firstBoxesEst}) and (\ref{eq:secondBoxesEst}) in
combination with (\ref{eq:ThetaDef}) yield
\begin{eqnarray*}
\int_{[-1,1]^{d-1}} E_{0,\omega} \bigl({\mathcal{B}}_{k,\omega} \vert
W^{(w)} \bigr) \,dw
&\leq&(15^d-1) L^{\beta-6\delta} + L^{(-\beta+ 6\delta+ 2\chi
)(d-1)} L^{\alpha+ \delta}\\
&\leq&15^d L^{\beta-6\delta},
\end{eqnarray*}
due to our choice of $\delta$.
\end{pf}

Because of the modifications in our construction of the auxiliary walk
in comparison to the one
in \cite{Be-09}, we give here a modified result concerning the density
of the path measures of $X$
with respect to $Y$.\vadjust{\goodbreak}
%
\begin{lemma}[(Lemma 6.5 of \cite{Be-09})] \label{lem:densityEst}
Let $(v_n)=(v_1, \ldots, v_{T^v_{L^{1+\delta}}})$ be a finite
nearest-neighbor path in $\Z^d$ starting in $0$
such that $T^v_{L^{1+\delta}} = \inf\{n \in\N\dvtx v_n \cdot e_1 >
L^{1+\delta}\}$.
Furthermore, for $k \in\{1, \ldots, \iota\}$ and $\omega\in\Omega
$, let
\[
Q_{k,\omega}(v) := \bigl\vert\bigl\{ z \in\mathcal{D}_{k,\omega} \dvtx
\bigl\{v_n \dvtx n \in\{1, \ldots, T^v_{L^{1+\delta}}\}\bigr\}
\cap\mathcal{P}(z,L_k) \not= \varnothing\bigr\} \bigr\vert\vspace*{3pt}
\]
and set
$
Q_\omega(v) := \sum_{k=1}^\iota Q_{k,\omega}(v).
$

Then for all $L$ large enough and all $\omega\in\Omega$ we have
%
%
\begin{equation} \label{eq:pathDensityEst}\qquad
\frac{P_{0,\omega}(X_j = v_j\ \forall j \in\{1, \ldots,
T^v_{L^{1+\delta}}\})}
{P_{0,\omega}(Y_j = v_j\ \forall j \in\{1, \ldots, T^v_{L^{1+\delta
}}\})}
\geq\frac12 \kappa^{3Q_\omega(v)\iota L^{9\psi/4} + 4\iota
L^{\beta-6\delta}}\vspace*{3pt}
\end{equation}
for all admissible choices of $(v_n)$.\vspace*{3pt}
\end{lemma}
\begin{pf}
Due to ellipticity, the numerator in (\ref{eq:pathDensityEst}) is
positive; therefore, it is sufficient to
consider those trajectories $(v_n)$ only for which the
probability in the denominator is positive as well.

To any such $(v_n)$ and environment $\omega$, there belong sequences
$(\zeta_n)$ and
$(\zeta_n')$ as in the definition of $Y$.
In fact, set $\zeta_0 := \zeta_0' := 0$ and $\zeta_1 := \zeta_1' :=
T^v_{\partial{\mathcal{P}}(0,L_1)}$. Given $\zeta_0, \ldots, \zeta
_n$ and $\zeta_0', \ldots, \zeta_{n-1}'$,
define $x_n := v_{\zeta_n}$, $\zeta_n' := \min\{l > \zeta_n \dvtx
v_{l-1} \cdot e_1 > x_n \cdot e_1 \}$ (only if $n > 1$)
as well as $x_n' := v_{\zeta_n'}$.
For $k(x_n')$ as in (\ref{eq:kDef}), if $k(x_n') > 0$, set $\zeta
_{n+1} := T^v_{\partial_+ {\mathcal{P}}^{(k)}(x_n')}$,\vspace*{2pt}
otherwise set $\zeta_{n+1} := T^v_{\partial_+ {\mathcal{P}}^{(1)}(x_n')}$.

Now to estimate the probability in the denominator from above, we only
consider the
contributions coming from $Y$ moving in good boxes in which it behaves
like the quenched walk $X$ conditioned
on leaving the box through its right boundary part:
\begin{eqnarray*}
&&P_{0,\omega} (Y_j = v_j\ \forall j \in\{1, \ldots,
T^v_{L^{1+\delta}}\})\\[1pt]
&&\qquad\leq\prod_{n \dvtx k(x_n') > 0} P_{x_n', \omega}
\bigl( X_l = v_{\zeta_n' + l}\\[1pt]
&&\hspace*{102pt} \forall l \in\{ 1, \ldots, \zeta
_{n+1} - \zeta_n'\} 
\vert T_{\partial{\mathcal{P}}^{(k(x_n'))}(x_n')} =
T_{\partial_+ {\mathcal{P}}^{(k(x_n'))}(x_n')} \bigr).\vspace*{3pt}
\end{eqnarray*}
To obtain a lower bound for the numerator, as a consequence of the
strong Markov property
we may decompose
it into movements within the corresponding boxes as follows:
\begin{eqnarray*}
&&P_{0, \omega} (X_j = v_j\ \forall j \in\{1, \ldots,
T^v_{L^{1+\delta}}\})\\[1pt]
&&\qquad\geq\prod_{n \dvtx k(v_{\zeta_n'}) > 0} P_{v_{\zeta_n'}, \omega}
\bigl( X_l = v_{ \zeta_n'+l}\\[1pt]
&&\hspace*{90pt} \forall l \in\{ 1, \ldots, \zeta_{n+1} -
\zeta_n'\}
\vert
T_{\partial{\mathcal{P}}^{(k(v_{\zeta_n'}))}(v_{\zeta_n'})} =
T_{\partial_+ {\mathcal{P}}^{(k(v_{\zeta_n'}))}(v_{\zeta_n'})}
\bigr)\\[1pt]
&&\qquad\quad{} \times\prod_{n \dvtx k(v_{\zeta_n'})>0} P_{v_{\zeta_n'}, \omega}
\bigl( T_{\partial{\mathcal{P}}^{(k(v_{\zeta_n'}))}(v_{\zeta_n'})}
= T_{\partial_+ {\mathcal{P}}^{(k(v_{\zeta_n'}))}(v_{\zeta_n'})}
\bigr)\\[1pt]
&&\qquad\quad{} \times(\kappa^{\iota L^{3\chi}})^{Q_\omega(v)}
(2^{-\iota} )^{Q_\omega(v)}\\[1pt]
&&\qquad\quad{} \times\prod_{n\dvtx\zeta_n' < T_{L^{1+\delta}}^v} \kappa^2
\kappa^{CL^{2\psi}} \prod_{n \dvtx k(v_{\zeta_n'}) = 0}
\kappa^{CL^{2\psi}}
\end{eqnarray*}
for $L$ large enough
as well as $k(v_{\zeta_n'})$ and $k'(v_{\zeta_n'})$ as defined
in~(\ref{eq:kDef}) and~(\ref{eq:kPrimeDef}).
Here, the first and second product on the right-hand side come from~$X$
moving in good boxes.
The third and fourth factor on the right-hand side
originate from the corrections in the case of moving in bad boxes. In
this case,
Remark \ref{rem:betaBound} tells us that each of the correcting
variables $\beta_{k,j}$ is bounded from above by $L^{3\chi}$.
Since each time such a correction occurs, the number of influencing
correcting variables
$\beta_{k,j}$ is bounded from above by $\iota$, we obtain the third factor.
The fourth factor originates from the conditioning on $D_k$ in (\ref
{eq:DkCond}),
the probability of
which can be estimated using Lemma \ref{lem:goodExitOfBigBox}.
The fifth factor follows from the fact that directly before each time
$\zeta_n'$ we force the walk
to do one step in the direction of $e_1$ and one step back,
while the last factor originates from the deterministic moves performed
within bad boxes of scale one.
Consequently, we obtain
%
%
\begin{eqnarray} \label{eq:badProd}
\qquad&&\frac{P_{0,\omega}(X_j = v_j\ \forall j \in\{1, \ldots,
T^v_{L^{1+\delta}}\})}
{P_{0,\omega}(Y_j = v_j\ \forall j \in\{1, \ldots, T^v_{L^{1+\delta
}}\})}\nonumber\\
\qquad&&\qquad\geq\prod_{n \dvtx k(x_n')>0} P_{x_n', \omega}
\bigl( T_{\partial{\mathcal{P}}^{(k(x_n'))}(x_n')} = T_{\partial_+
{\mathcal{P}}^{(k(x_n'))}(x_n')} \bigr)
\bigl(\kappa^{3Q_\omega(v)\iota L^{3\chi}}\bigr)\\
\qquad&&\qquad\quad{} \times\prod_{n\dvtx\zeta_n' < T_{L^{1+\delta}}^v} \kappa
^2\kappa^{CL^{2\psi}} \prod_{n \dvtx k(v_{\zeta_n'}) = 0}
\kappa^{CL^{2\psi}}
\nonumber
\end{eqnarray}
for $L$ large enough.
Since $k(v_{\zeta_n'})>0$ implies that ${\mathcal{P}}^{(k(v_{\zeta
_n'}))}(v_{\zeta_n'})$ is good,
from (\ref{eq:quenProbRightExitEst}) we infer that
the value of the first product on the right-hand side is bigger than $1/2$
uniformly in all $(v_n)$ we consider, for all $L$ large enough.
Due to the construction of the auxiliary walk $Y$, there are
at most $\sum_{k=1}^\iota M_k \leq2\iota L^{\beta-6\delta}$
stopping times $\zeta_n'$ such
that $\zeta_n' < T_{L^{1+\delta}}^v$.
Therefore, and due to the choice of $\delta$ and $\psi$, for $L$
large enough,
the total expression on the right-hand side is bounded from below by
$
\kappa^{3Q_\omega(v)\iota L^{9\psi/4} + 4\iota L^{\beta-6\delta}}
/ 2,
$
which finishes the proof.
\end{pf}

\subsection{\texorpdfstring{Proof of Proposition \protect\ref{prop:quenchedConeTypeExit}}{Proof of Proposition 2.1}}
With $\mathcal{B}_{k,\omega}$ as defined in (\ref{eq:NrBadBoxesHit})
and for $L$ large enough, Lemma \ref{lem:NrBadPointsVisitedBound} yields
\[
\int_{[-1,1]^{d-1}} E_{0,\omega} \Biggl(\sum_{k=1}^\iota\mathcal
{B}_{k,\omega} \big| W^{(w)} \Biggr) \,dw \leq15^d \iota
L^{\beta-6\delta}
\]
for $\omega\in\Theta_L$.
Hence, for such $\omega$ and $L$ fixed, we can find $w \in
[-1,1]^{d-1}$ such that
%
%
\begin{equation} \label{RDExpBd}
E_{0,\omega} \Biggl( \sum_{k=1}^\iota\mathcal{B}_{k,\omega}
\big|  W^{(w)} \Biggr) \leq15^d \iota L^{\beta- 6\delta}.
\end{equation}
Fix such $w$ and define
\[
\overline{W} := \Biggl\{ \sum_{k=1}^\iota\mathcal{B}_{k,\omega}
\leq2\cdot15^d \iota L^{\beta- 6\delta} \Biggr\} \cap W^{(w)}.
\]
Using (\ref{RDExpBd}),
Markov's inequality yields
\[
P_{0,\omega} \Biggl( \Biggl\{ \sum_{k=1}^\iota\mathcal{B}_{k,\omega}
\geq2\cdot15^d \iota L^{\beta- 6\delta} \Biggr\} \bigg\vert
W^{(w)} \Biggr)
\leq\frac12,
\]
whence we obtain
%
%
\begin{eqnarray}
\label{overlineWProbEst}
P_{0,\omega} (\overline{W})
&=& P_{0,\omega} \Biggl( \Biggl\{ \sum_{k=1}^\iota\mathcal
{B}_{k,\omega}
\leq2\cdot15^d \iota L^{\beta- 6\delta} \Biggr\} \bigg\vert
W^{(w)} \Biggr) P_{0,\omega} \bigl(W^{(w)}\bigr)\nonumber\\[-8pt]\\[-8pt]
&\geq&\frac12 P_{0,\omega} \bigl(W^{(w)}\bigr) \geq e^{-CL^{\beta- 6\delta}}
\nonumber
\end{eqnarray}
for $L$ large enough and
where the last inequality follows from Lemma \ref{lem:RDETotalProbBdLem}.

We now observe that there is a set
$V_{L,\omega}$ of paths such that $\overline{W} = \{ (Y_n) \in
V_{L,\omega}\}$
and in particular, for $(v_n) \in V_{L,\omega}$ we have
$
Q_\omega(v) \leq2 \cdot15^d \iota L^{\beta- 6\delta}.
$
Thus, as a consequence of (\ref{eq:psiDef}) and
Lemma \ref{lem:densityEst},
%
%
\begin{eqnarray} \label{eq:XGoodExitEst}
P_{0,\omega} \bigl((X_n) \in V_{L,\omega}\bigr) &\geq& e^{-L^{\beta- \delta}/2}
P_{0,\omega} \bigl((Y_n) \in V_{L,\omega}\bigr)\nonumber\\[-8pt]\\[-8pt]
&=& e^{-L^{\beta- \delta}/2} P_{0,\omega} (\overline{W})
\geq e^{-L^{\beta- \delta}}\nonumber
\end{eqnarray}
for $L$ large enough,
where
the first inequality follows from the fact that $\omega\in\Theta_L$
in combination with Lemma
\ref{lem:densityEst} and our choices of $\delta$ and $\psi$, while
the last estimate follows from (\ref{overlineWProbEst}).
Due to Remark \ref{rem:YProp},
we may and do choose~$V_{L,\omega}$ in such a way that it only contain
paths that start in $0$ and leave~$C_L$ through $\partial_+ C_L$.
We take the required family of events in Proposition~\ref
{prop:quenchedConeTypeExit}
as $\Xi_L := \Theta_L$, and observe that from
(\ref{eq:XGoodExitEst}) and
Lemma \ref{lem:GProbEstimate} we can infer that~$\Xi_L$ has the
desired properties.

\begin{appendix}\label{sec:Appendix}
\section*{\texorpdfstring{Appendix: Auxiliary results and proof of Proposition~\lowercase{\protect\ref{prop:closenessbase}}}
{Appendix: Auxiliary results and proof of Proposition 3.4}}

This section contains slight modifications of auxiliary results
proven in \cite{Be-09} as well as some further lemmas.
With respect to results to which the first point applies, this section
is very much based
on \cite{Be-09}.

In order to prove Proposition \ref{prop:closenessbase}, we will
proceed as outlined in
Remark~\ref{rem:PropProofStrat}.

\subsection{\texorpdfstring{Proof of Proposition \protect\ref{prop:closenessbase}\protect\hyperlink{item:quenchedBadExitEst}{(i)}}
{Proof of Proposition 3.4(i)}}
Set
\[
G_L^{\mbox{\hyperlink{item:quenchedBadExitEst}{\textup{\iitem}}}} := \Bigl\{ \omega\in
\Omega\dvtx\max_{z \in\tilde{{\mathcal{P}}}(0,L)}
P_{z, \omega}\bigl(T_{\partial{\mathcal{P}}(0,L)} \not= T_{\partial_+
{\mathcal{P}}(0,L)}\bigr) \leq
e^{- R_1(L)^\gamma} \Bigr\}.
\]
Then Markov's inequality in combination with Lemma \ref
{lem:annealedPosExitProbEst}
yields
\begin{eqnarray*}
\PP\bigl( {G_L^{\mbox{\hyperlink{item:quenchedBadExitEst}{\textup{\iitem}}}}}^c \bigr)
&\leq& e^{R_1(L)^\gamma}
\E\max_{z \in\tilde{{\mathcal{P}}}(0,L)} P_{z,\omega}\bigl(T_{\partial
{\mathcal{P}}(0,L)} \not= T_{\partial_+ {\mathcal{P}}(0,L)}\bigr)\\
&\leq& e^{R_1(L)^\gamma}
\sum_{z \in\tilde{{\mathcal{P}}}(0,L)} P_z\bigl(T_{\partial{\mathcal
{P}}(0,L)} \not= T_{\partial_+ {\mathcal{P}}(0,L)}\bigr)\\
&\leq& e^{R_1(L)^\gamma} Ce^{-C^{-1} R_2(L)^\gamma} \leq Ce^{-C^{-1}
R_2(L)^\gamma}.
\end{eqnarray*}
In combination with Remark \ref{rem:PropProofStrat},
this finishes the proof.

\subsection{\texorpdfstring{Auxiliary results for the proof of Proposition
    \protect\ref{prop:closenessbase}\protect\hyperlink{item:quenchedAnnealedExpDiff}{(ii)} and
    \protect\hyperlink{item:quenchedAnnealedProbDiff}{(iii)}}
{Auxiliary results for the proof of Proposition 3.4(ii) and (iii)}}

We need the following local CLT-type results.
\setcounter{theorem}{0}
\begin{claim}\label{claim:sum_der}
Let $(Y_i)_{i \in\N}$ be $\Z^d$-valued, 
independent random variables with finite $(m+1)$st moments for some $m
\geq3$.
Furthermore, assume that $(Y_i)_{i\geq2}$ are identically distributed
and that
there exists $v\in\Z^d$ such that $P(Y_2=v)>0$ and $P(Y_2=v+e_j)>0$
for all $j\in\{1,\ldots,d\}$.
Let $\Gamma$ denote the covariance matrix of $Y_2$ and
$S_n=\sum_{i=1}^n (Y_i - EY_i)$. Then there exists a~constant $C$
which is determined by the
distributions of $Y_1$ and $Y_2$ such that for all $n \in\N$ and all
$x$, $y$ and $z \in\Z^d$
with $\Vert x-y \Vert_1 =1$ and $z-y=y-x$:
\begin{longlist}[(a)]
\item[(a)]
%
%
\setcounter{equation}{0}
\begin{equation} \label{eq:sumHPupperBd}
P(S_n=x)\leq Cn^{-{d}/{2}},
\end{equation}

\item[(b)]
%
%
\begin{equation} \label{eq:sumHPfirstDer}
\vert P(S_n = x) - P(S_n = y)\vert
\leq Cn^{-({d+1})/{2}},
\end{equation}

\item[(c)]
%
%
\begin{equation} \label{eq:sumHPsecondDer}
\vert P(S_n=x) -2P(S_n=y) + P(S_n=z) \vert\leq Cn^{-({d+2})/{2}}.
\end{equation}

\item[(d)]
In addition, for all $w$, $x$, $y$ and $z$ such that there exist $i\neq
j$ with
$x-y=w-z=e_i$ and $x-w=y-z=e_j$,
%
%
\begin{equation} \label{eq:sumHPsecondMixedDer}\qquad
\vert P(S_n=x)+P(S_n=z)-P(S_n=y)-P(S_n=w) \vert<Cn^{-({d+2})/{2}}.
\end{equation}
\end{longlist}
\end{claim}
\begin{pf}
Display (\ref{eq:sumHPupperBd}) is essentially a consequence of the
local limit theorem, see, for example, Theorem 2.3.8 in
Lawler and Limic \cite{LaLi}.
Indeed, if $EY_2 \in\Z^d$, that source yields\vadjust{\goodbreak} that for $S_n' := \sum
_{k=2}^{n+1} (Y_k - EY_k)$
and $\Gamma$ the covariance matrix of~$Y_2$,
there exists a constant $C$ such that
%
%
\begin{eqnarray} \label{eq:LCLT}
\qquad&&\vert P(S_n' = x) - p_n(x) \vert\nonumber\\[-8pt]\\[-8pt]
\qquad&&\qquad\leq Cn^{-({d+1})/{2}} \bigl( ( \Vert x \Vert_1^m n^{-
{m}/{2}} + 1 )
e^{-({ x^T \Gamma^{-1} x})/({2 n})} + n^{-({m-2})/{2}} \bigr)\nonumber
\end{eqnarray}
for all $n \in\N$ and $x \in\Z^d$,
where
\[
p_n(x) := \frac{1}{(2\pi n)^{{d}/{2}} \sqrt{ \det\Gamma}}
e^{-({x^T \Gamma^{-1}x})/({2n})}
\]
denotes the heat-kernel.
Equality (\ref{eq:LCLT}) in particular implies $P(S_n' = x) \leq Cn^{-
d/2}$, which entails
(\ref{eq:sumHPupperBd}).
If $EY_2 \notin\Z^d$, then as one may check by redoing the proof,
(\ref{eq:LCLT}) holds true for all $n \in\N$ and $x \in\Z^d - n EY_2$,
with $P(S_n' = x)$ replaced by $P(S_n' = x + nEY_2)$, which again
implies (\ref{eq:sumHPupperBd}).

Now in order to prove (\ref{eq:sumHPfirstDer}), note that
the triangle inequality yields
\begin{eqnarray*}
&&\vert P(S_{n+1} = x) - P(S_{n+1} = y) \vert\\
&&\qquad\leq
\max_{z_1,z_2 \dvtx\Vert z_1-z_2 \Vert_1=1} \vert P(S'_n = z_1) - P(S'_n
= z_2) \vert\\
&&\qquad\leq\max_{z_1,z_2 \dvtx\Vert z_1-z_2 \Vert_1=1}
\vert P(S'_n = z_1) - p_n(z_1) \vert+ \vert p_n(z_1) - p_n(z_2) \vert
\\
&&\qquad\quad{} + \vert p_n(z_2) - P(S'_n = z_2)\vert.
\end{eqnarray*}
Then (\ref{eq:LCLT}) in combination with standard heat kernel
estimates yields the desired result.

In a similar manner, (\ref{eq:sumHPsecondDer}) and (\ref
{eq:sumHPsecondMixedDer}) can be deduced from Theorem
2.3.8 of \cite{LaLi}, which we will omit for the sake of conciseness.
\end{pf}

Using a decomposition according to regeneration times, the previous
claim can be employed
to prove the following lemma.
%
\begin{lemma} \label{lem:annExitDistDerivatives}
For $L$ and $x \in\tilde{\mathcal{P}}(0,L)$, let $\nu_{x,L}$ denote
either $P_x(X_{T_{L^2}} \in\cdot)$,
$P_x(X_{T_{\partial\mathcal{P}(0,L)}} \in\cdot)$, $\mu_{x,0}^L$
or $P_x(X_{T_{L^2}} \in \cdot \vert (X_n-x) \cdot e_1 \geq0\
\forall n \in\N)$.

\begin{longlist}[(a)]
\item[(a)] \hypertarget{item:hittingProbBd}
There exists a constant $C$ such that for all $L$, $x \in\tilde
{{\mathcal{P}}}(0,L)$ and \mbox{$y \in H_{L^2}$},
%
%
\begin{equation}
\nu_{x,L}(y) \leq CL^{-d+1}.
\end{equation}

\item[(b)] \hypertarget{item:firstDer}
There exists a constant $C$ such that for all $L$,
$x \in\tilde{{\mathcal{P}}}(0,L)$,
$y \in H_{L^2}$ and $j \in\{2, \ldots, d\}$,
\[
\vert\nu_{x,L}(y)- \nu_{x,L}(y \pm e_j)\vert<CL^{{-d}}.
\]

\item[(c)] \hypertarget{item:firstStartPtDer}
There exists a constant $C$ such that for all $L$ and
$x,y \in\tilde{{\mathcal{P}}}(0,L)$ with $\Vert x-y \Vert_1$ as
well as
$z \in H_{L^2}$,
\[
\vert\nu_{x,L}(z)- \nu_{y,L}(z) \vert<CL^{{-d}}.
\]

\item[(d)] \hypertarget{item:secondDer}
There exists a constant $C$ such that for all $L$,
$x \in\tilde{{\mathcal{P}}}(0,L)$ and
$w,y,z \in H_{L^2}$ such that $\Vert w -y \Vert_1 = 1$ and $w-y = y-z$,
\[
\vert\nu_{x,L}(w)
- 2\nu_{x,L}(y) + \nu_{x,L}(z) \vert\leq CL^{-d-1}.
\]

\item[(e)] \hypertarget{item:mixedSecondDer}
There exists a constant $C$ such that for all $L$,
$x \in\tilde{{\mathcal{P}}}(0,L)$ and
$v,w,y,\allowbreak z \in H_{L^2}$ such that $\Vert v -w \Vert_1 = 1$, $z-y =w-v$
and $z-w = y-v$,
\[
\bigl\vert\nu_{x,L}(z) - \nu_{x,L}(y)
- \bigl(\nu_{x,L}(w)- \nu_{x,L}(v)\bigr) \bigr\vert\leq CL^{-d-1}.
\]
\end{longlist}
\end{lemma}
\begin{pf}
The fact that the particular choice among the first three possibilities
for $\nu_{x,L}$
is irrelevant, is a direct consequence of
Lemma \ref{lem:annealedPosExitProbEst}. With respect to the case that
$\nu_{x,L} = P_x(X_{T_{L^2}} \in \cdot \vert (X_n-x) \cdot
e_1 \geq0\ \forall n \in\N)$,
the desired result follows analogously from what comes below in combination
with Corollary 1.5 of \cite{SzZe-99}.

We will give the proof for $\nu_{x,L} = P_x(X_{T_{L^2}} \in\cdot)$.\vspace*{1pt}

For $k, l \in\N$ we define the event $B(l,k):= \{X_{\tau_k} \cdot
e_1 =l\}$ as well as
$B(l):=\bigcup_{k=1}^l B(l,k)$ and
\[
\hat{B}(l):=B(l)\cap\bigcap_{j=l+1}^{L^2-1}B^c(j);
\]
that is, for $l < L^2$ and on $\hat{B}(l)$, one has that $l$ is the
$e_1$-coordinate at which the last renewal before reaching the
$e_1$-coordinate $L^2$ occurs.

\hyperlink{item:hittingProbBd}{(a)}
We have
%
%
\begin{equation} \label{eq:relevantlDecompA}
P_x(X_{T_{L^2}}=y)
\leq P_x(A_L^c) + \sum_{l=L^2 - R_2(L)}^{L^2} F_l
\end{equation}
with
$
F_l := P_x( X_{T_{L^2}} = y, \hat{B}(l)),
$
and furthermore
%
%
\begin{eqnarray} \label{eq:FlRep}
F_l &=& \sum_{k=0}^{L^2} \sum_{z \in H_l}
P_x\bigl(X_{\tau_k} = z, X_{T_{L^2}} = y, \hat{B}(l)\bigr) \nonumber\\
&=& \sum_{k=0}^{L^2} \sum_{z \in H_l}
P_x(X_{\tau_k} = z) P_x\bigl(X_{T_{L^2}} = y, \hat{B}(l) \vert
X_{\tau_k} = z\bigr) \\
&=& \sum_{z \in H_l} P_x\bigl(X_{T_{L^2}} = y, \hat{B}(l) \vert
X_{\tau_1} = z\bigr)
\sum_{k=0}^{L^2} P_x(X_{\tau_k} = z).\nonumber
\end{eqnarray}
In order to estimate the inner sum of (\ref{eq:FlRep}),
set
$
m := E_0(X_{\tau_2} - X_{\tau_1})
$
and for $l \in\{L^2 - R_2(L), \ldots, L^2\}$ fixed, define
$
l^*:= \lfloor\frac{l}{m \cdot e_1} \rfloor.
$
We now distinguish cases.

First, assume $k \geq l^*$.
Then
$
\{X_{\tau_k} \cdot e_1 = l\} \subset H^1 \cup H^2,
$
where
$
H^1 := \{ X_{\tau_{{\lfloor{k}/{2} \rfloor}}} \cdot e_1 \leq
l/2 \}
$
and
$
H^2 := \{(X_{\tau_k} - X_{\tau_{{\lfloor{k}/{2} \rfloor}}})
\cdot e_1
\leq l/2 \}.
$
We get
\[
P_x(X_{\tau_k} = z, H^1)
= \sum_{y\dvtx y \cdot e_1 \leq l/2}
P_x(X_{\tau_k} = z \vert X_{\tau_{{\lfloor{k}/{2} \rfloor}}} = y)
P_x(X_{\tau_{{\lfloor{k}/{2} \rfloor}}} = y),
\]
and uniformly in $y$ and $z$ we have
$
P_x(X_{\tau_k} = z \vert X_{\tau_{{\lfloor{k}/{2} \rfloor
}}} = y)
\leq Ck^{-{d}/{2}}
$
due to the independence of the renewals (cf. Corollary 1.5 in \cite
{SzZe-99}) and (\ref{eq:sumHPupperBd}).
Now observe that using standard estimates for centred random variables
with finite $2d$th moment
[note that $(X_{\tau_2}-X_{\tau_1})\cdot e_1$ has finite $2d$th
moment as a~consequence of the assumption $(T)_\gamma$],
there exists a constant $C$ such that uniformly in $k$ and $L$ we have
%
%
\begin{equation} \label{eq:CLTRenewalEst}
P_x(X_{\tau_{{\lfloor{k}/{2} \rfloor}}} \cdot e_1 \leq l/2)
\leq1 \wedge
C {k^d}(k- l^* )^{-2d}.
\end{equation}
We therefore get
\begin{eqnarray*}
\sum_{k=l^*}^{L^2} P_x(X_{\tau_k} = z, H^1)
&\leq& C{l}^{-{d}/{2}} \Biggl( \sqrt{l^*} + \sum_{j=1}^\infty\bigl(l^*
+(j+1)\sqrt{l^*}\bigr)^d \bigl(j \sqrt{l^*}\bigr)^{-2d} \Biggr)\\
&\leq& Cl^{({-d+1})/{2}}
\end{eqnarray*}
and analogously for $H^2$, whence
%
%
\begin{equation} \label{eq:largekA}
\sum_{k=l^*}^{L^2} P_x(X_{\tau_k} = z) \leq Cl^{({-d+1})/{2}}.
\end{equation}
Now assume $k < l^*$. Then in the same manner as above we obtain
%
%
\begin{equation} \label{eq:smallkA}
\sum_{k=l^*/2}^{l^*} P_x(X_{\tau_k} = z) \leq Cl^{({-d+1})/{2}}
\end{equation}
and furthermore (\ref{eq:CLTRenewalEst}) supplies us with
%
%
\begin{equation} \label{eq:tinykA}
\sum_{k=0}^{l^*/2} P_x(X_{\tau_k} = z)
\leq C \sum_{k=0}^{l^*/2} k^d(l^* - k)^{-2d} \leq Cl^{-d+1}.
\end{equation}
In order to deal with the outer sum of (\ref{eq:FlRep}), note that for fixed
$l$ as well as $y^* \in H_{L^2}$ and $z^* \in H_l$ we have
%
%
\begin{eqnarray}
\label{eq:sumStartEndPt}\qquad\quad
\sum_{z \in H_l} P_x \bigl(X_{T_{L^2}} = y^*, \hat{B}(l) \vert
X_{\tau_1} = z\bigr)
&=& \sum_{y \in H_{L^2}} P_x \bigl(X_{T_{L^2}} = y, \hat{B}(l) \vert
X_{\tau_1} = z^*\bigr) \nonumber\\[-8pt]\\[-8pt]
&=& P_x\bigl(\hat{B}(l) \vert X_{\tau_1} = z^* \bigr).
\nonumber
\end{eqnarray}
Using (\ref{eq:largekA}) to (\ref{eq:tinykA}) in combination with
(\ref{eq:FlRep}), we therefore deduce
that for all $l \in\{L^2-R_2(L), \ldots, L^2\}$,
\[
F_l \leq C P_x(\hat{B}(l)) L^{-d+1}.\vadjust{\goodbreak}
\]
Thus, in combination with
(\ref{eq:relevantlDecompA}) and Lemma \ref{lem:ANEst} we get
\[
P_x(X_{T_{L^2}}=y) \leq CL^{-d+1},
\]
which finishes the proof.

\hyperlink{item:firstDer}{(b)}
We have
%
%
\begin{equation} \label{eq:relevantlDecomp}\qquad
\vert P_x(X_{T_{L^2}}=y)- P_x(X_{T_{L^2}}=y \pm e_j)\vert
\leq2P_x(A_L^c) + \sum_{l=L^2 - R_2(L)}^{L^2} F_l
\end{equation}
with $\hat{B}(l)$ as defined before and
\[
F_l := P_x\bigl( X_{T_{L^2}} = y, \hat{B}(l)\bigr) - P_x \bigl(X_{T_{L^2}} = y \pm
e_j, \hat{B}(l)\bigr).
\]
We compute
%
%
\begin{eqnarray}\label{eq:Fl2Rep}\hspace*{28pt}
F_l & = &\sum_{k=0}^{L^2} \sum_{z \in H_l}
\bigl( P_x\bigl(X_{\tau_k} = z, X_{T_{L^2}} = y, \hat{B}(l)\bigr) \nonumber\\
&&\hspace*{36.74pt}{} - P_x\bigl(X_{\tau_k} = z \pm e_j, X_{T_{L^2}} = y \pm e_j, \hat
{B}(l)\bigr) \bigr) \nonumber\\
& = &\sum_{k=0}^{L^2} \sum_{z \in H_l}
\bigl( P_x(X_{\tau_k} = z) P_x\bigl(X_{T_{L^2}} = y, \hat{B}(l) \vert
X_{\tau_k} = z\bigr) \nonumber\\[-8pt]\\[-8pt]
&&\hspace*{37pt}{} - P_x(X_{\tau_k} = z \pm e_j) P_x\bigl(X_{T_{L^2}} = y \pm e_j,
\hat{B}(l) \vert X_{\tau_k} = z \pm e_j\bigr) \bigr) \nonumber\\
& = &\sum_{z \in H_l}
P_x\bigl(X_{T_{L^2}} = y, \hat{B}(l) \vert X_{\tau_1} = z\bigr)\nonumber\\
&&\hspace*{16.84pt}{}\times\sum_{k=0}^{L^2}
\vert P_x(X_{\tau_k} = z) - P_x(X_{\tau_k} = z \pm e_j) \vert,\nonumber
\end{eqnarray}
where to obtain the last line we used the translation invariance of $\PP$.
Fix $l \in\{L^2 - R_2(L), \ldots, L^2\}$ and let $m$ and $l^*$ as before.
Again we distinguish cases.

First, assume $k \geq l^*$.
Then
$
\{X_{\tau_k} \cdot e_1 = l\} \subset H^1 \cup H^2,
$
where
$
H^1
$
and
$
H^2
$
as before.
Then
\begin{eqnarray*}
&&\vert P_x(X_{\tau_k} = z, H^1) - P_x(X_{\tau_k} = z \pm e_j, H^1)
\vert\\
&&\qquad= \sum_{y\dvtx y \cdot e_1 \leq l/2}
\bigl\vert P_x ( X_{\tau_k} = z \vert X_{\tau_{{\lfloor
{k}/{2} \rfloor}}} = y )
- P_x ( X_{\tau_k} = z \pm e_j \vert X_{\tau_{{\lfloor
{k}/{2} \rfloor}}} = y ) \bigr\vert\\
&&\hspace*{72.54pt}{} \times P_x ( X_{\tau_{{\lfloor{k}/{2} \rfloor}}} = y ),
\end{eqnarray*}
and uniformly in $y$, we have
\[
\bigl\vert P_x ( X_{\tau_k} = z \vert X_{\tau_{{\lfloor
{k}/{2} \rfloor}}} = y )
- P_x ( X_{\tau_k} = z \pm e_j \vert X_{\tau_{{\lfloor
{k}/{2} \rfloor}}} = y ) \bigr\vert
\leq Ck^{({-d-1})/{2}}
\]
due to the independence of the renewals and (\ref{eq:sumHPfirstDer}).
Using (\ref{eq:CLTRenewalEst}), we get
\begin{eqnarray*}
&&\sum_{k=l^*}^{L^2} \vert P_x(X_{\tau_k} = z, H^1) - P_x(X_{\tau_k}
= z \pm e_j, H^1) \vert\\
&&\qquad\leq C{l}^{({-d-1})/{2}} \Biggl( \sqrt{l^*} + \sum_{j=1}^\infty
\bigl(l^* +(j+1)\sqrt{l^*}\bigr)^d \bigl(j \sqrt{l^*}\bigr)^{-2d} \Biggr)
\leq Cl^{-{d}/{2}}
\end{eqnarray*}
and analogously for $H^2$, whence
%
%
\begin{equation} \label{eq:largek2}
\sum_{k=l^*}^{L^2} \vert P_x(X_{\tau_k} = z) - P_x(X_{\tau_k} = z
\pm e_j) \vert\leq Cl^{-{d}/{2}}.
\end{equation}
Now assume $k < l^*$. Then in the same manner we obtain
%
%
\begin{equation} \label{eq:smallk2}
\sum_{k=L^2/2}^{l^*} \vert P_x(X_{\tau_k} = z) - P_x(X_{\tau_k} = z
\pm e_j) \vert\leq Cl^{-{d}/{2}}
\end{equation}
and furthermore (\ref{eq:CLTRenewalEst}) supplies us with
%
%
\begin{eqnarray} \label{eq:tinyk2}
\sum_{k=0}^{l^*/2} \vert P_x(X_{\tau_k} = z) - P_x(X_{\tau_k} = z
\pm e_j) \vert
&\leq& C \sum_{k=0}^{l^*/2} k^d(l^* - k)^{-2d}\nonumber\\[-8pt]\\[-8pt]
&\leq& Cl^{-d+1}.\nonumber
\end{eqnarray}
Using (\ref{eq:sumStartEndPt}) and
(\ref{eq:Fl2Rep}) to (\ref{eq:tinyk2}), we deduce that there exists
$C$ such that for all $l \in\{L^2-R_2(L), \ldots, L^2\}$,
\[
F_l \leq C P_x(\hat{B}(l)) L^{-d}.
\]
In combination with
(\ref{eq:relevantlDecomp}) and Lemma \ref{lem:ANEst}, we get
\[
\vert P_x(X_{T_{L^2}}=y)- P_x(X_{T_{L^2}}=y \pm e_j)\vert\leq CL^{-d},
\]
which finishes the proof.

Parts \hyperlink{item:firstStartPtDer}{(c)}, \hyperlink
{item:secondDer}{(d)} and \hyperlink{item:mixedSecondDer}{(e)} follow
from analogous
calculations using (\ref{eq:sumHPfirstDer}), (\ref
{eq:sumHPsecondDer}) and (\ref{eq:sumHPsecondMixedDer}), respectively.
For the sake of conciseness, we omit giving the corresponding proofs here.
\end{pf}

To prove parts \hyperlink{item:quenchedAnnealedExpDiff}{(ii)} and
\hyperlink{item:quenchedAnnealedProbDiff}{(iii)} of Proposition \ref
{prop:closenessbase},
we quote and reprove a~conditional Azuma-type inequality appearing in
\cite{Be-09}.

In this context, denote by $(M_k)_{k \in\N_0}$ a one-dimensional martingale
on a~probability space $(\Omega, \mathcal{F}, P)$ with
filtration $(\mathcal{F}_k)_{k \in\N_0}$ and $M_0 = 0$.
Set $\Delta_k := M_k - M_{k-1}$ and assume that the
$\vert\Delta_k\vert$ are uniformly bounded from above by a finite constant.
Define for any nonnegative random variable $X$ its \textit{conditional
essential supremum} with respect to $\mathcal{F}_k$
as
$
\esssup(X \vert \mathcal{F}_{k-1}) := \lim_{n \to\infty}
E(X^n \vert \mathcal{F}_{k-1})^{1/n},
$
where the right-hand side exists due to Jensen's inequality. Set
\[
\sigma_k := \esssup( \vert\Delta_k \vert \vert \mathcal{F}_{k-1}).
\]
Then the \textit{essential variance} of the martingale is defined as
\[
V_k := \esssup\Biggl( \sum_{j=1}^k \sigma_j^2 \Biggr).
\]
%
\begin{lemma} \label{lem:Azuma}
If the $\Delta_k$ are uniformly bounded, then for all $n \in\N$ and
$t >0$,
\[
P(\vert M_n \vert> t) \leq2e^{-{t^2}/({2V_n})}.
\]
Furthermore, if $M_n = (M_n^{(1)}, \ldots, M_n^{(d)})$ with
$M_n^{(j)}$ being one-dimensional martingales
such that the differences $\Delta_k$ are uniformly bounded and with~$V_n^{(j)}$ as essential
variance, then writing $V_n^{\max}:= \max_{j \in\{1, \ldots, d\}}
V_n^{(j)}$ one has
\[
P(\Vert M_n \Vert_\infty> t) \leq2de^{-{t^2}/({2 V_n^{\max}})}.
\]
\end{lemma}
\begin{pf}
First, observe that
the $d$-dimensional case
is a direct consequence of the one-dimensional case by considering its
components and a~standard union bound.
It is therefore sufficient to prove the one-dimensional case.

We start with showing that for each $k \in\{1, \ldots, n\}$,
%
%
\begin{equation} \label{eq:azumaIndStatement}
E ( e^{\sum_{j=k}^n \Delta_j} \vert \mathcal{F}_{k-1}
)
\leq e^{(1/2) \esssup( \sum_{j=k}^n \sigma_j^2 \vert
\mathcal{F}_{k-1})}.
\end{equation}
To establish this inequality in the case $k=n$, we first of all note that
%
%
\begin{equation} \label{condEssSup}
\lim_{m \to\infty} E( \vert\Delta_n \vert^m \vert \mathcal
{F}_{n-1})^{{1}/{m}} \mathbh{1}_A
\geq\esssup(\vert\Delta_n \vert\mathbh{1}_A) \mathbh{1}_A -
\varepsilon
\end{equation}
for all
$x \in[0, \infty)$, $\varepsilon> 0$ and
\[
A:= A_{x,\varepsilon} := \Bigl\{ \lim_{m \to\infty}
E( \vert\Delta_n \vert^m \vert \mathcal{F}_{n-1})^{
{1}/{m}} \in(x,x+\varepsilon] \Bigr\} \in\mathcal{F}_{n-1}.
\]
We then observe that for such $A$ and with $C_A:=\esssup( \vert\Delta
_n \vert\mathbh{1}_A )$ as well as
\[
h_A \dvtx[-C_A, C_A] \ni s \mapsto\frac{e^{C_A} + e^{-C_A}}{2} + \frac
{e^{C_A} - e^{-C_A}}{2} \frac{s}{C_A},
\]
we obtain
\begin{eqnarray*}
E(e^{\Delta_n} \vert \mathcal{F}_{n-1})
&\leq& E(h_A(\Delta_n ) \vert \mathcal{F}_{n-1}) \mathbh{1}_A\\
&=& h_A (E(\Delta_n \mathbh{1}_A \vert \mathcal{F}_{n-1}))
\mathbh{1}_A\\
&=& h_A(0) \mathbh{1}_A = \frac{e^{C_A} + e^{-C_A}}{2} \mathbh{1}_A\\
&=&
\cosh(C_A) \mathbh{1}_A.
\end{eqnarray*}
Since by comparison of the corresponding power series one has $\cosh
(x) \leq e^{x^2/2}$, we obtain with (\ref{condEssSup}) that
%
%
\begin{eqnarray} \label{ineqCondA}
&&E(e^{\Delta_n} \mathbh{1}_A \vert
\mathcal{F}_{n-1})\nonumber\hspace*{-25pt}\\[-8pt]\\[-8pt]
&&\qquad \leq
e^{C_A^2/2} \mathbh{1}_A
\leq\exp\biggl\{\frac12 \Bigl(\lim_{m \to\infty}
E( \vert\Delta_n \vert^m \vert \mathcal{F}_{n-1})^{
{1}/{m}} + \varepsilon\Bigr)^2 \biggr\} \mathbh{1}_A.\nonumber\hspace*{-25pt}
\end{eqnarray}
Summing (\ref{ineqCondA}) over all $A:= A_{x,\varepsilon}$ for $x=j
\varepsilon$, $j \in\N_0$ we get
\begin{eqnarray*}
E(e^{\Delta_n} \vert \mathcal{F}_{n-1})
&\leq&\exp\biggl\{\frac12 \lim_{m \to\infty} E( \vert\Delta_n
\vert^m \vert \mathcal{F}_{n-1})^{{2}/{m}} \biggr\}\\
&&{} \times\exp\{ \esssup\vert\Delta_n \vert^2 \varepsilon+
\varepsilon^2/ 2\}.
\end{eqnarray*}
Since $\Delta_n$ was assumed to be bounded, taking $\varepsilon
\downarrow0$ yields (\ref{eq:azumaIndStatement})
for $k=n$.

Now we assume (\ref{eq:azumaIndStatement}) to hold true for $k+1$ and
deduce its validity for~$k$:
\begin{eqnarray*}
E (e^{\sum_{j=k}^n \Delta_j} | \mathcal{F}_{k-1} )
&=&
E ( e^{ \Delta_k} E ({e^{\sum_{j=k+1}^n \Delta_j}} |
\mathcal{F}_{k} ) | \mathcal{F}_{k-1} )\\
&\leq&
E \bigl(e^{ \Delta_k}e^{({1/2}) \esssup(\sum_{j=k+1}^{n} \sigma
_j^2 | \mathcal{F}_k )} | \mathcal{F}_{k-1} \bigr)\\
&\leq&
E \bigl(e^{ \Delta_k}e^{({1/2}) \esssup(\sum_{j=k+1}^{n} \sigma
_j^2 \vert \mathcal{F}_{k-1} )} \vert \mathcal{F}_{k-1}
\bigr)\\
&=&
e^{({1/2}) \esssup(\sum_{j=k+1}^{n} \sigma_j^2 \vert \mathcal
{F}_{k-1} )} E (e^{ \Delta_k} \vert \mathcal{F}_{k-1} )\\
&\leq&
e^{({1/2}) \esssup(\sum_{j=k+1}^{n} \sigma_j^2 \vert \mathcal
{F}_{k-1} )}e^{({1/2}) \sigma_k^2}\\
&=&
e^{({1/2}) \esssup(\sum_{j=k}^n \sigma_j^2 \vert \mathcal
{F}_{k-1} )},
\end{eqnarray*}
where to obtain the second inequality we used that for any nonnegative
random variable $X$ we have
\[
\esssup(X \vert \mathcal{F}_k )
\leq\esssup(X \vert \mathcal{F}_{k-1} ).
\]
Altogether, this establishes (\ref{eq:azumaIndStatement}).

Inserting $k=1$ in (\ref{eq:azumaIndStatement}), we deduce
$
E e^{\lambda M_n}
\leq e^{({1/2})\lambda^2 V_n }
$
for any real $\lambda$.
This estimate in combination with the exponential Chebyshev inequality yields
\begin{eqnarray*}
P(\vert M_n \vert> t) &=& P(M_n > t) + P(M_n < -t)\\
&\leq& e^{-\lambda t} (Ee^{\lambda M_n} + Ee^{-\lambda M_n})\\
&\leq& 2e^{-\lambda t} e^{({1/2}) \lambda^2 V_n }
\end{eqnarray*}
for $\lambda> 0$.
Setting $\lambda: = t/V_n$, this finishes the proof.
\end{pf}

The following result appears as Lemma 3.3 in Berger and Zeitouni \cite
{BeZe-08} and will prove
helpful in the following.
%
\begin{lemma} \label{lem:hittingProbBd}
Let $d \geq3$ and let $(v_n)_{n \in\N}$ be i.i.d., $\Z^d$-valued
random variables such that
$P$-a.s. we have $v_1 \cdot e_1 \geq1$ as well as
$
E \Vert v_1 \Vert^r < \infty
$
for some
$r \in[2,d-1]$.
Furthermore, assume that for some $\delta> 0$,
\[
P( v_1 \cdot e_1 = 1) > \delta,
\]
and that for all $z \in\Z^d$ of the form $z=e_1 \pm e_j$, $j \in\{2,
\ldots, d\}$, one has
\[
P(v_1=z \vert v_1 \cdot e_1 = 1) > \delta.
\]
Set $S_n:= \sum_{i=1}^n v_i$.
Then there exists a constant $K > 0$ such that
for all $z \in\Z^d$,
\[
P(\exists n \in\N\dvtx S_n = z) \leq K \vert z \cdot e_1 \vert
^{-r(d-1)/(r+d-1)}.
\]
Furthermore, for all $l \in\N$,
\[
\sum_{z \in H_l} P(\exists n \in\N\dvtx S_n = z) \leq1.
\]
\end{lemma}

The following result guarantees that with positive probability with
respect to the annealed measure,
the trajectories of two independent RWRE do never intersect.
%
\begin{lemma} \label{lem:nonIntersectionProbs}
Let $d \geq4$.
Then there exists $M \in(0,\infty)$ such that for $x_1, x_2 \in\Z
^d$ with $(x_1-x_2)\cdot e_1 = 0$
and $\Vert x_1-x_2 \Vert_\infty> M$ we have
%
%
\begin{equation} \label{eq:infStartPts}
\mathbf{P}_{x_1, x_2} \bigl(\bigl\{X^{(1)}_n \dvtx n \in\N\bigr\} \cap\bigl\{X^{(2)}_n
\dvtx n \in
\N\bigr\} = \varnothing\bigr) > 0,
\end{equation}
where
\[
\mathbf{P}_{x_1, x_2}:=
P_{x_1} \bigl( \cdot \vert X^{(1)}_n \cdot e_1 \geq x_1 \cdot e_1\
\forall n \in\N\bigr) \otimes P_{x_2}
\bigl( \cdot \vert X^{(2)}_n \cdot e_2 \geq x_2 \cdot e_2\ \forall
n \in\N\bigr)
\]
and
$X^{(1)}$ and $X^{(2)}$ denote copies of the RWRE $X$ ``driven'' by
$P_{x_1}$ and $P_{x_2}$, respectively.

In particular, for all $l \in\N$,
\[
\inf_{x_1,x_2 \in H_l, x_1 \not= x_2} P_{x_1} \otimes P_{x_2} \bigl(\bigl\{
X^{(1)}_n \dvtx n \in\N\bigr\} \cap\bigl\{X^{(2)}_n \dvtx n \in\N\bigr\} = \varnothing\bigr)
> 0
\]
also.
\end{lemma}
\begin{pf}
Due to uniform ellipticity, the last statement is a direct consequence
of (\ref{eq:infStartPts}). Thus,
we prove (\ref{eq:infStartPts}) now.

The proof is inspired by the proof of Proposition 3.4 in \cite{BeZe-08}.
The translation invariance of $\PP$ implies that we can assume $x \cdot
e_1 = x \cdot e_2 = 0$ without loss of generality.
Denote by $m := E_{x_1} ({X^{(1)}})^{*(2)}$ the expectation of the
second renewal radius and for $N \in\N$ set
\[
B^{(j)}_N:= \Biggl\{\sum_{k=1}^{N/(4m)} \bigl\Vert\bigl({X^{(j)}}\bigr)^{*(k)} \bigr\Vert
_1 \leq N/2 \Biggr\}.
\]
For $j \in\{1,2\}$,
with respect to $\mathbf{P}_{x_1,x_2} ( \cdot \vert A_N(X^{(j)}))$,
\[
\Biggl( \sum_{k=1}^n \bigl\Vert\bigl({X^{(j)}}\bigr)^{*(k)} \bigr\Vert_1
- \mathbf{E}_{x_1,x_2} \bigl( \bigl\Vert\bigl({X^{(j)}}\bigr)^{*(k)} \bigr\Vert_1 \vert
A_N\bigl(X^{(j)}\bigr) \bigr) \Biggr)_{n \in\{0, \ldots, 2N^2\}}
\]
is a martingale with bounded increments. Therefore, applying Azuma's
inequality for $N \in4m\N$ large enough results in
%
%
\begin{eqnarray}\qquad
&&\mathbf{P}_{x_1, x_2} \bigl( \bigl(B^{(j)}_N\bigr)^c
\bigr)\nonumber\\[-1pt]
&&\qquad= \mathbf{P}_{x_1, x_2} \Biggl( \sum_{k=1}^{N/(4m)} \bigl\Vert
\bigl({X^{(j)}}\bigr)^{*(k)} \bigr\Vert_1 > N/2 \Biggr) \nonumber\\[-1pt]
&&\qquad\leq\mathbf{P}_{x_1, x_2} \Biggl( \Biggl[ \sum_{k=1}^{N/(4m)} \bigl\Vert
\bigl({X^{(j)}}\bigr)^{*(k)} \bigr\Vert_1 \\[-1pt]
&&\qquad\quad\hspace*{36.3pt}{} - \mathbf{E}_{x_1,x_2} \bigl( \bigl\Vert\bigl({X^{(j)}}\bigr)^{*(k)} \bigr\Vert_1
\vert A_N\bigl(X^{(j)}\bigr) \bigr) > N/4
\Biggr] \bigg\vert A_N\bigl(X^{(j)}\bigr) \Biggr) \nonumber\\[-1pt]
&&\qquad\quad{} + \mathbf{P}_{x_1, x_2} \bigl( A_N\bigl(X^{(j)}\bigr)^c \bigr) \nonumber\\
\label{eq:earlyIntersectionEst}
&&\qquad\leq2\exp\biggl\{-\frac{(N/4)^2}{2NR^2_2(N)/4} \biggr\}
+ \mathbf{P}_{x_1, x_2} \bigl( A_N\bigl(X^{(j)}\bigr)^c \bigr);
\end{eqnarray}
here, we took advantage of
\[
m \geq\mathbf{E}_{x_1,x_2} \bigl( \bigl\Vert\bigl({X^{(j)}}\bigr)^{*(k)} \bigr\Vert_1
\vert A_N\bigl(X^{(j)}\bigr) \bigr)
\]
for all $k$ and $N \in\N$.

Furthermore, for $j \in\{1, 2\}$, $n \in\N$ and $\nu\in(0,1)$
define the
random times $h_{j,n} := \max\{ k \in\N_0 \dvtx X^{(j)}_{\tau_k} \cdot
e_1 \leq n \}$
as well as the event
\[
T^{(j)}_{\nu,N}:= \bigcap_{n \geq N/(4m)} \bigl\{
{\bigl(X^{(j)}\bigr)}^{*(h_{j,n+1})} \leq(2mn)^\nu\bigr\}.
\]
Then $(T)_\gamma$ implies that for any $\nu> 0$ and $K > 0$ there
exists a constant $C>0$ such that for all $N$ we have
%
%
\begin{equation} \label{eq:largeRenewalRadius}
\mathbf{P}_{x_1, x_2} \bigl( \bigl(T^{(j)}_{\nu,N}\bigr)^c \bigr) \leq CN^{-K}.
\end{equation}

Now we distinguish the situations in which the trajectories of the two
walks could intersect
in order to explain the decomposition
in~(\ref{eq:nonIntersectionDec1})\break
and~(\ref{eq:nonIntersectionDec2}) below; for this purpose,
assume that $x_1$ and $x_2$ from the assumptions satisfy
%
%
\begin{equation} \label{eq:distanceCond}
\Vert x_1 - x_2 \Vert_\infty\geq N^4.
\end{equation}

\begin{longlist}[(a)]
\item[(a)]
If the walks intersect within the first $N/(4m)$ renewal times of both walks,
then due to (\ref{eq:distanceCond})
this event is a subset of
$(B_N^{(1)})^c \cup(B_N^{(2)})^c$. This yields the first summand in
(\ref{eq:nonIntersectionDec1}).

\item[(b)]
Otherwise, the intersection may occur on $(T^{(1)}_{\nu, N})^c \cup
(T^{(2)}_{\nu, N})^c$, which
yields the second summand in (\ref{eq:nonIntersectionDec1}).

\item[(c)]
It remains to consider intersections after
$N/(4m)$ renewal times for at least one walk on $B_N^{(1)} \cap
B_N^{(2)} \cap T^{(1)}_{\nu, N} \cap T^{(2)}_{\nu, N}$; note that due
to the restriction to $B_N^{(1)} \cap B_N^{(2)}$ and (\ref{eq:distanceCond}),
the intersection can take place in $H_n$
with $n \geq N/(4m)$ only.
In this case, since we restrict to $T^{(1)}_{\nu, N} \cap T^{(2)}_{\nu
, N}$,
if the trajectories intersect in the hyperplane~$H_n$,
there must have occurred a renewal for each of the walks in distance at
most $(2mn)^\nu$
from the point of intersection which implies that the two renewals must
occur at sites that have distance $2(2mn)^\nu$ at most
from each other.
Thus, (\ref{eq:nonIntersectionDec2}) corresponds to an intersection
after at least $N/(4m)$ renewals
for at least one walk,
on $B_N^{(1)} \cap B_N^{(2)} \cap T^{(1)}_{\nu, N} \cap T^{(2)}_{\nu, N}$.
\end{longlist}

Consequently, choosing $\nu> 0$ small enough,
we obtain using Lemma \ref{lem:hittingProbBd} with $r=2$ as well as
(\ref{eq:earlyIntersectionEst}) and (\ref{eq:largeRenewalRadius}), that
%
%
\begin{eqnarray}
&&\mathbf{P}_{x_1, x_2} \bigl( \bigl\{X^{(1)}_n \dvtx n \in\N\bigr\} \cap\bigl\{X^{(2)}_n
\dvtx
n \in\N\bigr\} \not= \varnothing\bigr) \nonumber\\[-2pt]
\label{eq:nonIntersectionDec1}
&&\qquad\leq2 \mathbf{P}_{x_1, x_2} \bigl( \bigl(B_N^{(1)}\bigr)^c \bigr)
+ 2\mathbf{P}_{x_1, x_2} \bigl( \bigl(T^{(1)}_\nu\bigr)^c \bigr) \\[-2pt]
\label{eq:nonIntersectionDec2}
&&\qquad\quad{} + \sum_{j \geq N/(4m)} \sum_{z \in H_{j}} \sum
_{z' \dvtx\Vert z - z'\Vert_1 \leq2(2mj)^\nu}
\mathbf{P}_{x_1, x_2} \bigl(\exists i \dvtx X^{(1)}_{\tau_i} =
z\bigr)\nonumber\\[-9pt]\\[-9pt]
&&\qquad\quad\hspace*{149.5pt}{}\times\mathbf{P} _{x_1, x_2}\bigl(\exists k \dvtx X^{(2)}_{\tau_k} = z'\bigr) \nonumber\\[-2pt]
&&\qquad\leq C \biggl( N^{-K} + \sum_{j \geq N/(4m)} \sum_{z \in H_{j}}
\mathbf{P}_{x_1, x_2} \bigl(\exists i \dvtx X^{(1)}_{\tau_i} = z\bigr) \\[-2pt]
&&\qquad\quad\hspace*{109pt}\hspace*{-62.2pt}{} \times\sum_{z' \dvtx\Vert z - z'\Vert_1 \leq2(2mj)^\nu}
\bigl(j - 2(2mj)^\nu\bigr)^{-2(d-1)/(d+1)} \biggr) \nonumber\\[-2pt]
\label{eq:probConvZero}
&&\qquad\leq C \bigl(N^{-K} + N^{\nu d} N^{-2(d-1)/(d+1) + 1} \bigr) \to0
\end{eqnarray}
as $N \to\infty$, and where for ease of notation we omitted to
emphasize that the renewal times $\tau_i$ refer to the process
that is evaluated at these times.
Choosing $M = N^4$ for some $N$ such that the term in (\ref
{eq:probConvZero}) is smaller than~1, this establishes the
lemma.\vspace*{-3pt}
\end{pf}

For $\omega\in\Omega$ and $z \in\Z^d$ we set
$
\mathbf{P}_{z,\omega} := P_{z,\omega} \otimes P_{z, \omega}
$
as well as
$
\mathbf{P}_z := \int_\Omega P_{z,\omega} \otimes P_{z, \omega} \PP
(d\omega),
$
where the RWRE ``driven'' by the first factor is denoted by $X^{(1)}$
and the one
driven by the second factor is denoted by $X^{(2)}$.
Using the previous lemma, we can bound the number of intersections of
two independent RWREs as follows.\vadjust{\goodbreak}
%
\begin{lemma}\label{lem:manyIntersectionsProbUpperBd}
There exists a positive constant $C$ such that for all $L$ large enough as
well as $z \in\tilde{\mathcal{P}}(0,L)$ and $m \in\N$,
\begin{eqnarray*}
\hspace*{-4pt}&&\mathbf{P}_{z} \bigl(
\bigl\vert
\bigl\{X_n^{(1)}\dvtx n \in\N\bigr\} \cap\bigl\{X_n^{(2)}\dvtx n \in\N\bigr\} \cap{\mathcal
{P}}(0,L) \bigr\vert>mR_2^{d+1}(L) \vert
A_L\bigl(X^{(1)}\bigr),\\
&&\qquad A_L\bigl(X^{(2)}\bigr)
\bigr)
 <e^{-Cm}.
\end{eqnarray*}
\end{lemma}
\begin{pf}
For $L$ large enough, any $k$ such that $k+R_2(L)<L$ and $j \in\{1, 2\}
$ we have
%
%
\begin{equation}\label{eq:inslab}
\qquad\qquad\mathbh{1}_{A_L(X^{(j)})} \cdot\bigl\vert\bigl\{x\in\bigl\{X^{(j)}_n \dvtx n
\in\N\bigr\}
\dvtx k< x \cdot e_1 <k+R_2(L) \bigr\} \bigr\vert< R_2^{d+1} (L)
\end{equation}
as well as
%
%
\begin{equation} \label{eq:inLeftPartSlab}
\mathbh{1}_{A_L(X^{(j)})} \cdot\bigl\vert\bigl\{x \in \bigl\{X^{(j)}_n \dvtx
n \in\N\bigr\}\dvtx x \cdot e_1 \leq0\bigr\} \bigr\vert< R_2^{d+1} (L).
\end{equation}

For every $k$, let $Q^-_k :={\mathcal{P}}(0,L)\cap\{x \dvtx x \cdot
e_1 < k R_2(L)\}$ and
$Q^+_k:={\mathcal{P}}(0,L) \cap\{x \dvtx x \cdot e_1\geq kR_2(L) \}$.
Due to Lemma \ref{lem:nonIntersectionProbs}, we can infer that there
exists $\rho>0$ such that for every $k$
and uniformly in $z \in\tilde{\mathcal{P}}(0,L)$,
%
%
\begin{eqnarray}\label{eq:noIntersBoxPosProb}
&&\mathbf{P}_{z} \bigl(\bigl\{ \bigl\{X_n^{(1)}\dvtx n \in\N\bigr\} \cap\bigl\{X_n^{(2)}\dvtx n
\in\N\bigr\} \cap
Q^+_{k+1}=\varnothing
\bigr\}
\vert\nonumber\\
&&\qquad A_L\bigl(X^{(1)}\bigr), A_L\bigl(X^{(2)}\bigr), \bigl\{X_n^{(1)}\dvtx n \in\N\bigr\} \cap
Q^-_{k},\\
&&\hspace*{129.64pt} \bigl\{X_n^{(2)}\dvtx n \in\N\bigr\} \cap Q^-_{k}
\bigr) >\rho.
\nonumber
\end{eqnarray}
Let
\[
J^{(\mathrm{even})}:=\bigl\{k \in2\N_0 \dvtx
\bigl\{X_n^{(1)}\dvtx n \in\N\bigr\} \cap\bigl\{X_n^{(2)}\dvtx n \in\N\bigr\} \cap
Q^+_{k}\cap
Q^-_{k+1}\neq\varnothing
\bigr\}
\]
and
\[
J^{(\mathrm{odd})}:=\bigl\{k \in2\N_0 +1\dvtx
\bigl\{X_n^{(1)}\dvtx n \in\N\bigr\} \cap\bigl\{X_n^{(2)}\dvtx n \in\N\bigr\} \cap
Q^+_{k}\cap
Q^-_{k+1}\neq\varnothing
\bigr\}.
\]

Then by (\ref{eq:noIntersBoxPosProb}), conditioned on $A_L(X^{(1)})
\cap A_L(X^{(2)})$, both~$J^{(\mathrm{even})}$\break and~$J^{(\mathrm
{odd})}$ are
stochastically dominated
by a geometric variable with parameter~$\rho$.

The lemma now follows when we remember that by (\ref{eq:inslab}) and
(\ref{eq:inLeftPartSlab}),
\begin{eqnarray*}
&&\mathbh{1}_{A_L(X^{(1)})} \mathbh{1}_{A_L(X^{(2)})} \cdot
\bigl\vert\bigl\{X_n^{(1)}\dvtx n \in\N\bigr\} \cap\bigl\{X_n^{(2)}\dvtx n \in\N\bigr\} \cap
{\mathcal{P}}(0,L) \bigr\vert\\
&&\qquad\leq R_2^{d+1}(L) \bigl(J^{(\mathrm{even})}+J^{(\mathrm{odd})} \bigr).
\end{eqnarray*}
\upqed\end{pf}

As a corollary of Lemma \ref{lem:manyIntersectionsProbUpperBd}, we
obtain the following estimate.
%
\begin{corollary}\label{cor:intersectionNumberExpUpperBd}
With the same notation as in Lemma \ref{lem:manyIntersectionsProbUpperBd},
%
%
\begin{eqnarray}
\label{eq:interquenched}\qquad
&&\PP\bigl( \exists z \in\tilde{\mathcal{P}}(0,L) \dvtx\nonumber\\
&&\qquad\mathbf{E}_{z,\omega} \bigl( \bigl\vert
\bigl\{X_n^{(1)}\dvtx n \in\N\bigr\} \\
&&\qquad\hspace*{113.5pt}\hspace*{-87.8pt}{}\cap\bigl\{X_n^{(2)}\dvtx n \in\N\bigr\} \cap{\mathcal
{P}}(0,L)\bigr\vert
\vert A_L\bigl(X^{(1)}\bigr), A_L\bigl(X^{(2)}\bigr) \bigr)
\geq R_3(L) \bigr)\nonumber
\end{eqnarray}
is contained in $\mathcal{S}(\N)$ as a function in $L$.
\end{corollary}
\begin{pf}
Set $Z:= \vert\{X_n^{(1)}\dvtx n \in\N\} \cap\{X_n^{(2)}\dvtx n \in\N\}
\cap{\mathcal{P}}(0,L) \vert$ and note
that on $A_L(X^{(1)})\cap A_L(X^{(2)})$, the variable $Z$ is bounded
from above by $\vert{\mathcal{P}}(0,L) \vert\leq(2L^2)^d$. Thus,
\begin{eqnarray*}
&&\PP\bigl( \mathbf{E}_{z,\omega} Z \geq R_3(L) \vert A_L\bigl(X^{(1)}\bigr),
A_L\bigl(X^{(2)} \bigr)\bigr)\\[-2pt]
&&\qquad\leq\PP\bigl( \mathbf{E}_{z,\omega} \bigl(Z, Z \geq n R_2^{d+1}(L)
\vert A_L\bigl(X^{(1)}\bigr), A_L\bigl(X^{(2)}\bigr) \bigr) \geq R_3(L)/2 \bigr)\\[-2pt]
&&\qquad\quad{} +\, \underbrace{\PP\bigl( \mathbf{E}_{z,\omega} \bigl( Z, Z \leq n
R_2^{d+1}(L) \vert A_L\bigl(X^{(1)}\bigr), A_L\bigl(X^{(2)}\bigr) \bigr)
\geq R_3(L)/2 \bigr) }_{= 0\ \mathrm{for}\ n = R_2(L)}\\[-2pt]
&&\qquad\leq(2L^2)^d \mathbf{P}_{z} \bigl( Z \geq nR_2^{d+1}(L) \vert
A_L\bigl(X^{(1)}\bigr), A_L\bigl(X^{(2)}\bigr) \bigr) \leq e^{-CR_2(L)}
\end{eqnarray*}
for $n = R_2(L)$ and $L$ large enough due to Lemma \ref
{lem:manyIntersectionsProbUpperBd}.
Taking the union bound for $z \in\tilde{\mathcal{P}}(0,L)$ finishes
the proof.
\end{pf}

We define $J(L) \subset\Omega$ to be the set of all $\omega$ such that
for every $z \in\tilde{\mathcal{P}}(0,L)$,
\[
\mathbf{E}_{z,\omega}\bigl( \bigl\vert
\bigl\{X_n^{(1)}\dvtx n \in\N\bigr\} \cap\bigl\{X_n^{(2)}\dvtx n \in\N\bigr\} \cap{\mathcal
{P}}(0,L) \bigr\vert \vert
A_L\bigl(X^{(1)}\bigr), A_L\bigl(X^{(2)}\bigr) \bigr) \leq R_3(L).
\]
Then by Corollary \ref{cor:intersectionNumberExpUpperBd},
%
%
\begin{equation} \label{eq:JNProb}
\PP(J(\cdot)^c) \in\mathcal{S}(\N),
\end{equation}
and for $\omega\in J(L)$ and $z \in\tilde{\mathcal{P}}(0,L)$,
%
%
\begin{equation} \label{eq:JHittingProbBd}
\sum_{x\in{\mathcal{P}}(0,L)} P_{z,\omega}(T_x<
\infty)^2<R_3(L).\vspace*{-3pt}
\end{equation}

\subsection{\texorpdfstring{Proof of Proposition
    \protect\ref{prop:closenessbase}\protect\hyperlink{item:quenchedAnnealedExpDiff}{(ii)}}
{Proof of Proposition 3.4(ii)}}

The following lemma will yield part~\hyperlink
{item:quenchedAnnealedExpDiff}{(ii)} of Proposition \ref{prop:closenessbase}.\vspace*{-3pt}
%
\begin{lemma}\label{lem:quenchedAnnealedExpDiff}
There exists a sequence of events $G^{\mbox{\hyperlink
{item:quenchedAnnealedExpDiff}{\textup{\ii}}}}_L \subset\Omega$
such that
\[
\PP\bigl( {G^{\mbox{\hyperlink{item:quenchedAnnealedExpDiff}{\textup{\ii}}}}_{\cdot
}}^c \bigr) \in\mathcal{S}(\N)
\]
and for every $\omega\in G^{\mbox{\hyperlink
{item:quenchedAnnealedExpDiff}{\textup{\ii}}}}_L$ and $z\in\tilde{{\mathcal
{P}}}(0,L)$,
\begin{eqnarray*}
&&\bigl\Vert
E_{z,\omega} \bigl( X_{T_{\partial{\mathcal{P}}(0,L)}}
\vert T_{\partial{\mathcal{P}}(0,L)} = T_{\partial_+ {\mathcal
{P}}(0,L)} \bigr)
-E_z \bigl( X_{T_{\partial{\mathcal{P}}(0,L)}}
\vert T_{\partial{\mathcal{P}}(0,L)} = T_{\partial_+ {\mathcal
{P}}(0,L)} \bigr)
\bigr\Vert_1\\[-2pt]
&&\qquad\leq R_4(L).\vspace*{-3pt}
\end{eqnarray*}
\end{lemma}
\begin{pf}
As a consequence of Lemma \ref{lem:ANEst}, Proposition \ref
{prop:closenessbase}\hyperlink{item:quenchedBadExitEst}{(i)} and (\ref
{eq:JNProb}),
it is sufficient to show that denoting
\[
U(\omega,z):=
\Vert
E_{z,\omega}
( X_{T_{L^2}}, A_L, J(L) )
-E_z ( X_{T_{L^2}}, A_L, J(L) )
\Vert_1,
\]
one has that
%
%
\begin{eqnarray} \label{eq:UupperBd}
\PP\biggl(
\bigcup_{z \in\tilde{\mathcal{P}}(0,L)}
\{\omega\dvtx U(\omega,z) > R_4(L)/2\}\biggr)\nonumber\\[-9pt]\\[-9pt]
&&\eqntext{\mbox{is contained
in } \mathcal{S}(\N) \mbox{ as a function in } L.}\vadjust{\goodbreak}
\end{eqnarray}
To this end, note that on $A_L$ the walk starting in $\tilde{\mathcal
{P}}(0,L)$ can visit sites in
%
%
\begin{eqnarray} 
S_L &:=& \{ x \in\Z^d \dvtx-R_2(L) -L^2/3\leq x \cdot e_1 <
L^2,\nonumber\\[-8pt]\\[-8pt]
&&\hspace*{67.4pt}\Vert\pi_{e_1^\bot}(x)\Vert_\infty\leq2L^2R_2(L)
\}\nonumber
\end{eqnarray}
only before hitting $H_{L^2}$. Order the vertices contained in $S_L$
lexicographically,
that is in increasing order of their first differing coordinate, as
$x_1,x_2,\ldots, x_n$.
Let $\mathcal{G}_0 := \{\Omega, \varnothing\}$ and for
$k \in\{1, \ldots, n\}$, let $\mathcal{G}_k$ be the $\sigma
$-algebra on
$\Omega$ that is generated by
$(\omega(x_j))_{j \in\{1, \ldots, k \}}$. Furthermore, define the martingale
\[
M_k:= E_{z} ( X_{T_{L^2}}, A_L, J(L) \vert \mathcal{G}_k
).
\]
Note that due to the independence structure of $\PP$,
taking the conditional expectation with respect to $\mathcal{G}_k$
is nothing else than taking the expectation with respect to the process
as well as over all those $\omega(x)$ for which
$x \notin\{x_1, \ldots, x_k\}$.
Thus,\vspace*{1pt}
$
E_{z,\omega} (X_{T_{L^2}}, A_L, J(L)) = E_z (X_{T_{L^2}}, A_L, J(L)
\vert \mathcal{G}_n)(\omega)
$
for $\PP$-a.a. $\omega$ as well as
$
E_z(X_{T_{L^2}}, A_L, J(L)) = E_z(X_{T_{L^2}}, A_L, J(L) \vert
\mathcal{G}_0)$.\vspace*{1pt}

Next, we estimate $\esssup(\Vert M_k-M_{k-1} \Vert_1 | \mathcal
{G}_{k-1})$
similarly to \cite{Be-09}, which again is based on ideas from
Bolthausen and Sznitman \cite{BoSz-02b}.
For $x \in\Z^d$, let
\[
B(x):= \{ y \in H_{x\cdot e_1 - 1} \dvtx\Vert x-y \Vert_1 \leq
R_2(L)+1 \}.
\]
Note that if $x$ is visited, then on $A_L$ the first visit to the
affine hyperplane~$H_{x \cdot e_1 -1}$
will occur at a point contained in $B(x)$.
Therefore,
%
%
\begin{eqnarray} \label{eq:UkEst}
U_k:\!&=&
\esssup( \Vert M_k-M_{k-1} \Vert_1 \vert \mathcal
{G}_{k-1})\nonumber\\
&=&\esssup\bigl(
\bigl\Vert E_{z} \bigl( X_{T_{L^2}}, A_L, J(L),T_{x_k} < \infty
\vert \mathcal{G}_k \bigr)\nonumber\\
&&\hspace*{37.1pt}{} - E_{z} \bigl( X_{T_{L^2}}, A_L, J(L), T_{x_k} < \infty \vert
\mathcal{G}_{k-1} \bigr) \bigr\Vert_1
\vert \mathcal{G}_{k-1} \bigr)\nonumber\\
&\leq&
R_2^2(L)P_z\bigl(T_{x_k} <\infty, A_L, J(L) \vert \mathcal
{G}_{k-1}\bigr)\nonumber\\[-8pt]\\[-8pt]
&\leq&
R_2^2(L) \sum_{y\in B(x_k) \cap S_L} P_z\bigl(X_{T_{y \cdot e_1}}=y, J(L)
\vert \mathcal{G}_{k-1}\bigr)\nonumber\\
&=&
R_2^2(L) \sum_{y\in B(x_k) \cap S_L} P_{z,\omega}\bigl(X_{T_{y\cdot
e_1}}=y, J(L)\bigr) \nonumber\\
&\leq&
R_2^2(L) \biggl(\sum_{y\in B(x_k) \cap\mathcal{P}(0,L)}
P_{z,\omega}\bigl(T_{y} <\infty, J(L)\bigr)
\nonumber\\
&&\hspace*{57.8pt}{}+ P_{z,\omega}\bigl(T_{\partial_+ \mathcal{P}(0,L)} \not= T_{\partial
\mathcal{P}(0,L)}\bigr) \biggr).
\nonumber
\end{eqnarray}
Here, the first equality follows since
\[
E_{z} \bigl( X_{T_{L^2}}, A_L, J(L), T_{x_k} = \infty \vert
\mathcal{G}_k \bigr)
- E_{z} \bigl( X_{T_{L^2}}, A_L, J(L), T_{x_k} = \infty \vert
\mathcal{G}_{k-1} \bigr) = 0,
\]
which is
due to the fact that the restriction to $T_{x_k} =\infty$ makes the
inner random variables
independent of the realization\vadjust{\goodbreak} of $\omega(x_k)$.
The first inequality follows since, if the walker hits $x_k$, then on
$A_L$ the site of
the subsequent renewal has distance at most $R_2(L)$ to $x_k$
and consequently, using standard coupling arguments, one obtains that
the values of
\[
E_{z} \bigl( X_{T_{L^2}}, A_L, J(L), T_{x_k} < \infty
\vert \mathcal{G}_k \bigr)
\]
as a function in $\omega(x_k)$ (and all other coordinates fixed)
lie within distance of $R_2^2(L)$ of each other.

Now for $\omega\in G^{\mbox{\hyperlink{item:quenchedBadExitEst}{\iitem}}}_L
\cap J(L)$,
remembering that $\vert B(x_k) \vert\leq(3R_2(L))^d$ and that every
$y$ is in $B(x)$ for at most
$(3R_2(L))^d$ different points $x$, using (\ref{eq:UkEst}) we infer
\begin{eqnarray*}
\sum_{k=1}^n U_k^2
&\leq&
\sum_{k=1}^n R_2^4(L) \biggl(\sum_{y\in B(x_k) \cap\mathcal{P}(0,L)}
P_{z,\omega}\bigl(T_{y} <\infty, J(L)\bigr)\\[-2pt]
&&\hspace*{74.7pt}{} + P_{z,\omega}\bigl(T_{\partial_+ \mathcal{P}(0,L)} \not= T_{\partial
\mathcal{P}(0,L)}\bigr) \biggr)^2\\[-2pt]
&\leq&
2(3R_2(L))^{2d}R_2^{4}(L)\sum_{k=1}^n \biggl(
\sum_{y\in B(x_k) \cap\mathcal{P}(0,L)}
P_{z,\omega}\bigl(T_{y} <\infty, J(L)\bigr)^2 \\[-2pt]
&&\hspace*{131.5pt}{}
+ P_{z,\omega}\bigl(T_{\partial_+ \mathcal{P}(0,L)} \not= T_{\partial
\mathcal{P}(0,L)}\bigr)^2
\biggr)\\[-2pt]
&\leq&
2(3R_2(L))^{4d} R_2^{4}(L) \biggl( \sum_{y\in\mathcal{P}(0,L)}
P_{z,\omega}\bigl(T_y < \infty, J(L)\bigr)^2\\[-2pt]
&&\hspace*{89.7pt}{}
+ P_{z,\omega}\bigl(T_{\partial_+ \mathcal{P}(0,L)} \not= T_{\partial
\mathcal{P}(0,L)}\bigr)^2 \biggr)\\[-2pt]
&\leq&4(3R_2(L))^{4d} R_2^{4}(L) R_3(L) \leq R_3^2(L)
\end{eqnarray*}
for $L$ large enough, where the fourth inequality is a consequence of
(\ref{eq:JHittingProbBd}) and
part \hyperlink{item:quenchedBadExitEst}{(i)} of Proposition \ref
{prop:closenessbase}.

Therefore, by Lemma \ref{lem:Azuma} applied to the $(d-1)$-dimensional
martingale~$M_k$,
\[
\PP\bigl(U(\omega,z) > R_4(L)/2\bigr)
<2de^{-{R_4^2(L)}/({8R_3^2(L)})} + \PP(J(L)^c) + \PP\bigl(
{G^{\mbox{\hyperlink{item:quenchedBadExitEst}{\iitem}}}_L}^c \bigr)
\]
and the right-hand side
is contained in $\mathcal{S}(\N)$ as a function in $L$ due to
Proposition \ref{prop:closenessbase}\hyperlink
{item:quenchedBadExitEst}{(i)}, Lemma \ref{lem:ANEst} and (\ref{eq:JNProb}).
Now since the above estimates and hence the last inclusion were uniform
in $z \in\tilde{\mathcal{P}}(0,L)$,
we infer that (\ref{eq:UupperBd}) holds,
which in
combination with Remark \ref{rem:PropProofStrat} finishes the proof.
\end{pf}

\subsection{\texorpdfstring{Auxiliary results for the proof of Proposition
    \protect\ref{prop:closenessbase}\protect\hyperlink{item:quenchedAnnealedProbDiff}{(iii)}}
{Auxiliary results for the proof of Proposition 3.4(iii)}}

The following lemma is the basis for proving Proposition \protect\ref
{prop:closenessbase}\protect\hyperlink
{item:quenchedAnnealedProbDiff}{(iii)}.

\begin{lemma} \label{lem:quenchedAnnealedLargeBoxes}
Fix $\vartheta\in(\frac{d-1}{d},1]$, let $C$ be a constant and
denote by $B^{\vartheta}(L)$ the set of those $\omega$ for which for\vadjust{\goodbreak}
all $M \in\{{\lfloor\frac25L^2 \rfloor}, \ldots, L^2\}$,
all $z \in\tilde{\mathcal{P}}(0,L)$ and all $(d-1)$-dimensional
hypercubes $Q$ of side length $\lceil L^{\vartheta} \rceil$ that are
contained in $H_M$,
one has
%
%
\begin{equation} \label{eq:quenchedAnnealedLargeBoxes}
\vert P_{z,\omega} (X_{T_M}\in Q)
-P_z (X_{T_M}\in Q)\vert
\leq CL^{(\vartheta-1)(d-1)}.
\end{equation}
Then for $C$ large enough,
$
\PP(B^{\vartheta}(L)^c)
$
\mbox{is contained in $\mathcal{S}(\N)$ as a function~in~$L$.}
\end{lemma}
\begin{pf}
Choose $\vartheta' \in(\frac{d-1}{d},\vartheta)$, set
$
U:={\lfloor L^{2\vartheta'} \rfloor},
$
fix
$
M \in\{{\lfloor\frac25 L^2 \rfloor}, \ldots, L^2\}
$
and
with $S_L$ as in the previous proof
set
$
S_L^M
:=S_L \cap\{x \in\Z^d \dvtx x \cdot e_1 \leq M\}.
$
Similarly to the proof of Lemma \ref{lem:quenchedAnnealedExpDiff}, we
let $x_1, x_2, \ldots, x_n$ be a lexicographic ordering of the
vertices in
$
S^M_L
$
and denote by $\mathcal{G}_k$ the $\sigma$-algebra on $\Omega$
generated by $\omega(x_1), \ldots, \omega(x_k)$.
For $v\in H_{M+U}$, we start with estimating
\[
\bigl\vert P_{z}\bigl(X_{T_{M+U}}=v, A_L, J(L) \vert \mathcal{G}_n \bigr)
- P_z\bigl(X_{T_{M+U}} = v, A_L, J(L)\bigr) \bigr\vert,
\]
and for this purpose consider the martingale
\[
M_k:= P_z\bigl(X_{T_{M+U}}=v, A_L, J(L) \vert \mathcal{G}_k \bigr).
\]
As in the proof of Lemma \ref{lem:quenchedAnnealedExpDiff},
from which we borrow the notation $B(x_k)$,
we are going to take advantage of Lemma \ref{lem:Azuma}, whence we
will need to bound
$
\Delta_k
:=\esssup( \vert M_{k}-M_{k-1} \vert \vert \mathcal{G}_{k-1}).
$
By part \hyperlink{item:firstStartPtDer}{(c)} of Lemma \ref
{lem:annExitDistDerivatives},
again in combination with standard coupling arguments, we obtain
\begin{eqnarray*}
\Delta_k
&\leq& C R_2^2(L) P_z\bigl( T_{x_k} < \infty, A_L, J(L) \vert \mathcal
{G}_{k-1} \bigr) \cdot U^{-d/2}\\
&\leq& CU^{-d/2} R_2^2(L) \sum_{y\in B(x_k) \cap S_L^M} P_{z,\omega
}\bigl(X_{T_{y\cdot e_1}}=y, A_L, J(L)\bigr)\\
&\leq&
C U^{-d/2} R_2^2(L)
\biggl(\sum_{y\in B(x_k) \cap\mathcal{P}(0,L)}
P_{z,\omega}\bigl(T_{y} <\infty, J(L)\bigr) \\
&&\hspace*{96pt}{}
+ P_{z,\omega}\bigl(T_{\partial_+ \mathcal{P}(0,L)} \not= T_{\partial
\mathcal{P}(0,L)}\bigr) \biggr).
\end{eqnarray*}
Therefore, for $\omega\in J(L) \cap G^{\mbox{\hyperlink
{item:quenchedBadExitEst}{\iitem}}}_L$,
and based on the same calculation as in the proof of Lemma \ref
{lem:quenchedAnnealedExpDiff},
%
%
\begin{equation} \label{eq:differenceBound2}
\esssup\Biggl( \sum_{k=1}^{n} \Delta_k^2 \Biggr)
\leq R_4(L)U^{-d}.
\end{equation}
Indeed, continuing the previous chain, for $\omega\in J(L) \cap
G^{\mbox{\hyperlink{item:quenchedBadExitEst}{\iitem}}}_L$ we have
\begin{eqnarray*}
\sum_{k=1}^{n} \Delta_k^2
&\leq&
CU^{-d} R_2^4(L)
\sum_{k=1}^{n} \biggl(\sum_{y\in B(x_k) \cap\mathcal{P}(0,L)}
P_{z,\omega}\bigl(T_{y} <\infty, J(L)\bigr)\\
&&\hspace*{103.1pt}{} + P_{z,\omega}\bigl(T_{\partial_+ \mathcal{P}(0,L)} \not=
T_{\partial\mathcal{P}(0,L)}\bigr) \biggr)^2 \\
&\leq&
CU^{-d} R_2^{4}(L)
(3R_2(L))^d \sum_{k=1}^{n}
\biggl(\sum_{y\in B(x_k) \cap\mathcal{P}(0,L)}
P_{z,\omega}\bigl(T_{y} <\infty, J(L)\bigr)^2 \\
&&\hspace*{151.1pt}{} + P_{z,\omega}\bigl(T_{\partial_+ \mathcal{P}(0,L)} \not=
T_{\partial\mathcal{P}(0,L)}\bigr)^2 \biggr) \\
&\leq&
C
U^{-d} R_2^{4}(L)
(3R_2(L))^{2d} \biggl(\sum_{y\in\mathcal{P}(0,L)} P_{z,\omega}\bigl(T_y <
\infty, J(L)\bigr)^2 \\
&&\hspace*{106.9pt}\hspace*{7.9pt}{} + nP_{z,\omega} \bigl(T_{\partial\mathcal{P}(0,L)} \not=
T_{\partial_+ \mathcal{P}(0,L)}\bigr)^2 \biggr)\\
&\leq& R_3^2(L) U^{-d}
\end{eqnarray*}
for $L$ large enough and
where to obtain the second line we applied the Cauchy--Schwarz
inequality in combination with $\vert B(x_k) \vert\leq(3R_2(L))^d$.
Therefore, using (\ref{eq:differenceBound2}) and Lemma \ref
{lem:Azuma}, for each $v\in H_{M+U}$ we have
\begin{eqnarray*}
&&\PP\bigl( \bigl\vert P_z\bigl( X_{T_{M+U}}=v, A_L, J(L) \vert \mathcal{G}_n\bigr)
- P_z\bigl(X_{T_{M+U}}=v, A_L, J(L)\bigr) \bigr\vert> L^{1-d}/4 \bigr)\\
&&\qquad\leq2e^{-{U^{2\eta}}/({32R_4(L)})}
\end{eqnarray*}
with $\eta:=\frac{d-({d-1})/{\vartheta'}}{2}>0$ [here we use the
assumption $\vartheta' > (d-1)/d$ to
guarantee the positivity of $\eta$].
We define the subset
\[
T(L):=
\mathop{\mathop{\bigcap_{M \in\{{\lfloor(2/5)L^2 \rfloor},
\ldots, L^2\}
,}}_{v \in H_{M+U},}}_{z \in\tilde{\mathcal{P}}(0,L)}
\bigl\{ \bigl\vert
P_z (X_{T_{M+U}}=v \vert
\mathcal{G}_n)
- P_z(X_{T_{M+U}}=v)
\bigr\vert\leq L^{1-d}/2 \bigr\}
\]
of $\Omega$.
Now for any of these
choices of $M$, $v$ and $z$, we obtain
\begin{eqnarray*}
&&\PP\bigl( \bigl\vert
P_z (X_{T_{M+U}}=v \vert
\mathcal{G}_n)
- P_z(X_{T_{M+U}}=v)
\bigr\vert> L^{1-d}/2 \bigr)\\
&&\qquad\leq\PP\bigl(\bigl\vert P_z\bigl(X_{T_{M+U}}=v, A_L , J(L) \vert
\mathcal{G}_n\bigr)\\
&&\hspace*{47.6pt}{}-P_z\bigl(X_{T_{M+U}}=v, A_L, J(L)\bigr)
\bigr\vert> L^{1-d}/4 \bigr)\\
&&\qquad\quad{} + \PP\bigl( P_z\bigl(A_L^c \cup J(L)^c \vert \mathcal{G}_n\bigr) >
L^{1-d}/8 \bigr) \\
&&\qquad\quad{}+ \PP\bigl( P_z\bigl(A_L^c \cup J(L)^c\bigr) > L^{1-d}/8 \bigr).
\end{eqnarray*}
Thus, in combination with Lemma \ref{lem:ANEst}
and Proposition \ref{prop:closenessbase}\hyperlink
{item:quenchedBadExitEst}{(i)},
and since the previous bounds were uniform in the (at most polynomially
many) admissible choices of $M$,
$v$ and $z$,\setcounter{footnote}{3}\footnote{More precisely, in order to have only
polynomially many choices for $v$,
we restrict $v$ to be contained in the union of all admissible
hypercubes $Q^{(2)}$ appearing in (\ref{eq:Q2Def}).}
we get that
%
%
\begin{equation} \label{eq:TProb}
\PP(T(\cdot)^c) \in\mathcal{S}(\N).
\end{equation}
Now in order to estimate
\[
\vert
P_{z,\omega} (X_{T_M}\in Q)
-P_z (X_{T_M}\in Q)
\vert,
\]
we denote by $c(Q)$ the centre of the cube $Q$ and set
$
c'(Q):=c(Q)+\frac{U}{{\hat{v}} \cdot e_1} {\hat{v}}.
$
Furthermore, set
\[
Q^{(1)} := \bigl\{y\in H_{M+U} \dvtx \Vert y-c'(Q)\Vert_\infty<
(0.9)^{1/(d-1)}L^{\vartheta}/2 \bigr\}
\]
and
%
%
\begin{equation} \label{eq:Q2Def}
Q^{(2)} := \bigl\{ y\in H_{M+U} \dvtx \Vert y-c'(Q)\Vert_\infty<
(1.1)^{1/(d-1)}L^{\vartheta}/2 \bigr\}.
\end{equation}
Then by standard annealed estimates there exists $\varphi\in\mathcal
{S}(\N)$ such that for all
$z \in\tilde{\mathcal{P}}(0,L)$,
%
%
\begin{eqnarray} \label{eq:smallAnnBoxEst}
P_z\bigl(X_{T_{M+U}}\in Q^{(1)}\bigr)&<&P_z(X_{T_{M}}\in Q)+ \varphi(L),
\\
\label{eq:largeAnnBoxEst}
P_z\bigl(X_{T_{M+U}}\in Q^{(2)}\bigr)&>&P_z(X_{T_{M}}\in Q)-\varphi(L),
\\
\label{eq:smallQuenchedBoxEst}
P_z\bigl(X_{T_{M+U}}\in Q^{(1)} \vert \mathcal{G}_n\bigr)&<& P_{z,\omega}
(X_{T_{M}}\in Q)+\varphi(L)
\end{eqnarray}
and
%
%
\begin{equation}\label{eq:largeQuenchedBoxEst}
P_z\bigl(X_{T_{M+U}}\in Q^{(2)} \vert \mathcal{G}_n\bigr)> P_{z,\omega}
(X_{T_{M}}\in Q)-\varphi(L)
\end{equation}
for $\omega\in\mathcal{A}_L$.
Indeed, in order to prove equation (\ref{eq:smallAnnBoxEst}) note that
\begin{eqnarray*}
P_z\bigl(X_{T_{M+U}} \in Q^{(1)}\bigr)
&\leq&
P_z\bigl(X_{T_M} \in Q, X_{T_{M+U}} \in Q^{(1)}\bigr)\\
&&{}
+ P_z\bigl(X_{T_M} \notin Q, X_{T_{M+U}} \in Q^{(1)}\bigr).
\end{eqnarray*}
By Lemma \ref{lem:ANEst} and restricting to $A_L$, using Azuma's
inequality one can show that
\[
\sup_{z \in\tilde{\mathcal{P}}(0,L)} P_z\bigl(X_{T_M} \notin Q,
X_{T_{M+U}} \in Q^{(1)}\bigr)
\]
is contained in $\mathcal{S}(\N)$ as a function in $L$;
this then implies (\ref{eq:smallAnnBoxEst}).
The remaining inequalities are shown in similar ways.

In order to make use of (\ref{eq:smallAnnBoxEst}) to (\ref
{eq:largeQuenchedBoxEst}), note that for $\omega\in T(L) \cap\mathcal{A}_L$
we get with Lemma \ref{lem:annExitDistDerivatives}\hyperlink
{item:hittingProbBd}{(a)}
that
\begin{eqnarray*}
&&\bigl\vert P_{z} \bigl(X_{T_{M+U}} \in Q^{(1)} \vert \mathcal{G}_n\bigr) -
P_z\bigl(X_{T_{M+U}} \in Q^{(2)}\bigr) \bigr\vert\\
&&\qquad\leq\bigl\vert P_{z} \bigl(X_{T_{M+U}} \in Q^{(1)} \vert \mathcal{G}_n\bigr)
- P_z\bigl(X_{T_{M+U}} \in Q^{(1)}\bigr) \bigr\vert\\
&&\qquad\quad{}+ \biggl\vert\sum_{v \in Q^{(2)} \setminus Q^{(1)}} P_z(X_{T_{M+U}} =
v) \biggr\vert\\
&&\qquad\leq\bigl\vert P_{z} \bigl(X_{T_{M+U}} \in Q^{(1)} \vert \mathcal{G}_n\bigr)
- P_z\bigl(X_{T_{M+U}} \in Q^{(1)}\bigr) \bigr\vert
+ \bigl\vert Q^{(2)} \setminus Q^{(1)} \bigr\vert C L^{1-d}\\
&&\qquad \leq C L^{(\vartheta-1) (d-1)}
\end{eqnarray*}
for some constant $C$.
If $P_{z,\omega}(X_{T_M} \in Q) \leq P_z(X_{T_M} \in Q)$, then this
estimate in combination
with (\ref{eq:largeAnnBoxEst}) and (\ref{eq:smallQuenchedBoxEst})
yields (\ref{eq:quenchedAnnealedLargeBoxes}).
Otherwise, again for $\omega\in T(L) \cap\mathcal{A}_L$ we compute
\begin{eqnarray*}
&&\bigl\vert P_{z, \omega} \bigl(X_{T_{M+U}} \in Q^{(2)}\bigr) - P_z\bigl(X_{T_{M+U}} \in
Q^{(1)}\bigr) \bigr\vert\\
&&\qquad\leq\bigl\vert P_{z} \bigl(X_{T_{M+U}} \in Q^{(2)} \vert \mathcal{G}_n\bigr)
- P_z\bigl(X_{T_{M+U}} \in Q^{(2)}\bigr) \bigr\vert
+ \bigl\vert Q^{(2)} \setminus Q^{(1)} \bigr\vert C L^{1-d}\\
&&\qquad\leq C L^{(\vartheta-1) (d-1)},
\end{eqnarray*}
which in combination
with (\ref{eq:smallAnnBoxEst}) and (\ref{eq:largeQuenchedBoxEst})
again implies (\ref{eq:quenchedAnnealedLargeBoxes}).

Thus, for $C$ large enough,
and since the bounds we derived so far were uniform in the admissible choices
of $M$, $z$ and $Q$,
it follows that $T(L)\cap\mathcal{A}_L \subseteq B^{\vartheta}(L)$.
Therefore, employing (\ref{eq:TProb}), we get $\PP(B^{\vartheta
}(\cdot)^c)\in\mathcal{S}(\N)$.
\end{pf}

Departing from Lemma \ref{lem:quenchedAnnealedLargeBoxes}, due to the
following result,
for a large set of environments we can bound from above the quenched
probability of hitting a hyperplane in a hypercube of side length
$\lceil L^\vartheta\rceil$ for any $\vartheta\in(0,1]$.
%
\begin{lemma}\label{lem:alltheta}
For $\vartheta\in(0,1]$ and $h \in\N$, denote by $\overline
{B}{}^{\vartheta}_h(L)$ the set of
those~$\omega$ for which
for all $z\in{\tilde{\mathcal{P}}(0,L)}$, all $M \in\{{\lfloor
\frac12 L^2 \rfloor}, \ldots, L^2\}$
and all $(d-1)$-dimensional hypercubes $Q$ of side length $\lceil
L^{\vartheta} \rceil$ which are contained in $H_M$,
%
%
\begin{equation}\label{eq:alltheta}
P_{z,\omega} (X_{T_M}\in Q)
\leq R_h(L) L^{(\vartheta-1)(d-1)}.
\end{equation}
Then there exists
$h(\vartheta) \in\N$ such that $\PP(\overline{B}{}^{\vartheta
}_{h(\vartheta)} (L)^c)$ is contained
in $\mathcal{S}(\N)$ as a function in $L$.
\end{lemma}
\begin{pf}
We prove the lemma by descending induction on $\vartheta$.
Lemma \ref{lem:quenchedAnnealedLargeBoxes} in combination with Lemma
\ref{lem:annExitDistDerivatives}\hyperlink{item:hittingProbBd}{(a)}
implies that $\PP(\overline{B}{}^{\vartheta}_1 (\cdot)^c)
\in\mathcal{S}(\N)$ for each
$\vartheta\in(\frac{d-1}{d}, 1]$. For the induction step, assume
that the statement of the lemma holds for some $\vartheta'$
and choose $\vartheta$ such that $\rho:=\frac{\vartheta}{\vartheta
'} \in(\frac{d-1}{d}, 1]$.
Set $h':=h(\vartheta')$. For $z \in\Z^d$, define
the canonical shift on $\Z^d$ via $\sigma_z\dvtx\Z^d \ni x \mapsto x+z
\in\Z^d$
and let
\[
G := B^{\rho}(L) \cap G_{{\lfloor L^\rho\rfloor}}^* \cap
\bigcap_{z\in{\mathcal{P}}(0,L)}
\sigma_z (\overline{B}{}^{\vartheta'}_{h'} ({\lfloor L^\rho
\rfloor})),
\]
where $B^{\rho}(L)$ as in Lemma \ref{lem:quenchedAnnealedLargeBoxes} and
\[
G_{{\lfloor L^\rho\rfloor}}^* := \bigcap_{x\in{\mathcal{P}}(0,L)}
\bigcap
_{y \in\tilde{\mathcal{P}}(x,{\lfloor L^\rho\rfloor})}
\bigl\{
P_{y,\omega} \bigl( T_{\partial{\mathcal{P}}(x, {\lfloor L^\rho
\rfloor})}
\not= T_{\partial_+{\mathcal{P}}(x, {\lfloor L^\rho\rfloor} )}
\bigr)<e^{-R_2({\lfloor L^\rho\rfloor})^\gamma}
\bigr\}.
\]
The translation invariance of $\PP$ implies that
\[
\PP( \sigma_z ( \overline{B}{}^{\vartheta'}_{h'} ({\lfloor
L^\rho\rfloor}) ) )
=\PP( \overline{B}{}^{\vartheta'}_{h'}( {\lfloor L^\rho\rfloor
}) ),
\]
and therefore,
as a consequence of Proposition \ref{prop:closenessbase}\hyperlink
{item:quenchedBadExitEst}{(i)},
$\PP(G^c)$ is contained in~$\mathcal{S}(\N)$ as a function
in $L$. Thus, it is sufficient to show that for some~$h$ and all $L$
large enough, we have that
$G \subseteq\overline{B}{}^{\vartheta}_{h}(L)$.
To this end we fix \mbox{$\omega\in G$}, \mbox{$z \in\tilde{\mathcal{P}}(0,L)$},
$M \in\{{\lfloor\frac12L^2 \rfloor}, \ldots, L^2\}$
and
a cube $Q$ of side length $\lceil L^{\vartheta} \rceil$ in ${\mathcal
{P}}(0,\allowbreak L)\cap H_M$. Let~$c(Q)$ be the centre of $Q$
and $x'$ be an element of $\Z^d$ closest to~$c(Q)- \frac{{\lfloor
L^\rho\rfloor}^2}{\hat{v} \cdot e_1} \hat{v}$.
Due to the strong Markov property and the fact that $\omega\in
G_{{\lfloor L^\rho\rfloor}}^*$,
%
%
\begin{eqnarray}\label{eq:markov}
&&\biggl\vert P_{z,\omega} (X_{T_M}\in Q)
-\sum_{v \in H_{M-{\lfloor L^\rho\rfloor}^2} \cap{\mathcal
{P}}(x', {\lfloor L^\rho\rfloor})}
P_{z,\omega} (X_{T_{M-{\lfloor
L^\rho\rfloor}^2}}=v)\nonumber\\[-8pt]\\[-8pt]
&&\hspace*{187.4pt}{}\times
P_{v,\omega} (X_{T_M}\in Q) \biggr\vert\nonumber
\end{eqnarray}
is contained in $\mathcal{S}(\N)$ as a function in $L$.
To estimate the second factor of the sum, observe that since
\[
\omega\in\bigcap_{z\in{\mathcal{P}}(0,L)}
\sigma_z(\overline{B}{}^{\vartheta'}_{h'}( {\lfloor L^\rho
\rfloor} )),
\]
we get that for every $v\in H_{M-{\lfloor L^\rho\rfloor}^2}$,
%
%
\begin{equation}\label{eq:katan}
P_{v,\omega} (X_{T_M}\in Q)
<R_{h'}(L) L^{\rho(\vartheta'-1)(d-1)}=R_{h'}(L) L^{(\vartheta-\rho)(d-1)}.
\end{equation}
With respect to the first factor of the sum, for $L$ large enough,
$H_{M-{\lfloor L^\rho\rfloor}^2} \cap{\mathcal{P}}(x', {\lfloor
L^\rho\rfloor})$ is
the union of less than
$R_7(L)$ cubes of side length ${\lfloor L^\rho\rfloor}$.
Since $\omega\in B^{\rho}(L)$, we deduce that for every cube $Q'$ of
side length ${\lfloor L^\rho\rfloor}$
that is contained in $H_{M-{\lfloor L^\rho\rfloor}^2} \cap{\mathcal
{P}}(0,L)$, one has
%
%
\begin{equation} \label{eq:gadol}
P_{z,\omega} (X_{T_{M-{\lfloor L^\rho\rfloor}^2}}\in Q')
<R_2(L)L^{(\rho-1)(d-1)}
\end{equation}
for $L$ large enough.
Combining (\ref{eq:markov}), (\ref{eq:katan}) and (\ref{eq:gadol}),
we infer that
\begin{eqnarray*}
P_{z,\omega} (X_{T_M}\in Q)
&\leq& R_7(L)R_{h'}(L) L^{(\vartheta-\rho)(d-1)}\cdot
R_2(L)L^{(\rho-1)(d-1)}\\
&\leq& R_h(L)L^{(\vartheta-1)(d-1)}
\end{eqnarray*}
for $h=\max\{7,h'\}+1$ and $L$ large enough.

Noting that the above estimates are uniform in the (at most
polynomially many) admissible choices
of $z$, $M$ and $Q$, this finishes the proof.
\end{pf}

The next result employs the previous lemmas to yield
bounds on the difference of certain annealed and semi-annealed hitting
probabilities.
%
\begin{lemma} \label{lem:distu}
Let $\mathcal{G}$ be the $\sigma$-algebra\vspace*{1pt} in $\Omega$ generated by
the functions
$\{\Omega\ni\omega\mapsto\omega(z) \dvtx z \cdot e_1 \leq L^2\}$.
Let $\eta\in(0, \frac{6}{d-1} \wedge 1)$, $U:={\lfloor L^\eta\rfloor}$ and denote
by~$B(L,\eta)$\vadjust{\goodbreak}
the set of those $\omega$ for which for all
$z\in{\tilde{\mathcal{P}}}(0,L)$ and all $v\in H_{L^2+U}$, one has
\[
\vert
P_{z} (X_{T_{L^2+U}}=v \vert \mathcal{G})
- P_z(X_{T_{L^2+U}}=v)
\vert
\leq L^{1-d}U^{({1-d})/{3}}.
\]
Then $\PP(B(L,\eta)^c)$ is contained in $\mathcal{S}(\N)$
as a function in $L$.
\end{lemma}
\begin{pf}
Let $v\in H_{L^2+U}$ and let $\vartheta>0$ be such that $\vartheta<
\frac{1}{12}\eta$. Define $K_L$ to be the natural
number such that $2^{-K_L}{L^2}>U\geq
2^{-K_L-1}{L^2}$, and for $k \in\{1, \ldots, K_L-1 \}$ we set
\[
{\mathcal{P}}^{(k)}:={\mathcal{P}}(0,L)\cap\{x \dvtx L^2 -
2^{-k}{L^2} < x \cdot e_1 \leq L^2 - 2^{-k-1}{L^2} \}.
\]
In addition, we take
\begin{eqnarray*}
{\mathcal{P}}^{(K_L)}&:=&{\mathcal{P}}(0,L)\cap\{x \dvtx L^2 - 2^{-K_L}L^2
< x \cdot e_1 \leq{L^2}\},
\\
{\mathcal{P}}^{(0)}&:=&{\mathcal{P}}(0,L)\cap\{x \dvtx x \cdot e_1 \leq
L^2/2\}
\end{eqnarray*}
and
\[
F(v):= \{ x\in{\mathcal{P}}(0,L) \dvtx\Vert x-u(v,x) \Vert_1 \leq
R_7(L) \Vert(v-x) \cdot e_1 \Vert_1^{1/2} \},
\]
where $u(v,x):=v+ \frac{(x-v)\cdot e_1}{\hat{v} \cdot e_1} \hat{v}$.
Then for $k \in\{0, \ldots, K_L\}$ we define
\[
{\mathcal{P}}^{(k)}(v) :={\mathcal{P}}^{(k)}\cap F(v)
\]
and
\[
\hat{{\mathcal{P}}}^{(k)}(v):= \bigl\{y \in\Z^d\dvtx\exists{x\in
{\mathcal{P}}^{(k)}(v)} \mbox{ such that } \Vert x-y \Vert_1
<R_2(L) \bigr\}.
\]
See Figure \ref{fig6} for an illustration.

%
%
\begin{figure}

\includegraphics{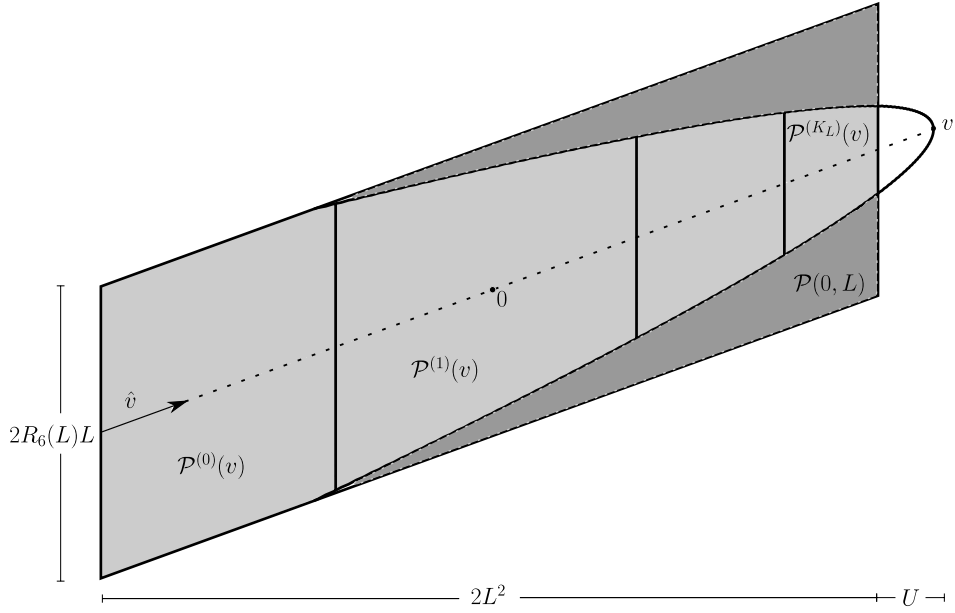}

\caption{The sets ${\mathcal{P}}^{(k)}(v)$ contained in $\mathcal{P}(0,L)$.}
\label{fig6}
\end{figure}

Similarly to the previous, we use a
lexicographic enumeration $x_1, x_2, \ldots, x_n$ of
%
%
\begin{equation} \label{eq:PHatDef}
\hat{\mathcal{P}} := \bigcup_{k=0}^{K_L} \hat{\mathcal{P}}^{(k)}
\end{equation}
and the corresponding filtration $\{\mathcal{G}_i\}_{i \in\{0, \ldots ,
n\}}$. We consider the martingale $ M_i:=P_z(X_{T_{L^2+U}}=v, A_L, J(L)
\vert \mathcal{G}_i)$. Again, in order to use
Lemma \ref{lem:Azuma}, we need to bound $ U_i :=\esssup( \vert
M_{i}-M_{i-1} \vert \vert\mathcal{G}_{i-1}). $ With the same reasoning
as in the proof of Lemma \ref {lem:quenchedAnnealedLargeBoxes} and with
Lemma \ref{lem:annExitDistDerivatives}\hyperlink
{item:firstStartPtDer}{(c)}, we obtain for $i$ such that $x_i
\in{\mathcal{P}}^{(k)}(v)$:
\[
U_i
\leq CR_2(L) P_z\bigl(T_{x_i} < \infty, A_L, J(L) \vert \mathcal{G}_{i-1}\bigr)
\cdot L^{-d}2^{(k+1)({d}/{2})}.
\]
To obtain a useful upper bound
for $U_i$ with $k \in\{0, \ldots, K_L\}$ and $\omega\in\overline
{B}{}^{\vartheta}_{h(\vartheta)} \cap J(L)$, we will estimate
\[
V_\omega(k)=\sum_{x\in{\mathcal{P}}^{(k)}(v)} P_{z,\omega} \bigl(T_x <
\infty,A_L, J(L)\bigr)^2.
\]
Using (\ref{eq:JHittingProbBd}), we get for $\omega\in J(L)$ that
\[
V_\omega(0) \leq R_3(L).
\]
Now choose $h(\vartheta) \geq8$ such that the implication of Lemma
\ref{lem:alltheta}
holds true; then for
$k>0$ as well as $B(x)$ as in the proof of Lemma \ref
{lem:quenchedAnnealedLargeBoxes},
%
%
\begin{eqnarray}\label{eq:bcbt}
V_\omega(k)&=&\sum_{x\in{\mathcal{P}}^{(k)}(v)} P_{z,\omega} (T_x <
\infty,A_L)^2 \nonumber\\
&\leq&\sum_{x\in{\mathcal{P}}^{(k)}(v)} \biggl(
\sum_{y \in B(x)}P_{z,\omega} (X_{T_{y \cdot e_1}}=y)
\biggr)^2\nonumber\\
&\leq&(3R_2(L))^{d-1} \sum_{y\in\hat{{\mathcal{P}}}^{(k)}(v)}
P_{z,\omega} (X_{T_{y \cdot e_1}}=y)^2 \\
&\leq&(3R_2(L))^{d-1} \sum_{y\in\hat{{\mathcal{P}}}^{(k)}(v)}
R^2_{h(\vartheta)} (L)L^{2(\vartheta-1)(d-1)} \nonumber\\
&\leq& R_{h(\vartheta)+1}(L)L^{2 ( ({d+1})/{2}+(\vartheta
-1)(d-1) )}2^{-k {\lfloor({d+1})/{2} \rfloor} } \nonumber
\end{eqnarray}
for $L$ large enough,
where inequality (\ref{eq:bcbt}) follows from the fact that $\omega
\in\overline{B}{}^{\vartheta}_{h(\vartheta)}(L)$.

Therefore, we get that for $\omega\in J(L) \cap\overline
{B}{}^\vartheta_{h(\vartheta)}$ we have
\begin{eqnarray*}
\esssup\Biggl( \sum_{i=1}^{n} U_i^2 \Biggr)
&\leq& CR_2^2(L) \sum_{k=0}^{K_L} V_\omega(k)L^{-2d}2^{kd}\\
&\leq& CR_{h(\vartheta)+1}(L)L^{-2d}\\
&&{} + CR_{h+1}(L)L^{2 ( ({d+1})/{2}+(\vartheta-1)(d-1) )
-2d}\sum_{k=1}^{K_L} 2^{kd-k({d+1})/{2}}\\
&\leq& CR_{h(\vartheta)+1}(L)\bigl(L^{-2d}+L^{ 3-3d + 2(d-1)\vartheta
}2^{K_L ({d-1})/{2}}\bigr)\\
&\leq& CR_{h(\vartheta)+1}(L)\bigl(L^{-2d}+L^{ 2-2d + 2(d-1)\vartheta
}U^{-({d-1})/{2}}\bigr)\\
&\leq& L^{ 2-2d}U^{({1-d})/{3}}/2
\end{eqnarray*}
for $L$ large enough and where the penultimate inequality follows from
the definition
of $K_L$, while the last inequality is due to the choice of $\vartheta$.
Thus, Lemma \ref{lem:Azuma} yields that
%
%
\begin{equation} \label{eq:quenchedAveragedDist}
\PP\bigl( \bigl\vert P_{z}(X_{T_{L^2+U}}=v, A_L \vert \mathcal{G}_n)
-P_z(X_{T_{L^2+U}}=v, A_L) \bigr\vert
> L^{1-d}U^{({1-d})/{3}} /2
\bigr)\hspace*{-32pt}
\end{equation}
is contained in $\mathcal{S}(\N)$ as a function in $L$,
uniformly in the admissible choices of $z$ and $v$.

Observe furthermore that with $\hat{\mathcal{P}}$ defined in (\ref
{eq:PHatDef}),
for example, by Azuma's inequality,
%
\[
P_{z} (X_{T_{L^2+U}}=v, A_L, T_{\partial\hat{\mathcal{P}}} <
T_v )
\]
is contained in $\mathcal{S}(\N)$,
and thus, due to Markov's inequality, so is
\[
\PP\bigl( P_{z} (X_{T_{L^2 + U}} = v, A_L, T_{\partial\hat
{\mathcal{P}}} < T_v \vert \mathcal{G} ) \geq
L^{1-d}U^{({1-d})/{3}} /2 \bigr).
\]
In combination with (\ref{eq:quenchedAveragedDist})
and the fact that
\begin{eqnarray*}
&&\bigl\{ \omega\dvtx\bigl\vert P_{z} (X_{T_{L^2+U}}=v, A_L \vert
\mathcal{G} )
- P_{z} (X_{T_{L^2+U}}=v, A_L \vert \mathcal{G}_n)
\bigr\vert\geq L^{1-d}U^{({1-d})/{3}} /2 \bigr\} \\
&&\qquad{}\subset
\bigl\{ \omega\dvtx P_{z} (X_{T_{L^2+U}}=v, A_L, T_{\partial\hat
{\mathcal{P}}} < T_v \vert \mathcal{G} ) \geq
L^{1-d}U^{({1-d})/{3}} /2 \bigr\},
\end{eqnarray*}
this supplies us with the fact that
\[
\PP\bigl(
\bigl\vert
P_{z} ( X_{T_{L^2+U}}=v, A_L \vert \mathcal{G})
-P_z(X_{T_{L^2+U}}=v, A_L )
\bigr\vert
> L^{1-d}U^{({1-d})/{3}}
\bigr)
\]
is contained in $\mathcal{S}(\N)$ also.

A union bound in combination with Lemma \ref{lem:ANEst} and
Lemma \ref{lem:annealedPosExitProbEst} completes the proof of the lemma.
\end{pf}
%
\begin{lemma}\label{lem:goodbound}
For any $\vartheta\in(0,\frac{6}{d-1} \wedge 1)$ denote by $D^{\vartheta}(L)$ the set of
those~$\omega$ for which
for all $z\in\tilde{{\mathcal{P}}}(0,L)$ and all
$(d-1)$-dimensional\vadjust{\goodbreak}
hypercubes~$Q$ of side length $\lceil L^{\vartheta} \rceil$ that are
contained in $\partial_+{\mathcal{P}}(0,L)$,
%
%
\begin{eqnarray}\label{eq:goodbound}
&&\bigl\vert
P_{z,\omega} \bigl(X_{T_{\partial{\mathcal{P}}(0,L)}}\in Q \vert
T_{\partial\mathcal{P}(0,L)} = T_{\partial_+ \mathcal{P}(0,L)}
\bigr)\nonumber\\[-2pt]
&&\quad\hspace*{0pt}{} - P_z \bigl( X_{T_{\partial{\mathcal{P}}(0,L)}}\in Q \vert
T_{\partial\mathcal{P}(0,L)} = T_{\partial_+ \mathcal{P}(0,L)}
\bigr)
\bigr\vert\\[-2pt]
&&\qquad\leq
L^{(\vartheta-1)(d-1)-\vartheta({d-1})/({d+1})
}.\nonumber\vspace*{-2pt}
\end{eqnarray}
Then $\PP(D^{\vartheta}(L)^c)$ is contained in $\mathcal
{S}(\N)$ as a function in $L$.\vspace*{-3pt}
\end{lemma}
\begin{pf}
Choose $\vartheta' \in(\frac34\vartheta, \vartheta)$ and
$U:={\lfloor L^{{4\vartheta'}/({d+1})} \rfloor}$.
Then by Lemma~\ref{lem:distu} and Proposition \ref
{prop:closenessbase}\hyperlink{item:quenchedBadExitEst}{(i)},
we know that
$\PP( B(\cdot,{\frac{4\vartheta'}{d+1}})^c \cup{G^{\mbox{\hyperlink
{item:quenchedBadExitEst}{\iitem}}}_\cdot}^c ) \in\mathcal{S}(\N)$ whence
it is sufficient to show that $B(L,{\frac{4\vartheta'}{d+1}}) \cap
G^{\mbox{\hyperlink{item:quenchedBadExitEst}{\iitem}}}_L \subset D^{\vartheta}(L)$;
this we will do similarly to the last step of the proof of Lemma \ref
{lem:quenchedAnnealedLargeBoxes}.
We denote by $c(Q)$
one of those elements of $\Z^d$ closest to the centre of $Q$ and let
$x' \in H_{L^2 + U}$ be one of the lattice points closest to
$c(Q) +\frac{U}{{\hat{v}} \cdot e_1} {\hat{v}}$.
Furthermore, let $Q^{(1)}$ and~$Q^{(2)}$ be $(d-1)$-dimensional hypercubes
that are contained in $H_{L^2+U}$ and are centred in $x'$, such that
the side length of $Q^{(1)}$ is
${\lfloor L^{\vartheta}-R_6(L)\sqrt{U} \rfloor}$ and the side length of
$Q^{(2)}$ is $\lceil L^{\vartheta}+R_6(L)\sqrt{U} \rceil$.
Then due to Lemma \ref{lem:distu}, on $B(L,{\frac{4\vartheta
'}{d+1}})$ for $i \in\{1,2\}$,
%
%
\begin{equation} \label{eq:hitOne}
\bigl\vert
P_{z} \bigl( X_{T_{L^2+U}}\in Q^{(i)} \vert \mathcal{G} \bigr)
-P_z \bigl( X_{T_{L^2+U}}\in Q^{(i)} \bigr)
\bigr\vert
\leq\bigl\vert Q^{(i)} \bigr\vert L^{1-d}U^{({1-d})/{3}}\hspace*{-28pt}\vspace*{-2pt}
\end{equation}
for all corresponding $z$ and $Q$.
Now similarly to (\ref{eq:smallAnnBoxEst}) to (\ref
{eq:largeQuenchedBoxEst}), there exists $\varphi\in\mathcal{S}(\N)$ such
that for all such $z$ and $Q$,
%
%
\begin{eqnarray}
&&
P_z \bigl( X_{T_{L^2 + U}}\in Q^{(1)} \bigr) - \varphi(L) \nonumber\\[-10pt]\\[-10pt]
&&\qquad< P_z (X_{T_{L^2}}\in Q ) < P_z \bigl( X_{T_{L^2 + U}}\in Q^{(2)} \bigr)
+ \varphi(L)\nonumber\vspace*{-2pt}
\end{eqnarray}
as well as
%
%
\begin{eqnarray}\label{eq:hitTwo}
&& P_z\bigl( X_{T_{L^2 + U}}\in Q^{(1)} \vert \mathcal{G} \bigr) -
\varphi(L)\nonumber\\[-10pt]\\[-10pt]
&&\qquad < P_{z,\omega} (X_{T_{L^2}}\in Q )
< P_z \bigl( X_{T_{L^2 + U}}\in Q^{(2)} \vert \mathcal{G}
\bigr)+\varphi(L).\nonumber\vspace*{-2pt}
\end{eqnarray}
Proposition \ref{prop:closenessbase}\hyperlink
{item:quenchedBadExitEst}{(i)} and Lemma
\ref{lem:annealedPosExitProbEst} imply that for $\omega\in
G^{\mbox{\hyperlink{item:quenchedBadExitEst}{\iitem}}}_L$,
\[
\bigl\vert P_{z,\omega} \bigl( X_{T_{\partial\mathcal{P}(0,L)}} \in
Q \vert T_{\partial\mathcal{P}(0,L)} = T_{\partial_+ \mathcal
{P}(0,L)} \bigr)
- P_{z,\omega} ( X_{T_{L^2}} \in Q) \bigr\vert\vspace*{-2pt}
\]
and
\[
\bigl\vert P_z \bigl( X_{T_{\partial\mathcal{P}(0,L)}} \in Q \vert
T_{\partial\mathcal{P}(0,L)} = T_{\partial_+ \mathcal{P}(0,L)}
\bigr)
- P_{z}(X_{T_{L^2}} \in Q ) \bigr\vert\vspace*{-2pt}
\]
are both contained in $\mathcal{S}(\N)$ as functions in $L$.
Therefore, for $\omega\in B(L,\frac{4\vartheta'}{d+1}) \cap
G^{\mbox{\hyperlink{item:quenchedBadExitEst}{\iitem}}}_L$,
using (\ref{eq:hitOne}) to (\ref{eq:hitTwo})
and as a consequence of Lemma \ref
{lem:annExitDistDerivatives}\hyperlink{item:hittingProbBd}{(a)},
\begin{eqnarray*}
&&\bigl\vert
P_{z,\omega} \bigl(X_{T_{\partial{\mathcal{P}}(0,L)}}\in Q
\vert T_{\partial\mathcal{P}(0,L)} = T_{\partial_+ \mathcal
{P}(0,L)} \bigr) \\[-2pt]
&&\quad{} -P_z \bigl( X_{T_{\partial{\mathcal{P}}(0,L)}} \in Q \vert
T_{\partial\mathcal{P}(0,L)} = T_{\partial_+ \mathcal{P}(0,L)}
\bigr)
\bigr\vert\\[-2pt]
&&\qquad\leq\vert P_{z,\omega} (X_{T_{L^2}} \in Q )
-P_z (X_{T_{L^2}} \in Q ) \vert+ \varphi(L) \\[-2pt]
&&\qquad\leq
\bigl\vert Q^{(2)} \bigr\vert L^{1-d}U^{({1-d})/{3}}
+C\bigl( \bigl\vert Q^{(2)} \bigr\vert- \bigl\vert Q^{(1)} \bigr\vert\bigr)L^{1-d}+
\varphi(L)\\[-2pt]
&&\qquad\leq
C\bigl(
L^{(\vartheta-1)(d-1)}U^{({1-d})/{3}} + R_6(L)\sqrt{U}\lceil
L^\vartheta\rceil^{d-2} L^{1-d}
\bigr)\\[-2pt]
&&\qquad\leq L^{(\vartheta-1)(d-1) - \vartheta({d-1})/({d+1})}
\end{eqnarray*}
for $L$ large enough and some $\varphi\in\mathcal{S}(\N)$.
Here, we used that
$U={\lfloor L^{{4\vartheta'}/({d+1})} \rfloor}$ and $\vartheta'
\in(\frac
34 \vartheta, \vartheta)$ to obtain the last line.
\end{pf}

\subsection{\texorpdfstring{Proof of Proposition \protect\ref{prop:closenessbase}
    \protect\hyperlink{item:quenchedAnnealedProbDiff}{(iii)}}
{Proof of Proposition 3.4(iii)}}

Denote by $D^\vartheta(L)$ the set of all $\omega$ such that
\begin{eqnarray*}
&&\max_{z \in\tilde{{\mathcal{P}}}(0,L)} \max_Q
\bigl\vert P_{z,\omega} \bigl( X_{T_{\partial{\mathcal{P}}(0,L)}}
\in Q \vert T_{\partial\mathcal{P}(0,L)} = T_{\partial_+
\mathcal{P}(0,L)} \bigr)\\
&&\quad\hspace*{45.6pt}{} - P_z \bigl( X_{T_{\partial{\mathcal{P}}(0,L)}} \in Q \vert
T_{\partial\mathcal{P}(0,L)} = T_{\partial_+ \mathcal{P}(0,L)}
\bigr) \bigr\vert\\
&&\qquad\hspace*{45.6pt}< L^{(\vartheta- 1)(d-1) - \vartheta({d-1})/({d+1})}
\end{eqnarray*}
holds,
where the maximum in $Q$ is taken over all
$(d-1)$-dimensional hypercubes $Q \subset\partial_+ {\mathcal{P}}(0,L)$
of side length $\lceil{L^\vartheta} \rceil$.
Then for $\vartheta \in (0, \frac{6}{d-1} \wedge 1)$,
Lem\-ma~\ref{lem:goodbound} is applicable and
yields that
$\PP(D^\vartheta(L)^c)$ is contained in $\mathcal{S}(\N)$ as a~function in $L$.
In combination with Remark \ref{rem:PropProofStrat}, this finishes the proof.

\subsection{Further auxiliary results} \label{sec:furtherAuxRes}
The principal purpose of this subsection is to prove Lemma~\ref
{lem:sumOfClose} that
has been employed in step \hyperlink{item:goodExitOfBigBox}{(B)} in
the construction of
the auxiliary random walk on page \pageref{aaaa}.

We start with proving some further auxiliary results, parts of which
have been stated and employed above already.
\begin{pf*}{Proof of Lemma \ref{lem:bddFluctExp}} \label{pr:bddFluctExp}
We observe that due to Lemma
\ref{lem:annealedPosExitProbEst}, it is sufficient to establish (\ref
{eq:bddFluctExp}).
With
$\hat{v}_{L}$
as defined in (\ref{eq:vLDef}),
we obtain
%
%
\begin{eqnarray} \label{eq:bddFluctExpDecomp}
&&\biggl\Vert E_x X_{T_{\partial{\mathcal{P}}(0,L)}} -x - \frac{L^2-
x\cdot e_1}{\hat{v} \cdot e_1} \hat{v}
\biggr\Vert_1\nonumber\\
&&\qquad\leq\biggl\Vert E_x X_{T_{\partial{\mathcal{P}}(0,L)}} -x
- \frac{L^2 -x \cdot e_1}{\hat{v}_L \cdot e_1} \hat{v}_L \biggr\Vert
_1\\
&&\qquad\quad{} + \biggl\Vert\frac{L^2 - x \cdot e_1}{\hat{v}_L \cdot e_1} \hat
{v}_L - \frac{L^2 -x \cdot e_1}{\hat{v}\cdot e_1} \hat{v} \biggr\Vert_1.
\nonumber
\end{eqnarray}
To estimate the first summand on the right-hand side of (\ref
{eq:bddFluctExpDecomp}),
note that
for $H:= \inf\{n \in\N\dvtx\sum_{j=1}^n (X_{\tau_j}-X_{\tau_{j-1}})
\cdot e_1 \geq L^2\}$,
we can infer from Lemma~\ref{lem:annealedPosExitProbEst} and Lemma~\ref{lem:ANEst} that
%
%
\begin{equation} \label{eq:stoppingTimeRep}
\Vert E_x X_{T_{\partial\mathcal{P}(0,L)}}
- E_x (X_{\tau_H}, A_L) \Vert_1 \leq2R_2(L)
\end{equation}
for $L$ large enough.
Now
$
(\sum_{j=1}^n X_{\tau_j} - X_{\tau_{j-1}} - E_x(X_{\tau
_{j}}-X_{\tau_{j-1}} \vert A_L)
)_{n \in\{1, \ldots, 2L^2\}}
$
is a zero-mean martingale with respect to $P_x ( \cdot\vert A_L)$,
whence the optimal stopping theorem implies
%
%
\begin{eqnarray} \label{eq:stoppingThm}
E_x ( X_{\tau_H} - x \vert A_L ) &=& \bigl(E_x (H \vert A_L) -1\bigr)
\cdot E_x(X_{\tau_2}-X_{\tau_1} \vert A_L)\nonumber\\[-8pt]\\[-8pt]
&&{} + E_x(X_{\tau_1} -x
\vert A_L).\nonumber
\end{eqnarray}
But as a consequence of the conditioning on $A_L$, we have
%
%
\begin{equation} \label{eq:e1HittingDist}
\bigl\Vert E_x ( X_{\tau_H} - x \vert A_L ) \cdot e_1 - (L^2 - x
\cdot e_1) \bigr\Vert_1 \leq R_2(L).
\end{equation}
Since furthermore
%
%
\begin{equation} \label{eq:renewalExpDiff}
\bigl\Vert E_x (X_{\tau_2} - X_{\tau_1} \vert A_L) - E_x (X_{\tau
_1} -x \vert A_L) \bigr\Vert_1 \leq2R_2(L)
\end{equation}
using (\ref{eq:stoppingThm}) to (\ref{eq:renewalExpDiff}), we get
\[
\biggl\Vert E_x ( X_{\tau_H} - x \vert A_L )
- \frac{L^2 - x \cdot e_1}{\hat{v}_L \cdot e_1} \hat{v}_L \biggr\Vert
_1 \leq3R_2(L).
\]
Combining this with (\ref{eq:stoppingTimeRep}) we obtain
that the first summand on the right-hand side of (\ref
{eq:bddFluctExpDecomp}) is bounded from above by
$5R_2(L)$.

Furthermore,
the second summand on the right-hand side of (\ref{eq:bddFluctExpDecomp})
is contained in $\mathcal{S}(\N)$ as a function in $L$
due to Lemma \ref{lem:directionsDistance}.
This finishes the proof.
\end{pf*}

In the following, we will sometimes consider distributions $\mu
_{0,0}^{\sqrt{j}L}$ for \mbox{$j \in\N$},
and in particular, $\sqrt{j}L$ is not necessarily a natural number
anymore. However, as one may
check, this does not lead to any complications.
%
\begin{claim}\label{claim:annealedConvErr}
For $j \in\{1, \ldots, {\lfloor L^\chi\rfloor}^2\}$,
let $U$ be distributed according to the convolution $\mu_{0,0}^L * \mu
_{0,0}^{\sqrt{j-1}L}$.
Then $U$ can be represented as $U=\hat{U}+U'$ such that $\hat{U} \sim
\mu_{0,0}^{\sqrt{j}L}$
and
\[
P \bigl(\Vert U' \Vert_1 > 2R_2(L)\bigr) \leq C e^{-C^{-1}R_2(L)^{\gamma}}
\]
for some constant $C$ independent
of $j$ and $L$.
\end{claim}
\begin{pf}
Since we assume all appearing probability spaces to be large enough, it
is sufficient to construct $U, \hat{U}$ and $U'$
as desired.
First, observe that
for
\[
A_{k,N} := \bigl\{X^{*(n)} < R_2(N)\ \forall n \in\{1, \ldots, k\}
\bigr\},
\]
the same reasoning as in the proof of Lemma \ref{lem:ANEst} yields that
\[
P_0 \bigl( A_{({\lfloor L^\chi\rfloor}L)^2,L}^c \bigr) \leq C \exp
\{ -C^{-1}
R_2(L)^\gamma\}.
\]
This in combination with Azuma's inequality, Lemma \ref
{lem:annealedPosExitProbEst} and Lem\-ma~\ref{lem:bddFluctExp},
yields that for $L$ large enough we have
\begin{eqnarray*}
&&P_0 \bigl(T_{\partial\mathcal{P}(0,\sqrt{j-1}L)} \ne T_{\partial
_+ \mathcal{P}(0,\sqrt{j-1}L)} \bigr) \\
&&\qquad\leq P_0 \bigl( A_{({\lfloor L^\chi\rfloor}L)^2,L}^c \bigr) + 2d
({\lfloor L^\chi\rfloor}L)^2 \exp\biggl\{ - \frac
{R_6(L)^2}{4R_2(L)^2} \biggr\}\\
&&\qquad\leq C \exp\{-C^{-1} R_2(L)^\gamma\}.
\end{eqnarray*}
Now for $l \in\N$, let $n(l)$ be the unique natural number such that
$\tau_{n(l)-1} < T_l \leq\tau_{n(l)}$.
Then due to the above, in combination with Lemma \ref{lem:ANEst},
\[
P_0 \bigl( \Vert X_{\tau_{n((j-1)L^2)}} - X_{T_{\partial\mathcal
{P}(0,\sqrt{j-1}L)}} \Vert_1 \geq R_2(L) \bigr)
\leq C \exp\{-C^{-1} R_2(L)^\gamma\}.
\]
Now let $Z$ denote a RWRE coupled to $X$ in such a way that
$Z_0!=\!X_{T_{\partial_+ \mathcal{P}(0, \sqrt{j-1}L)}}$
and
\[
Z_{\tau^Z_1 + \cdot} = X_{\tau_{n((j-1)L^2)} + \cdot} - X_{\tau
_{n((j-1)L^2)}},
\]
whereas between times $0$ and $\tau^Z_1$ it evolves independently of $X$.
Then
%
%
\begin{eqnarray} \label{eq:manyProbEst}\quad
&&P_0\Bigl( \bigl\{ T_{\partial\mathcal{P}(0,\sqrt{j}L)} \not=
T_{\partial_+ \mathcal{P}(0,\sqrt{j}L)} \bigr\}\nonumber\\
&&\qquad\hspace*{0pt}{}
\cup\bigl\{T_{\partial\mathcal{P}(0,\sqrt{j-1}L)} \not=
T_{\partial_+ \mathcal{P}(0,\sqrt{j-1}L)} \bigr\}\cup\bigl\{T^Z_{\partial\mathcal{P}(Z_0,L)}
\not= T^Z_{\partial_+ \mathcal{P}(Z_0,L)} \bigr\}\nonumber\\[-8pt]\\[-8pt]
&&\qquad\hspace*{66.8pt}{}
\cup\Bigl\{\max_{0 \leq n \leq\tau_1^Z} \Vert Z_n - Z_0 \Vert_1
\geq R_2(L) \Bigr\}
\cup A_{({\lfloor L^\chi\rfloor}L)^2,L}^c\Bigr) \nonumber\\
&&\qquad\quad\leq C \exp\{-C^{-1} R_2(L)^\gamma\}\nonumber
\end{eqnarray}
for $C$ large enough and all $L$; restricted to the complement of the
event on the left-hand side
of (\ref{eq:manyProbEst}),
%
%
\begin{equation} \label{eq:XZbd}\qquad\qquad
\Vert X_{T_{\partial_+ \mathcal{P}(0,\sqrt{j-1}L)} } +
( Z_{T^Z_{\partial_+ \mathcal{P}(Z_0,L)}}
- Z_0 )
- X_{T_{\partial_+ \mathcal{P}(0,\sqrt{j}L)} } \Vert_1 \leq2R_2(L).
\end{equation}
Furthermore, with respect to
\[
P_0 \bigl( \cdot \vert T_{\partial\mathcal{P}(0,\sqrt
{j-1}L)} = T_{\partial_+ \mathcal{P}(0,\sqrt{j-1}L)},
T^Z_{\partial\mathcal{P}(Z_0,L)}
= T^Z_{\partial_+ \mathcal{P}(Z_0,L)} \bigr),
\]
the variable
\[
U:=X_{T_{\partial_+ \mathcal{P}(0,\sqrt{j-1}L)} } +
Z_{T^Z_{\partial_+ \mathcal{P}(Z_0,L)}}
- Z_0
\]
is distributed according to $\mu_{0,0}^{\sqrt{j-1}L} * \mu
_{0,0}^{L}$, while with respect to
\[
P_0 \bigl( \cdot \vert T_{\partial\mathcal{P}(0,\sqrt
{j}L)} = T_{\partial_+ \mathcal{P}(0,\sqrt{j}L)} \bigr),
\]
the variable
$
\hat{U} := X_{T_{\partial\mathcal{P}(0,\sqrt{j}L)} }
$
is distributed as
$
\mu_{0,0}^{\sqrt{j}L}.
$
Therefore, setting $U' := U - \hat{U}$,
in combination with (\ref{eq:manyProbEst}) and (\ref{eq:XZbd}) we
deduce the desired result.\vadjust{\goodbreak}~%
\end{pf}

The following lemma is essentially a discrete second-order Taylor expansion.
%
\begin{lemma} \label{lem:taylor}
Let $\mu$ be a finite signed measure on $\Z^d$ and let $f\dvtx\Z^d\to\R
$. Choose $m$, $k$, $J$, $N \in\N$
and $\varrho\in\Z^d$ such that:
\begin{longlist}[(a)]
\item[(a)] \hypertarget{item:TaylorFirstDerBd} for every $x,y \in\Z
^d$ such that $\Vert x-y \Vert_1 =1$,
we have $\vert f(x)-f(y) \vert<m$;
\item[(b)] \hypertarget{item:TaylorSecondDerBd} for every $x,y,z,w
\in\Z^d$ and $i,j \in\{1, \ldots, d\}$ such that $x-y=\break z-w=e_i$ and
$x-z=y-w=e_j$,
we have that $ \vert f(x)+f(w)-f(y)-f(z) \vert<k$ (note that if $i=j$
then $y=z$ and this is the discrete second derivative,
while if $i\neq j$ it is a discrete mixed second derivative);
\item[(c)] \hypertarget{item:TaylorMeasureBd} $\sum_x \mu(x)=0$;
\item[(d)] \hypertarget{item:TaylorExpBd} $\Vert{\sum_x x}\mu
(x)\Vert_1 \leq N$;
\item[(e)] \hypertarget{item:TaylorVarBd} $\sum_x \Vert x-\varrho
\Vert_1^2 \cdot\vert\mu(x) \vert< J$.
\end{longlist}
Then
\[
\biggl\vert\sum_{x} f(x) \mu(x) \biggr\vert\leq mN + kJ.
\]
\end{lemma}
\begin{pf}
From \hyperlink{item:TaylorMeasureBd}{(c)}, we infer that
$\sum_{x} f(x) \mu(x) = \sum_{x} (f(x)+c) \mu(x) $ for every $c \in
\R$.
Therefore, without loss of generality, we may assume that \mbox{$f(\varrho)=0$}.
Let $g\dvtx\Z^d\to\R$ be the affine function
characterized by
%
%
\begin{equation} \label{eq:gDef}\qquad\quad
g(\varrho)=f(\varrho)=0 \quad\mbox{and}\quad
g(\varrho+e_i)=f(\varrho+e_i)\qquad \forall i \in\{1,\ldots,d\}.
\end{equation}
Then
for any $x \in\Z^d$,
%
%
\begin{equation} \label{eq:finalBd}
\vert f(x)-g(x) \vert< k \Vert x-\varrho\Vert_1^2.
\end{equation}
In fact, setting $h:= f-g$ we get for $B(x,\varrho):= \{y \in\Z^d \dvtx
x_i \wedge\varrho_i \leq y_i \leq x_i \vee\varrho_i\allowbreak
\forall i \in\{1, \ldots, d\} \}$ that
%
%
\begin{equation} \label{eq:fnBd}
\vert f(x) - g(x) \vert\leq\vert h(\varrho) \vert+ \mathop{\max
_{i \in\{1, \ldots, d\},}}_{y \in B(x,\varrho)}
\biggl\vert\frac{\partial}{\partial e_i} h(y) \biggr\vert\cdot\Vert
x - \varrho
\Vert_1,
\end{equation}
where $\frac{\partial}{\partial e_i}h(y) := h(y+e_i) - h(y)$.

In addition, for $\frac{\partial^2}{\partial e_j \,\partial e_i} h(y)
:= h(y+e_i) - h(y) - (h(y+e_i +e_j) - h(y+e_j))$
we get for $y \in B(x,\varrho)$ that
%
%
\begin{equation} \label{eq:derBd}\quad
\biggl\vert\frac{\partial}{\partial e_i} h(y) \biggr\vert\leq
\biggl\vert\frac{\partial}{\partial e_i} h(\varrho) \biggr\vert
+ \mathop{\max_{j \in\{1, \ldots, d\},}}_{z \in B(x,\varrho)}
\biggl\vert\frac{ \partial^2}{\partial e_j\, \partial e_i} h(z) \biggr\vert
\cdot\Vert y - \varrho\Vert_1.
\end{equation}
Noting that $h(\varrho) = \frac{\partial}{\partial e_i} h(\varrho)
= 0$
as well as
$
\frac{\partial^2}{\partial e_j\, \partial e_i} h = \frac{\partial
^2}{\partial e_j \,\partial e_i} f,
$
and plugging~(\ref{eq:derBd}) into (\ref{eq:fnBd}), \hyperlink
{item:TaylorSecondDerBd}{(b)} yields (\ref{eq:finalBd}).

Now \hyperlink{item:TaylorVarBd}{(e)} in combination with (\ref
{eq:finalBd}) results in
\[
\biggl\vert\sum_{x} f(x) \mu(x) -\sum_{x} g(x) \mu(x) \biggr\vert
\leq\sum_x \vert f(x)-g(x) \vert\cdot\vert\mu(x)\vert\leq kJ.
\]
In addition, since $g$ is affine, $g - g(0)$ is linear and hence (\ref
{eq:gDef}) in combination with
\hyperlink{item:TaylorFirstDerBd}{(a)} and \hyperlink
{item:TaylorExpBd}{(d)} yields
\[
\biggl\vert\sum_{x} g(x) \mu(x) \biggr\vert= \biggl\vert g\biggl(\sum
_x x \mu(x) \biggr) - g(0) + \sum_{x} g(0) \mu(x) \biggr\vert
\leq mN.
\]
Due to the triangle inequality, these two estimates imply the statement
of the lemma.
\end{pf}
\begin{pf*}{Proof of Lemma \ref{lem:sumOfClose}} \label{pr:lem:sumOfClose}
We will construct a coupling that establishes the desired closeness.
For each $k \in\{1,\ldots,n\}$, conditioned on $\Delta_1,\ldots
,\Delta_{k-1}$, the distribution of $\Delta_k$ is
$(\lambda,K)$-close to $\mu_{0,0}^L$ by assumption, whence a coupling
as defined in Definition \ref{def:closeness} exists.
As mentioned in Remark \ref{rem:couplingSpace}, the coupling can be
constructed on the (possibly extended) probability
space the variables $\Delta_k$ are defined on, with $\Delta_k$
playing the role of $Z_2$ of that definition. We will assume such
couplings to be given.
Thus, for each such $k$ we still denote the variable corresponding to
$Z_2$ in Definition \ref{def:closeness}
by $\Delta_k$; the variable corresponding to $Z_0$ will be denoted by $Y_k$.
Without loss of generality, due to the fact that the $\Delta_k$'s and
$\mu_{0,0}^L$ are supported on $\partial_+ \mathcal{P}(0,L),$ we may
assume that the $Y_k$'s take values in $\partial_+ \mathcal{P}(0,L)$
only. Again,
without loss of generality, we assume all these couplings to be defined
on one
common probability space $(\Omega, \mathcal{F},P)$.
Thus, using the notation
$\mathcal{F}_{k-1} := \sigma(\Delta_1,\ldots,\Delta_{k-1})$ for $k
\in\{2, \ldots, n\}$
and $\mathcal{F}_0 := \{ \varnothing, \Omega\}$,
the following hold $P$-a.s.:
\begin{longlist}[(a$'$)]
\item[(a$'$)] \hypertarget{item:indAssProbDiffBound} $\sum_x \vert
P(Y_k=x \vert \mathcal{F}_{k-1})-\mu_{0,0}^L(x) \vert\leq
\lambda$;
\item[(b$'$)] \hypertarget{item:indAssDistBound} $P(\Vert Y_k-\Delta
_k \Vert_1 \leq K \vert \mathcal{F}_{k-1})=1$;
\item[(c$'$)] \hypertarget{item:indAssExpIdent} $E(Y_k \vert
\mathcal{F}_{k-1})=E_{\mu_{0,0}^L}$;
\item[(d$'$)] \hypertarget{item:indAssVarBd} $\sum_{x} \Vert
x-E_{\mu_{0,0}^L} \Vert_1^2 \cdot
\vert P(Y_k=x \vert \mathcal{F}_{k-1})-\mu_{0,0}^L(x) \vert
\leq\lambda\Var_{\mu_{0,0}^L}$.
\end{longlist}
To prove the desired result, it is sufficient to show that there exists
a random variable $Y'$ defined on the same
probability space such that:
\begin{longlist}[(a$''$)]
\item[(a$''$)] \hypertarget{item:indStepProbDiffBound} $\sum_x \vert
P(Y'=x)-\mu_{0,0}^{\sqrt{n}L}(x) \vert\leq\lambda R_{9}(L)$;
\item[(b$''$)] \hypertarget{item:indStepDistBound} $P (\Vert
Y'-S_n \Vert_1 <4nK )=1$;
\item[(c$''$)] \hypertarget{item:indStepExpIdent} $EY'=E_{\mu
_{0,0}^{\sqrt{n}L}}$;\vspace*{-2pt}
\item[(d$''$)] \hypertarget{item:indStepVarBd} $\sum_{x} \Vert
x-E_{\mu_{0,0}^{\sqrt{n}L}} \Vert_1^2 \cdot
\vert P(Y'=x)-\mu_{0,0}^{\sqrt{n}L}(x) \vert
\leq\lambda R_{9}(L) \Var_{\mu_{0,0}^{\sqrt{n}L}}$.
\end{longlist}
To this end, set
%
%
\begin{equation} \label{eq:sumSjDef}
S^{(j)} :=\sum_{k=j}^nY_k.
\end{equation}
Using descending induction, we start with showing that for all $j \in\{
1, \ldots, n\}$
the following holds:
\begin{longlist}[(IS)]
\item[(IS)] \hypertarget{item:IS}
Conditioned on $\Delta_1, \ldots, \Delta_{j-1}$, we can write
$S^{(j)}=Y^{(j)}+Z^{(j)}$ for some $Y^{(j)}$ and $Z^{(j)}$ such that
$\Vert Z^{(j)}\Vert_1 \leq(n-j)R_3(L)$
a.s. and such that with respect to
$P( \cdot \vert \mathcal{F}_{j-1})$, the variable $Y^{(j)}$
is distributed as $\mu_{0,0}^{\sqrt{n-j+1}L} +D_2^{(j)}$,
where $D_2^{(j)}$ is a signed measure
the variational norm
$\Vert D_2^{(j)} \Vert_{\mathrm{TV}}$ of which
is bounded from above by
$\lambda^{(j)}$
with $\lambda^{(n)}=\lambda$ and
\[
\lambda^{(j)} := \lambda^{(j+1)} + C\lambda R_6\bigl(\sqrt{n-j} L\bigr)(n-j)^{-1}
\]
for $j<n$
and some constant $C$.
\end{longlist}
For $j=n$, the statement holds true due to the assumptions with
$Z^{(n)}=0$. We now assume that the statement holds for $j+1$ and prove
it for $j$.

Setting
$H:=Y_j+Y^{(j+1)}$, for each $z$ we have
\[
P(H=z \vert \mathcal{F}_{j-1})
=\sum_{x} P(Y_j=x \vert \mathcal{F}_{j-1}) P\bigl(Y^{(j+1)}=z-x
\vert Y_j=x,\mathcal{F}_{j-1}\bigr).
\]
With $\hat{\mu}_{n,j,Y_j}^L$ defined as the convolution $P(Y_j \in
\cdot \vert \mathcal{F}_{j-1}) * \mu_{0,0}^{\sqrt{n-j}L}$,
this yields
that
%
%
\begin{eqnarray} \label{eq:inhar}\qquad
&&\sum_z \bigl\vert P(H=z \vert \mathcal{F}_{j-1})- \hat{\mu
}_{n,j,Y_j}^L (z) \bigr\vert\nonumber\\
&&\qquad\leq
\sum_z\sum_x P(Y_j=x \vert \mathcal{F}_{j-1}) \nonumber\\
&&\qquad\quad\hspace*{28.3pt}{} \times\bigl\vert P\bigl(Y^{(j+1)}=z-x \vert Y_j=x,
\mathcal{F}_{j-1}\bigr)-\mu_{0,0}^{\sqrt{n-j}L}(z-x) \bigr\vert\\
&&\qquad=
\sum_{x,y}P(Y_j=x \vert \mathcal{F}_{j-1}) \cdot\bigl\vert P
\bigl(Y^{(j+1)}=y \vert Y_j=x, \mathcal{F}_{j-1}\bigr)
-\mu_{0,0}^{\sqrt{n-j}L}(y) \bigr\vert\nonumber\\
&&\qquad\leq\bigl\Vert D_2^{(j+1)} \bigr\Vert_{\mathrm{TV}} \leq\lambda^{(j+1)}
\nonumber
\end{eqnarray}
holds a.s.

Next, we set $\hat{\mu}^{L}_{1,n-j} := \mu_{0,0}^L * \mu
_{0,0}^{\sqrt{n-j}L}$
and will bound
%
%
\begin{eqnarray} \label{eq:convMeasDiffs}
&&\vert\hat{\mu}^{L}_{1,n-j} (z)- \hat{\mu}^{L}_{n,j, Y_j}(z)
\vert\nonumber\\[-8pt]\\[-8pt]
&&\qquad= \biggl\vert\sum_x \mu_{0,0}^{\sqrt{n-j}L} (x) \bigl(P(Y_j=z-x
\vert \mathcal{F}_{j-1})-\mu_{0,0}^L(z-x) \bigr) \biggr\vert\nonumber
\end{eqnarray}
from above.

For this purpose, for given $z$,
we will apply Lemma \ref{lem:taylor} to the function~$\mu
_{0,0}^{\sqrt{n-j}L}$
with the corresponding measure $\mu$ given by $P(Y_j \in \cdot
\vert \mathcal{F}_{j-1} )-\mu_{0,0}^L$
(note that $\Vert\mu\Vert_{\mathrm{TV}} \leq\lambda$).

We now determine the parameters $k,m,J$ and $N$ of the assumptions of
Lemma \ref{lem:taylor}.
Parts \hyperlink{item:secondDer}{(d)} and \hyperlink{item:mixedSecondDer}{(e)}
of Lemma \ref{lem:annExitDistDerivatives} yield that we can choose
$k \leq C(\sqrt{n-j}L)^{-d-1}$.
Furthermore, as a consequence of \hyperlink
{item:indAssExpIdent}{(c$'$)} we can choose~$N$ equal to $0$, whence
the exact value of $m$ does not matter
($m=1$ works).
In addition, \hyperlink{item:indAssVarBd}{(d$'$)} yields that $J$ can
be chosen equal to $2\lambda\Var_{\mu_{0,0}^L}$,
with $\varrho$ equal to one of the elements of $\Z^d$ closest to
$E_{\mu_{0,0}^L}$.
Thus, Lemmas \ref{lem:VarBd} and \ref{lem:taylor} in combination
with (\ref{eq:convMeasDiffs}) yield that a.s.,
%
%
\begin{equation} \label{eq:compk}
\vert\hat{\mu}^{L}_{1,n-j} (z)- \hat{\mu}^{L}_{n,j, Y_j}(z)
\vert
\leq C \lambda L^{1-d}(n-j)^{({-d-1})/{2}}.
\end{equation}
Note that for $z$ such that
$\Vert z-E_{\hat{\mu}_{1,n-j}^{L}} \Vert_\infty> 4dR_6(\sqrt{n-j}
L)\sqrt{n-j}L$,
the terms $\hat{\mu}^{L}_{1,n-j}(z)$ and $\hat{\mu}^L_{n,j,Y_j}(z)$ vanish.
Thus, using (\ref{eq:inhar}) and (\ref{eq:compk}),
the triangle inequality implies that a.s.,
\begin{eqnarray*}
&&\sum_z \bigl\vert P(H=z  \vert \mathcal{F}_{j-1} ) - \hat{\mu
}_{1,n-j}^L (z) \bigr\vert\\
&&\qquad\leq\sum_z \bigl\vert P(H=z \vert \mathcal{F}_{j-1} ) - \hat
{\mu}_{n,j,Y_j}^L (z) \bigr\vert\\
&&\qquad\quad{} + \mathop{\sum_{z \in H_{(n-j)L^2},}}_{\Vert z - E_{\hat{\mu
}_{1,n-j}^{L}} \Vert_1 \leq4dR_6(\sqrt{n-j} L)\sqrt{n-j}L }
\vert\hat{\mu}_{n,j,Y_j}^L (z) - \hat{\mu}_{1,n-j}^L(z)
\vert\\
&&\qquad \leq\lambda^{(j+1)} + C\lambda R_6\bigl(\sqrt{n-j} L\bigr)^{d-1} (n-j)^{-1}.
\end{eqnarray*}
Consequently, we get that the distribution of $H$ can be written as
$\hat{\mu}_{1,n-j}^L + \overline{D}{}^{(j)}_2$ for a signed measure
$\overline{D}_2^{(j)}$ with
%
%
\begin{equation} \label{eq:auxProbBd}
\bigl\Vert\overline{D}_2^{(j)} \bigr\Vert_{\mathrm{TV}} \leq\lambda^{(j+1)} +
C\lambda R_6\bigl(\sqrt{n-j} L\bigr)^{d-1} (n-j)^{-1}.
\end{equation}
By Claim \ref{claim:annealedConvErr}, there exists $Z'(j)$ such that
\[
P\bigl( \Vert Z'(j) \Vert_1 >R_3(L)\bigr) \leq Ce^{-C^{-1}R_{2}(L)^\gamma},
\]
and
such that the distribution of $H+Z'(j)$ is
$\mu_{0,0}^{\sqrt{n-j+1}L}+\overline{D}_2^{(j)}$.
Let
\[
\overline{H} :=H+Z'(j)\cdot\mathbh{1}_{\Vert Z'(j)\Vert_1 \leq R_3(L)}.
\]
Then due to (\ref{eq:auxProbBd}), the distribution of $\overline{H}$
equals $\mu_{0,0}^{\sqrt{n-j+1}L} +\check{D}_2^{(j)}$
for some signed measure $\check{D}_2^{(j)}$ such that
\[
\bigl\Vert\check{D}_2^{(j)} \bigr\Vert_{\mathrm{TV}}\!\leq\!\bigl\Vert\overline{D}_2^{(j)}
\bigr\Vert_{\mathrm{TV}}\,{+}\,Ce^{-C^{-1}R_2(L)^\gamma}
\!\leq\!\lambda^{(j+1)}\,{+}\, C\lambda R_6\bigl(\sqrt{n\,{-}\,j} L\bigr)^{d-1} (n\,{-}\,j)^{-1}.
\]
We let
\[
Z^{(j)}:=Z^{(j+1)}+Z'(j)\cdot\mathbh{1}_{\Vert Z'(j) \Vert_1 \leq R_3(L)}
\]
and
\[
Y^{(j)}:=S^{(j)}-Z^{(j)}.
\]
Then we infer that
%
%
\begin{equation} \label{eq:ZjBd}
\bigl\Vert Z^{(j)} \bigr\Vert_1 \leq(n-j)R_{3}(L)
\end{equation}
and the distribution of $Y^{(j)}$ is
$\mu_{0,0}^{\sqrt{n-j+1}L}+D_2^{(j)}$ where $D_2^{(j)}$ is a signed
measure such that
$\Vert D_2^{(j)} \Vert_{\mathrm{TV}} \leq\lambda^{(j)}$ with
\[
\lambda^{(j)}\leq\lambda^{(j+1)} + C\lambda R_6\bigl(\sqrt{n-j}
L\bigr)^{d-1}(n-j)^{-1}.
\]
This establishes \hyperlink{item:IS}{(IS)}.

Using \hyperlink{item:indAssExpIdent}{(c$'$)} and (\ref
{eq:sumSjDef}), the expectation of $Y^{(1)}$ is computed via
%
%
\begin{equation}
EY^{(1)}=ES^{(1)}-EZ^{(1)}=nE_{\mu_{0,0}^L}-EZ^{(1)}.
\end{equation}
Therefore, in combination with (\ref{eq:ZjBd}), we get
\begin{eqnarray*}
\bigl\Vert EY^{(1)}-E_{\mu_{0,0}^{\sqrt{n}L}} \bigr\Vert_1
&\leq&\Vert nE_{\mu_{0,0}^L} -E_{\mu_{0,0}^{\sqrt{n}L}}
\Vert_1 + \bigl\Vert EZ^{(1)} \bigr\Vert_1\\
&\leq& CnR_2\bigl(\sqrt{n}L\bigr)+nR_3(L)\\
&\leq&2nR_{3}(L)
\end{eqnarray*}
for $L$ large enough, since $n \leq L$ by assumption;
indeed, with the help of Lem\-ma~\ref{lem:bddFluctExp} one deduces
\begin{eqnarray*}
\Vert n E_{\mu_{0,0}^L} - E_{\mu_{0,0}^{\sqrt{n} L}} \Vert_1
&\leq& n \biggl\Vert E_{\mu_{0,0}^L} - \frac{L^2}{\hat{v} \cdot e_1}
\hat{v} \biggr\Vert_1
+ \biggl\Vert\frac{nL^2}{\hat{v} \cdot e_1} \hat{v} - E_{\mu
_{0,0}^{\sqrt{n}L}} \biggr\Vert_1\\
&\leq& C \bigl( nR_2(L) + R_2\bigl(\sqrt{n}L\bigr) \bigr).
\end{eqnarray*}
As in the proof of Corollary \ref{cor:closeness}, we can find a
variable $U$ which is independent of all
the variables we have seen so far, and such that $\Vert U \Vert_1 \leq
2nR_{3}(L)$, $U \in H_0$ almost surely and
\[
EU=E_{\mu_{0,0}^{\sqrt{n}L}} -EY^{(1)}.
\]
We define
\[
Y':=Y^{(1)}+U,
\]
which directly yields that \hyperlink{item:indStepExpIdent}{(c$''$)} holds.
To check \hyperlink{item:indStepDistBound}{(b$''$)},
note that in combination with \hyperlink{item:IS}{(IS)} and the
definition of $S^{(1)}$ we get
%
%
\begin{eqnarray} \label{eq:SYPrimeTriangle}
\qquad\Vert S_n - Y' \Vert_1
&\leq&\bigl\Vert S_n - S^{(1)} \bigr\Vert_1
+ \bigl\Vert S^{(1)} - Y^{(1)} \bigr\Vert_1
+ \bigl\Vert Y^{(1)} - Y' \bigr\Vert_1\nonumber\\[-8pt]\\[-8pt]
\qquad&\leq& nK + nR_3 (L) + 2n R_{3} (L) \leq4nK
\nonumber
\end{eqnarray}
since $K \geq R_3(L)$.
Now from \hyperlink{item:IS}{(IS)} it follows that
$\lambda^{(1)} \leq C\lambda R_6(L^2)^{d-1} \log(n) \leq C\lambda
R_7(L)$ for $L$ large enough.
Thus,
\hyperlink{item:indStepProbDiffBound}{(a$''$)}
is a consequence of
%
%
\begin{eqnarray} \label{eq:YPrimeProbDiffBd}\quad
&&\sum_x \bigl\vert P(Y' = x) - \mu_{0,0}^{\sqrt{n}L} (x) \bigr\vert\nonumber\\[-2pt]
&&\qquad\leq2\sum_{x \in\partial_+ \mathcal{P}(0,\sqrt{n}L)}
\biggl\vert\mathop{\sum_{y \in H_0,}}_{
\Vert y \Vert_1 \leq2nR_{3}(L) } P(U=y) P\bigl(Y^{(1)} = x-y\bigr)
- \mu_{0,0}^{\sqrt{n}L}(x) \biggr\vert\nonumber\\[-2pt]
&&\qquad\leq2\sum_{x \in\partial_+ \mathcal{P}(0,\sqrt{n}L)}
\biggl( \biggl\vert\mathop{\sum_{y \in H_0,}}_{
\Vert y \Vert_1 \leq2nR_{3}(L) } P(U=y) \bigl( P\bigl(Y^{(1)} = x-y\bigr)
\nonumber\\[-10pt]\\[-10pt]
&&\qquad\quad\hspace*{175.8pt}{} - \mu_{0,0}^{\sqrt{n}L}(x-y) \bigr) \biggr\vert\nonumber\\[-2pt]
&&\qquad\quad\hspace*{157.1pt}{} + Cn R_{3}(L) \bigl(\sqrt{n}L\bigr)^{-d} \biggr)\nonumber\\[-2pt]
&&\qquad\leq C\lambda R_7(L)\nonumber
\end{eqnarray}
for $L$ large enough and where we used Lemma \ref
{lem:annExitDistDerivatives}\hyperlink{item:firstDer}{(b)} to obtain
the second inequality,
and also the fact that $\lambda\geq nL^{-1}$.

The remaining part of the proof consists of establishing that $Y'$ also
satisfies~\hyperlink{item:indStepVarBd}{(d$''$)}.
Denoting by $D_{2}$ be the signed measure such that
$Y' \sim\mu_{0,0}^{\sqrt{n}L} +D_{2}$, this amounts to showing that
%
%
\begin{equation} \label{eq:dBound}
\sum_{x} \Vert x-E_{\mu_{0,0}^{\sqrt{n}L}} \Vert_1^2
\cdot\vert D_2(x) \vert
\leq\lambda R_{9}(L) \Var_{\mu_{0,0}^{\sqrt{n}L}}
\end{equation}
holds.

To start with, note that
%
%
\begin{eqnarray} \label{eq:VarCoordDec}
&&\sum_x \Vert x-E_{\mu_{0,0}^{\sqrt{n}L}} \Vert_1^2 \cdot
\vert D_{2}(x) \vert\nonumber\\[-2pt]
&&\qquad\leq(d-1)\sum_x \Vert x-E_{\mu_{0,0}^{\sqrt{n}L}} \Vert
_2^2 \cdot\vert D_{2}(x) \vert\\[-2pt]
&&\qquad= (d-1)\sum_{i=2}^d \sum_x \bigl( (x-E_{\mu_{0,0}^{\sqrt
{n}L}} ) \cdot e_i \bigr)^2 \cdot\vert D_{2}(x) \vert.
\nonumber
\end{eqnarray}
To proceed, we write $D_{2}= D_2^+ - D_2^-$ for the Jordan
decomposition of $D_2$
and estimate
%
%
\begin{eqnarray}
&&\sum_x \bigl( ( x-E_{\mu_{0,0}^{\sqrt{n}L}} ) \cdot e_i
\bigr)^2 \cdot\vert D_{2}(x) \vert\nonumber\\[-2pt]
\label{eq:JordanFirstSum}
&&\qquad\leq2\sum_x \bigl( ( x-E_{\mu_{0,0}^{\sqrt{n}L}} ) \cdot
e_i \bigr)^2 \cdot D_{2}^-(x) \\[-2pt]
\label{eq:JordanSecondSum}
&&\qquad\quad{} +\biggl\vert\sum_x \bigl( (x-E_{\mu_{0,0}^{\sqrt{n}L}}
) \cdot e_i \bigr)^2 \cdot D_{2} (x) \biggr\vert.
\end{eqnarray}
To bound (\ref{eq:JordanFirstSum}) from above, note that $D_{2}^-(x)
\leq\mu_{0,0}^{\sqrt{n}L} (x)$ for all $x$.
Combined with the fact that $\Vert D_{2}^-\Vert_{\mathrm{TV}} \leq\Vert D_{2}
\Vert_{\mathrm{TV}} \leq C\lambda R_7(L)$ [due
to (\ref{eq:YPrimeProbDiffBd})],
we obtain
%
%
\begin{equation}\label{eq:boundforminus}
\sum_x \bigl( (x-E_{\mu_{0,0}^{\sqrt{n}L}}) \cdot e_i \bigr)^2
D_{2}^- (x) \leq\lambda R_8(L) n L^2
\end{equation}
for $L$ large enough, since $\mu_{0,0}^{\sqrt{n}L}$ is supported on
$\partial_+ \mathcal{P}(0,\sqrt{n}L)$.

In order to estimate (\ref{eq:JordanSecondSum}), note that
due to \hyperlink{item:indStepExpIdent}{(c$''$)} we have $\sum_x x
D_2(x) =0$, and hence (\ref{eq:JordanSecondSum}) equals
$\vert{\Var}(D_{2},i) \vert$ with
%
%
\begin{equation} \label{eq:VarD2Def}
\Var(D_{2},i) := \sum_x (x \cdot e_i)^2 D_{2}(x) = \Var(Y' \cdot
e_i) - \Var(W \cdot e_i),
\end{equation}
where $W$ denotes a random variable distributed
according to $\mu_{0,0}^{\sqrt{n}L}\hspace*{-0.8pt}$.
By Claim~\ref{claim:annealedConvErr}, there exists a random variable
$W'$ such that $W' \sim(\mu_{0,0}^L)^{*n}$, with $(\mu
_{0,0}^L)^{*n}$ denoting the $n$-fold convolution of $\mu_{0,0}^L$, and
such that
\[
P\bigl(\Vert W-W' \Vert_1 >nR_3(L)\bigr) < Cne^{-C^{-1}R_2(L)^\gamma}
\]
for $L$ large enough.
Then
%
%
\begin{eqnarray} \label{eq:VarDiffTriangle}
\vert{\Var}(D_{2},i) \vert
&\leq&\vert{\Var}(W \cdot e_i) - \Var(W' \cdot e_i) \vert\nonumber\\
&&{} + \bigl\vert{\Var}(W' \cdot e_i) - \Var\bigl(S^{(1)} \cdot e_i\bigr) \bigr\vert\\
&&{} +\bigl\vert{\Var}\bigl(S^{(1)} \cdot e_i\bigr) - \Var(Y' \cdot e_i) \bigr\vert.
\nonumber
\end{eqnarray}
Now for $L$ large enough,
%
%
\begin{eqnarray} \label{eq:wprimew}
&&\vert{\Var}(W \cdot e_i) - \Var( W' \cdot e_i) \vert\nonumber\\
&&\qquad\leq\Var\bigl((W-W') \cdot e_i\bigr) \nonumber\\
&&\qquad\quad{} + 2 \bigl\vert\Cov\bigl((W-W')\cdot e_i, W' \cdot e_1\bigr)
\bigr\vert\nonumber\\[-8pt]\\[-8pt]
&&\qquad\leq\esssup\bigl(( W - W') \cdot e_i \bigr)^2 P\bigl(\Vert W-W' \Vert_1
>nR_3(L)\bigr) \nonumber\\
&&\qquad\quad{} + 2n^2R_3(L)^2+2nR_3(L)\sqrt{{\Var(W')}} \nonumber\\
&&\qquad\leq Cn^{3/2}R_3(L)L,\nonumber
\end{eqnarray}
where among others we used $\sqrt{n} \leq L$ and Lemma \ref{lem:VarBd}.
Furthermore,
%
%
\begin{eqnarray}
\nonumber
\bigl\vert{\Var}(Y' \cdot e_i) - \Var\bigl(S^{(1)} \cdot e_i \bigr) \bigr\vert
&=& \bigl\vert{\Var}\bigl(\bigl(S^{(1)}+U'\bigr)\cdot e_i\bigr) - \Var\bigl(S^{(1)} \cdot e_i\bigr)
\bigr\vert\nonumber\\
&\leq&2\esssup(\Vert U'\Vert_1) \sqrt{\Var\bigl(S^{(1)}\bigr)} + \esssup
(\Vert U'\Vert_1)^2\nonumber\\
&\leq& Cn^{3/2}R_{3}(L)L+4 n^2R_{3}^2(L)\\
\label{eq:yprimey}
&\leq& Cn^{3/2}R_{3}(L)L,
\end{eqnarray}
where we used $n \leq L$ and
that by the definition of $Y'$, we know that $U':=Y'-S^{(1)}$ satisfies
$\Vert U' \Vert_1 \leq3nR_{3}(L)$.

To estimate the remaining summand, note that $\Cov(Y_j,Y_k)=0$ for
$j\not= k$, and hence
%
%
\begin{equation} \label{eq:SWPrimeEst}\qquad\qquad
\bigl\vert{\Var}\bigl(S^{(1)} \cdot e_i\bigr) - \Var( W' \cdot e_i) \bigr\vert
\leq\sum_{j=1}^n \bigl\vert{\Var}(Y_j \cdot e_i)-\Var_{\mu_{0,0}^L
\circ(\cdot\cdot e_i)^{-1}} \bigr\vert.
\end{equation}
Furthermore, since $E_{\mu_{0,0}^L} = EY_j$ for every $j$, from
\hyperlink{item:indAssVarBd}{(d$'$)} we infer that
\begin{eqnarray*}
&&\bigl\vert{\Var}( Y_j \cdot e_i )-\Var_{\mu_{0,0}^L \circ(\cdot
\cdot{e_i})^{-1}} \bigr\vert\\
&&\qquad= \biggl\vert\sum_x \bigl( x \cdot e_i - (E_{\mu_{0,0}^L} \cdot e_i)
\bigr)^2 \bigl(P(Y_j=x)-\mu_{0,0}^L(x)\bigr) \biggr\vert\\
&&\qquad\leq\sum_x \bigl( x \cdot e_i -(E_{\mu_{0,0}^L } \cdot e_i)
\bigr)^2 \vert P(Y_j=x)-\mu_{0,0}^L(x) \vert\\
&&\qquad\leq\lambda\Var_{\mu_{0,0}^L},
\end{eqnarray*}
and as a consequence
%
%
\begin{equation} \label{eq:VarS1WPrimeBd}
\bigl\vert{\Var}\bigl(S^{(1)} \cdot e_i\bigr) - \Var( W' \cdot e_i) \bigr\vert\leq
n\lambda\Var_{\mu_{0,0}^L}.
\end{equation}
Using (\ref{eq:VarDiffTriangle}) to (\ref{eq:VarS1WPrimeBd}) in
combination with
Lemma \ref{lem:VarBd}, we deduce that
\[
\vert{\Var}(D_2,i) \vert\leq C\lambda n L^2, \label{eq:vard2}
\]
whence in combination with (\ref{eq:JordanFirstSum}) to (\ref
{eq:VarD2Def}) we have
\[
\sum_x ( x \cdot e_i -E_{\mu_{0,0}^{\sqrt{n}L}} \cdot e_i
)^2 \vert D_{2} (x) \vert
\leq2R_8(L)\lambda nL^2
\]
for $L$ large enough.
Therefore, (\ref{eq:VarCoordDec}) yields
\[
\sum_x \Vert x-E_{\mu_{0,0}^{\sqrt{n}L}} \Vert_1^2 \cdot
\vert D_{2}(x) \vert
\leq2d^2 R_8(L) \lambda nL^2
\]
for $L$ large enough. In combination with Lemma \ref{lem:VarBd}, we
deduce that
(\ref{eq:dBound}) holds and thus \hyperlink
{item:indStepVarBd}{(d$''$)} is fulfilled.
\end{pf*}

We now prove the previously employed Lemma \ref{lem:lbound}.
\begin{pf*}{Proof of Lemma \ref{lem:lbound}} \label{pr:lem:lbound}
We continue to use the
notation $B(l,k)$\break and~$B(l)$ introduced in the proof of Lemma \ref
{lem:annExitDistDerivatives},
from which this proof draws its strategy.\vadjust{\goodbreak}
Again, denote by $\Gamma$ the covariance matrix of $X_{\tau_2} -
X_{\tau_1}$ with respect
to $P_0$ and set $m:= E_0 (X_{\tau_2} - X_{\tau_1})$.
Using (\ref{eq:LCLT})
and the fact that, since $\Gamma^{-1}$ is positive definite, the corresponding
quadratic form induces a~norm,
we infer that for any $C > 0$ there exists a constant $c>0$ such that
for~$k$ large enough and $y \in H_{L^2}$
with $\Vert y - x - km \Vert_1 \leq C \sqrt{k}$, we have
%
%
\begin{equation} \label{eq:LCLTLowerEst}
P_x \bigl( X_{T_{L^2}} = y, B(L^2,k) \bigr) \geq c k^{-{d}/{2}}.
\end{equation}
Setting $l^* :=\frac{{L^2} - x \cdot e_1}{m \cdot e_1}$
and $C := 4C'$,
for $k \in\{{\lfloor l^* - \sqrt{l^*} \rfloor}, \ldots, \lceil l^*
+ \sqrt
{l^*} \rceil\}$
and~$x,y$ as in the assumptions,
we have
\[
\Vert y- x - k m \Vert_1 \leq\Vert\tilde{\pi}_{\hat{v}^\bot}
(y-x) \Vert_1
+ \Vert l^* m - k m\Vert_1 \leq CL
\]
for $L$ large enough.
Then,
using
(\ref{eq:LCLTLowerEst}) with $y \in\partial_+ \mathcal{P}(0,L)$,
uniformly in $x \in\tilde{\mathcal{P}}(0,L)$
we have
\begin{eqnarray*}
P_x (X_{T_{L^2}}=y )
&\geq& P_x \bigl(X_{T_{L^2}} =y, B({L^2})\bigr)\\
&\geq&\sum_{k={\lfloor{l^*}-\sqrt{{l^*}} \rfloor}}^{\lceil
{l^*}+\sqrt
{{l^*}}\rceil}
P_x \bigl(X_{T_{L^2}} =y, B({L^2},k) \bigr)\\
&\geq& cL^{1-d},
\end{eqnarray*}
which due to (\ref{eq:supAnnealedNotFrontExitEst}) finishes the proof.
\end{pf*}
\end{appendix}

\section*{Acknowledgment} We thank A.-S. Sznitman for introducing this
topic to us and for useful discussions.

The above forms part of the first author's Ph.D. thesis.

The final presentation of this article has profited from beneficial
hints of the referee.

%

%
\printaddresses

\end{document}